\documentclass[11pt]{amsart}
\usepackage{stmaryrd}		
\usepackage{amscd}		
   \title{Hilbert schemes, polygraphs, and the Macdonald positivity
	  conjecture}
                                                                   
   \author{Mark Haiman}                                            

   \address{Dept.\ of Mathematics\\
            U.C. San Diego\\
            La Jolla, CA, 92093-0112}

   \email{mhaiman@macaulay.ucsd.edu}
                                                                   
   \date{December 5, 2000}                                  



\thanks{Research supported in part by N.S.F. Mathematical Sciences
grant DMS-9701218.}

\DeclareMathOperator{\Hilb}{Hilb}
\DeclareMathOperator{\Spec}{Spec}
\DeclareMathOperator{\Proj}{Proj}
\DeclareMathOperator{\Hom}{Hom}
\DeclareMathOperator{\soc}{soc}
\DeclareMathOperator{\ch}{ch}
\DeclareMathOperator{\rk}{rk}
\DeclareMathOperator{\im}{im}
\DeclareMathOperator{\coker}{coker}
\DeclareMathOperator{\depth}{depth}
\DeclareMathOperator{\SL}{SL}

\newcommand{\Ocal}{{\mathcal O}}

\newcommand{\Jcal}{{\mathcal J}}
\newcommand{\Fcal}{{\mathcal F}}
\newcommand{\Hcal}{{\mathcal H}}
\newcommand{\Dcal}{{\mathcal D}}
\newcommand{\Lcal}{{\mathcal L}}

\renewcommand{\AA}{{\mathbb A}}
\newcommand{\PP}{{\mathbb P}}
\newcommand{\ZZ}{{\mathbb Z}}
\newcommand{\QQ}{{\mathbb Q}}
\newcommand{\NN}{{\mathbb N}}
\newcommand{\CC}{{\mathbb C}}
\newcommand{\FF}{{\mathbb F}}
\newcommand{\TT}{{\mathbb T}}

\renewcommand{\aa}{{\mathbf a}}
\newcommand{\bb}{{\mathbf b}}
\newcommand{\xx}{{\mathbf x}}
\newcommand{\yy}{{\mathbf y}}
\newcommand{\zz}{{\mathbf z}}

\newcommand{\rad}[1]{\sqrt{\vphantom (}#1} 

\newcommand{\icolon}{\,{:}\,}

\newcommand{\ext}{\wedge}

\newcommand{\Ghilb}[2]{#2{\sslash }#1}

\theoremstyle{plain}
\newtheorem{thm}{Theorem}

\newtheorem{prop}{Proposition}[subsection]
\newtheorem{cor}[prop]{Corollary}
\newtheorem{lem}[prop]{Lemma}
\newtheorem{conj}[prop]{Conjecture}
\newtheorem{prob}[prop]{Problem}

\theoremstyle{definition}
\newtheorem{defn}[prop]{Definition}
\newtheorem{conv}[prop]{Convention}
\begin{document}

\begin{abstract}
We study the {\it isospectral Hilbert scheme} $X_{n}$, defined as the
reduced fiber product of $(\CC ^{2})^{n}$ with the Hilbert scheme
$H_{n}$ of points in the plane $\CC ^{2}$, over the symmetric power
$S^{n}\CC ^{2} = (\CC ^{2})^{n}/S_{n}$.  By a theorem of Fogarty,
$H_{n}$ is smooth.  We prove that $X_{n}$ is normal, Cohen-Macaulay,
and Gorenstein, and hence flat over $H_{n}$.  We derive two important
consequences.

(1) We prove the strong form of the $n!$ {\it conjecture} of Garsia
and the author, giving a representation-theoretic interpretation of
the Kostka-Macdonald coefficients $K_{\lambda \mu }(q,t)$.  This
establishes the {\it Macdonald positivity conjecture}, namely that
$K_{\lambda \mu }(q,t)\in \NN [q,t]$.

(2) We show that the Hilbert scheme $H_{n}$ is isomorphic to the {\it
$G$-Hilbert scheme} $\Ghilb{S_{n}}{(\CC ^{2})^{n}}$ of Nakamura, in
such a way that $X_{n}$ is identified with the universal family over
$\Ghilb{S_{n}}{(\CC ^{2})^{n}}$.  From this point of view, $K_{\lambda
\mu }(q,t)$ describes the fiber of a {\it character sheaf} $C_{\lambda
}$ at a torus-fixed point of $\Ghilb{S_{n}}{(\CC ^{2})^{n}}$
corresponding to $\mu $.

The proofs rely on a study of certain subspace arrangements
$Z(n,l)\subseteq (\CC ^{2})^{n+l}$, called {\it polygraphs}, whose
coordinate rings $R(n,l)$ carry geometric information about $X_{n}$.
The key result is that $R(n,l)$ is a free module over the polynomial
ring in one set of coordinates on $(\CC ^{2})^{n}$.  This is proven by
an intricate inductive argument based on elementary commutative
algebra.
\end{abstract}

\maketitle

\tableofcontents


\section{Introduction}
\label{s:intro}

The Hilbert scheme of points in the plane $H_{n} = \Hilb ^{n}(\CC
^{2})$ is an algebraic variety which parametrizes finite subschemes
$S$ of length $n$ in $\CC ^{2}$.  To each such subscheme $S$
corresponds an $n$-element multiset, or unordered $n$-tuple with
possible repetitions, $\sigma (S) = \llbracket P_{1},\ldots,P_{n}
\rrbracket$ of points in $\CC ^{2}$, where the $P_{i}$ are the points
of $S$, repeated with appropriate multiplicities.  There is a variety
$X_{n}$, finite over $H_{n}$, whose fiber over the point of $H_{n}$
corresponding to $S$ consists of all ordered $n$-tuples
$(P_{1},\ldots,P_{n})\in (\CC ^{2})^{n}$ whose underlying multiset is
$\sigma (S)$.  We call $X_{n}$ the {\it isospectral Hilbert scheme}.

By a theorem of Fogarty \cite{Fog68}, the Hilbert scheme $H_{n}$ is
irreducible and nonsingular.  The geometry of $X_{n}$ is more
complicated, but also very special.  Our main geometric result,
Theorem~\ref{t:main} in \S \ref{ss:hilbert/main}, is that $X_{n}$ is
normal, Cohen-Macaulay and Gorenstein.

Earlier investigations by the author \cite{Hai99} unearthed
indications of a far-reaching correspondence between the geometry and
sheaf cohomology of $H_{n}$ and $X_{n}$ on the one hand, and the
theory of {\it Macdonald polynomials} on the other.  The Macdonald
polynomials
\begin{equation}\label{e:intro/P-mu}
P_{\mu }(x;q,t)
\end{equation}
are a basis of the algebra of symmetric functions in variables $x =
x_{1},x_{2},\ldots \,$, with coefficients in the field $\QQ (q,t)$ of
rational functions in two parameters $q$ and $t$.  They were
introduced in 1988 by Macdonald \cite{Mac88} to unify the two
well-known one-parameter bases of the algebra of symmetric functions,
namely, the {\it Hall-Littlewood polynomials} and the {\it Jack
polynomials} (for a thorough treatment see \cite{Mac95}).  It promptly
became clear that the discovery of Macdonald polynomials was
fundamental and sure to have many ramifications.  Developments in the
years since have borne this out, notably, Cherednik's proof of the
Macdonald constant-term identities \cite{Cher96} and other discoveries
relating Macdonald polynomials to representation theory of quantum
groups \cite{EtKi94} and affine Hecke algebras \cite{KiNo98, Kno97,
Mac96}, the Calogero--Sutherland model in particle physics
\cite{LaVi98}, and combinatorial conjectures on diagonal harmonics
\cite{BeGaHaTe99, GaHa96, Hai94}.


The link between Macdonald polynomials and Hilbert schemes comes from
work by Garsia and the author on the {\it Macdonald positivity
conjecture}.  The Schur function expansions of Macdonald polynomials
lead to transition coefficients $K_{\lambda \mu }(q,t)$, known as {\it
Kostka-Macdonald coefficients}.  As defined, they are rational
functions of $q$ and $t$, but conjecturally they are polynomials in
$q$ and $t$ with nonnegative integer coefficients:
\begin{equation}\label{e:intro/MPK}
K_{\lambda \mu }(q,t)\in \NN [q,t].
\end{equation}
The positivity conjecture has remained open since Macdonald formulated
it at the time of his original discovery.  For $q=0$ it reduces to the
positivity theorem for $t$-Kostka coefficients, which has important
algebraic, geometric and combinatorial interpretations \cite{BrHa98,
DCPr81, GaPr92, HoSp77, Kat82, Kra81, LaSc78, Lus81, Lus83, Spr78}.
Only recently have several authors independently shown that the
Kostka-Macdonald coefficients are polynomials, $K_{\lambda \mu
}(q,t)\in \ZZ [q,t]$, but these results do not establish the
positivity \cite{GaRe96, GaTe96, KiNo98, Kno97, Sah96}.

In \cite{GaHa93}, Garsia and the author conjectured an interpretation
of the Kostka-Macdonald coefficients $K_{\lambda \mu }(q,t)$ as graded
character multiplicities for certain doubly graded $S_{n}$-modules
$D_{\mu }$.  The module $D_{\mu }$ is the space of polynomials in $2n$
variables spanned by all derivatives of a certain simple determinant
(see \S \ref{ss:n!/n!} for the precise definition).  The conjectured
interpretation implies the Macdonald positivity conjecture.  It also
implies, in consequence of known properties of the $K_{\lambda \mu
}(q,t)$, that for each partition $\mu $ of $n$, the dimension of
$D_{\mu }$ is equal to $n!$.  This seemingly elementary assertion has
come to be known as the {\it $n!$ conjecture}.

It develops that these conjectures are closely tied to the geometry of
the isospectral Hilbert scheme.  Specifically, in \cite{Hai99} we were
able show that the Cohen-Macaulay property of $X_{n}$ is equivalent to
the $n!$ conjecture.  We further showed that the Cohen-Macaulay
property of $X_{n}$ implies the stronger conjecture interpreting
$K_{\lambda \mu }(q,t)$ as a graded character multiplicity for $D_{\mu
}$.  Thus the geometric results in the present article complete the
proof of the Macdonald positivity conjecture.


Another consequence of our results, equivalent in fact to our main
theorem, is that the Hilbert scheme $H_{n}$ is equal to the {\it
$G$-Hilbert scheme} $\Ghilb{G}{V}$ of Ito and Nakamura \cite{ItNa96},
for the case $V = (\CC ^{2})^{n}$, $G = S_{n}$.  The $G$-Hilbert
scheme is of interest in connection with the generalized {\it McKay
correspondence}, which says that if $V$ is a complex vector space, $G$
is a finite subgroup of $\SL (V)$ and $Y\rightarrow V/G$ is a
so-called {\it crepant} resolution of singularities, then the sum of
the Betti numbers of $Y$ equals the number of conjugacy classes of
$G$.  In many interesting cases \cite{BrKiRe99, Nak00}, the
$G$-Hilbert scheme turns out to be a crepant resolution and an
instance of the McKay correspondence.  By our main theorem, this holds
for $G = S_{n}$, $V = (\CC ^{2})^{n}$.

We wish to say a little at this point about how the discoveries
presented here came about.  It has long been known \cite{HoSp77,
Spr78} that the $t$-Kostka coefficients $K_{\lambda \mu }(t) =
K_{\lambda \mu }(0,t)$ are graded character multiplicities for the
cohomology rings of Springer fibers.  Garsia and Procesi \cite{GaPr92}
found a new proof of this result, deriving it directly from an
elementary description of the rings in question.  In doing so, they
hoped to reformulate the result for $K_{\lambda \mu }(t)$ in a way
that might generalize to the two-parameter case.  Shortly after that,
Garsia and the author began their collaboration and soon found the
desired generalization, in the form of the $n!$ conjecture.  Based on
Garsia and Procesi's experience, we initially expected that the $n!$
conjecture itself would be easy to prove and that the difficulties
would lie in the identification of $K_{\lambda \mu }(q,t)$ as the
graded character multiplicity.  To our surprise, however, the $n!$
conjecture stubbornly resisted elementary attack.

In the spring of 1992, we discussed our efforts on the $n!$ conjecture
with Procesi, along with another related conjecture we had stumbled
upon in the course of our work.  The modules involved in the $n!$
conjecture are quotients of the ring $R_{n}$ of coinvariants for the
action of $S_{n}$ on the polynomial ring in $2n$ variables.  This ring
$R_{n}$ is isomorphic to the space of diagonal harmonics.
Computations suggested that its dimension should be $(n+1)^{n-1}$ and
that its graded character should be related to certain well-known
combinatorial enumerations (this conjecture is discussed briefly in \S
\ref{ss:other/bases} and at length in \cite{GaHa96, Hai94}).  Procesi
suggested that the Hilbert scheme $H_{n}$ and what we now call the
isospectral Hilbert scheme $X_{n}$ should be relevant to the
determination of the dimension and character of $R_{n}$.
Specifically, he observed that there is a natural map from $R_{n}$ to
the ring of global functions on the scheme-theoretic fiber in $X_{n}$
over the origin in $S^{n}\CC ^{2}$.  With luck, this map might be an
isomorphism, and---as we are now able to confirm---$X_{n}$ might be
flat over $H_{n}$, so that its structure sheaf would push down to a
vector bundle on $H_{n}$.  Then $R_{n}$ would coincide with the global
sections of the vector bundle over the zero-fiber in $H_{n}$, and it
might be possible to compute its character using the Atiyah--Bott
Lefschetz formula.

The connection between $X_{n}$ and the $n!$ conjecture became clear
when the author sought to carry out the computation Procesi had
suggested, assuming the validity of the needed but unproven geometric
hypotheses.  More precisely, it became clear that the spaces in the
$n!$ conjecture should be the fibers of Procesi's vector bundle at
distinguished torus-fixed points in $H_{n}$ (a fact we prove in \S
\ref{ss:hilbert/n!}).  This led to a conjectured formula for the
character of $R_{n}$ in terms of Macdonald polynomials, which turned
out to be correct up to the limit of practical computation ($n\leq
7$).  In \cite{GaHa96} Garsia and the author further showed that the
series of combinatorial conjectures in \cite{Hai94} would follow from
the conjectured master formula.  Thus we had strong indications that
Procesi's proposed picture was indeed valid, and that a geometric
study of $X_{n}$ should ultimately lead to a proof of the $n!$ and
Macdonald positivity conjectures, as is borne out here.

By now the reader should expect the geometric study of $X_{n}$ also to
yield a proof of the character formula for diagonal harmonics and the
$(n+1)^{n-1}$ conjecture.  This subject will be taken up in a separate
article.


The remainder of the paper is organized as follows.  In Section
\ref{s:n!} we give the relevant definitions concerning Macdonald
polynomials and state the positivity, $n!$ and graded character
conjectures.  Hilbert scheme definitions and the statement and proof
of the main theorem, Theorem~\ref{t:main}, are in Section
\ref{s:hilbert}, along with the equivalence of the main theorem to the
$n!$ conjecture.  In \S \ref{ss:hilbert/character} we review the proof
from \cite{Hai99} that the main theorem implies the conjecture of
Garsia and the author on the character of the space $D_{\mu }$, and
hence implies the Macdonald positivity conjecture.

The proof of the main theorem uses a technical result,
Theorem~\ref{t:polygraphs-are-free}, that the coordinate ring of a
certain type of subspace arrangement we call a {\it polygraph} is a
free module over the polynomial ring generated by some of the
coordinates.  Section \ref{s:polygraphs} contains the definition and
study of polygraphs, culminating in the proof of
Theorem~\ref{t:polygraphs-are-free}.  At the end, in Section
\ref{s:other}, we discuss other implications of our results, including
the connection with $G$-Hilbert schemes, along with related
conjectures and open problems.


\section{The $n!$ and Macdonald positivity conjectures}
\label{s:n!}

\subsection{Macdonald polynomials}
\label{ss:n!/macdonald}

We work with the {\it transformed integral forms} $\tilde{H}_{\mu
}(x;q,t)$ of the Macdonald polynomials, indexed by integer partitions
$\mu $, and homogeneous of degree $n = |\mu |$.  These are defined as
in \cite{Hai99}, eq.~(2.18) to be
\begin{equation}\label{e:n!/Htilde}
\tilde{H}_{\mu }(x;q,t) = t^{n(\mu )}J_{\mu }[X/(1-t^{-1});q,t^{-1}],
\end{equation}
where $J_{\mu }$ denotes Macdonald's integral form as in \cite{Mac95},
Chapter VI, eq.~(8.3), and $n(\mu )$ is the partition statistic
\begin{equation}\label{e:n!/n(mu)}
n(\mu ) = \sum _{i}(i-1)\mu _{i}
\end{equation}
(not to be confused with $n = |\mu |$).

The square brackets in \eqref{e:n!/Htilde} stand for {\it plethystic
substitution}.  We pause briefly to review the definition of this
operation (see \cite{Hai99} for a fuller discussion).  Let $\FF [[x]]$
be the algebra of formal series over the coefficient field $\FF = \QQ
(q,t)$, in variables $x = x_{1},x_{2},\ldots\;$.  For any $A\in \FF
[[x]]$, we denote by $p_{k}[A]$ the result of replacing each
indeterminate in $A$ by its $k$-th power.  This includes the
indeterminates $q$ and $t$ as well as the variables $x_{i}$.  The
algebra of symmetric functions $\Lambda _{\FF }$ is freely generated
as an $\FF $-algebra by the power-sums
\begin{equation}\label{e:n!/power-sum}
p_{k}(x) = x_{1}^{k}+x_{2}^{k}+\cdots .
\end{equation}
Hence there is a unique $\FF $-algebra homomorphism
\begin{equation}\label{e:n!/ev-A}
\operatorname{ev}_{A}\colon \Lambda _{\FF }\rightarrow \FF [[x]]\quad
\text{defined by}\quad p_{k}(x)\mapsto  p_{k}[A].
\end{equation}
In general we write $f[A]$ for $\operatorname{ev}_{A}(f)$, for any
$f\in \Lambda _{\FF }$.  With this notation goes the convention that
$X$ stands for the sum $X = x_{1}+x_{2}+\cdots $ of the variables, so
we have $p_{k}[X] = p_{k}(x)$ and hence $f[X] = f(x)$ for all $f$.
Note that a plethystic substitution like $f\mapsto f[X/(1-t^{-1})]$,
such as we have on the right-hand side in \eqref{e:n!/Htilde},
yields again a symmetric function.


There is a simple direct characterization of the transformed Macdonald
polynomials $\tilde{H}_{\mu }$.

\begin{prop}[\cite{Hai99}, Proposition 2.6]\label{p:n!/H-tilde-triang}
The $\tilde{H}_{\mu }(x;q,t)$ satisfy
\begin{itemize}
\item [(1)] $\tilde{H}_{\mu }(x;q,t)\in \QQ (q,t)\{
s_{\lambda}[X/(1-q)] : \lambda \geq \mu \}$,
\item [(2)] $\tilde{H}_{\mu }(x;q,t)\in \QQ (q,t)\{
s_{\lambda}[X/(1-t)] : \lambda \geq \mu' \}$, and
\item [(3)] $\tilde{H}_{\mu }[1;q,t] = 1$,
\end{itemize}
where $s_{\lambda }(x)$ denotes a Schur function, $\mu '$ is the
partition conjugate to $\mu $, and the ordering is the dominance
partial order on partitions of $n = |\mu |$.  These conditions
characterize $\tilde{H}_{\mu }(x;q,t)$ uniquely.
\end{prop}

We set $\tilde{K}_{\lambda \mu }(q,t) = t^{n(\mu )}K_{\lambda \mu
}(q,t^{-1})$, where $K_{\lambda \mu }(q,t)$ is the Kostka-Macdonald
coefficient defined in \cite{Mac95}, Chapter VI, eq.~(8.11).  This is
then related to the transformed Macdonald polynomials by
\begin{equation}\label{e:n!/K-tilde}
\tilde{H}_{\mu }(x;q,t) = \sum _{\lambda } \tilde{K}_{\lambda \mu
}(q,t) s_{\lambda }(x).
\end{equation}
It is known that $K_{\lambda \mu }(q,t)$ has degree at most $n(\mu )$
in $t$, so the positivity conjecture \eqref{e:intro/MPK} from the
introduction can be equivalently formulated in terms of
$\tilde{K}_{\lambda \mu }$.

\begin{conj}[Macdonald positivity conjecture]\label{conj:n!/MPK}
We have $\tilde{K}_{\lambda \mu }(q,t)\in \NN [q,t]$.
\end{conj}


\subsection{The $n!$ and graded character conjectures}
\label{ss:n!/n!}

Let $\CC [\xx ,\yy ] = \CC [x_{1},y_{1},\ldots,x_{n},y_{n}]$ be the
polynomial ring in $2n$ variables.  To each $n$-element subset
$D\subseteq \NN \times \NN $, we associate a polynomial $\Delta
_{D}\in \CC [\xx ,\yy ]$ as follows.  Let
$(p_{1},q_{1}),\ldots,(p_{n},q_{n})$ be the elements of $D$ listed in
some fixed order.  Then we define
\begin{equation}\label{e:n!/Delta-D}
\Delta _{D} = \det \left(x_{i}^{p_{j}} y_{i}^{q_{j}} \right)_{1\leq
i,j\leq n}.
\end{equation}
If $\mu $ is a partition of $n$, its {\it diagram} is the set
\begin{equation}\label{e:n!/diagram-of-mu}
D(\mu ) = \{(i,j)\colon j<\mu _{i+1} \}\subseteq \NN \times \NN .
\end{equation}
(Note that in our definition the rows and columns of the diagram
$D(\mu )$ are indexed starting with zero.)  In the case where $D =
D(\mu )$ is the diagram of a partition, we abbreviate
\begin{equation}\label{e:n!/Delta-mu}
\Delta _{\mu } = \Delta _{D(\mu )}.
\end{equation}
The polynomial $\Delta _{\mu }(\xx ,\yy )$ is a kind of bivariate
analog of the Vandermonde determinant $\Delta (\xx )$, which occurs as
the special case $\mu = (1^{n})$.

Given a partition $\mu $ of $n$, we denote by
\begin{equation}\label{e:n!/Dmu}
D_{\mu } = \CC [\partial \xx ,\partial \yy ]\Delta _{\mu }
\end{equation}
the space spanned by the partial derivatives of all orders of the
determinant $\Delta _{\mu }$ in \eqref{e:n!/Delta-mu}.  In
\cite{GaHa93} Garsia and the author proposed the following conjecture,
which we will prove as a consequence of Proposition
\ref{p:hilbert/equiv} and Theorem~\ref{t:main}.

\begin{conj}[$n!$ conjecture]\label{conj:n!/n!}
The dimension of $D_{\mu }$ is equal to $n!$.
\end{conj}


The $n!$ conjecture arose as part of a stronger conjecture relating
the Kostka-Macdonald coefficients to the character of $D_{\mu }$ as a
doubly graded $S_{n}$-module.  The symmetric group $S_{n}$ acts by
$\CC $-algebra automorphisms of $\CC [\xx ,\yy ]$ permuting the
variables:
\begin{equation}\label{e:n!/Sn-action}
wx_{i} = x_{w(i)},\; wy_{i} = y_{w(i)}\quad \text{for}\quad w\in
S_{n}.
\end{equation}
The ring $\CC [\xx ,\yy ] = \bigoplus _{r,s} \CC [\xx ,\yy ]_{r,s}$ is
doubly graded, by degree in the $\xx $ and $\yy $ variables
respectively, and the $S_{n}$ action respects the grading.  Clearly
$\Delta _{\mu }$ is $S_{n}$-alternating, {\it i.e.}, we have $w\Delta
_{\mu } = \epsilon (w)\Delta _{\mu }$ for all $w\in S_{n}$, where
$\epsilon $ is the sign character.  Note that $\Delta _{\mu }$ is also
doubly homogeneous, of $x$-degree $n(\mu )$ and $y$-degree $n(\mu ')$.
It follows that the space $D_{\mu }$ is $S_{n}$-invariant and has a
double grading
\begin{equation}\label{e:n!/Dmu-grading}
D_{\mu } = \bigoplus_{r,s} (D_{\mu })_{r,s}
\end{equation}
by $S_{n}$-invariant subspaces $(D_{\mu })_{r,s} = D_{\mu }\cap \CC
[\xx ,\yy ]_{r,s}$.


We write $\ch V$ for the character of an $S_{n}$-module $V$, and
denote the irreducible $S_{n}$ characters by $\chi ^{\lambda }$, with
the usual indexing by partitions $\lambda $ of $n$.  The following
conjecture clearly implies the Macdonald positivity conjecture.

\begin{conj}[\cite{GaHa93}]\label{conj:n!/GH}
We have
\begin{equation}\label{e:n!/GH}
\tilde{K}_{\lambda \mu }(q,t) = \sum _{r,s}t^{r}q^{s}\langle \chi
^{\lambda }, \ch (D_{\mu })_{r,s} \rangle.
\end{equation}
\end{conj}

Macdonald had shown that $K_{\lambda \mu }(1,1)$ is equal to $\chi
^{\lambda }(1)$, the degree of the irreducible $S_{n}$ character $\chi
^{\lambda }$, or the number of standard Young tableaux of shape
$\lambda $.  Conjecture \ref{conj:n!/GH} therefore implies that
$D_{\mu }$ affords the regular representation of $S_{n}$.  In
particular, it implies the $n!$ conjecture.

In \cite{Hai99} the author showed that Conjecture \ref{conj:n!/GH}
would follow from the Cohen-Macaulay property of $X_{n}$.  We
summarize the argument proving Conjecture \ref{conj:n!/GH} in \S
\ref{ss:hilbert/character}, after the relevant geometric results have
been established.


\section{The isospectral Hilbert scheme}
\label{s:hilbert}

\subsection{Preliminaries}
\label{ss:hilbert/prelim}

In this section we define the isospectral Hilbert scheme $X_{n}$, and
deduce our main theorem, Theorem~\ref{t:main} (\S
\ref{ss:hilbert/main}).  We also define the the Hilbert scheme $H_{n}$
and the {\it nested Hilbert scheme} $H_{n-1,n}$, and develop some
basic properties of these various schemes in preparation for the proof
of the main theorem.

The main technical device used in the proof of Theorem~\ref{t:main} is
a theorem on certain subspace arrangements called {\it polygraphs},
Theorem~\ref{t:polygraphs-are-free} in \S \ref{ss:poly/first-defs}.
The proof of the latter theorem is lengthy and logically distinct from
the geometric reasoning leading from there to Theorem~\ref{t:main}.
For these reasons we have deferred Theorem~\ref{t:polygraphs-are-free}
and its proof to the separate Section \ref{s:polygraphs}.

Throughout this section we work in the category of schemes of finite
type over the field of complex numbers, $\CC $.  All the specific
schemes we consider are quasiprojective over $\CC $.  We use
classical geometric language, describing open and closed subsets of
schemes, and morphisms between reduced schemes, in terms of closed
points.  A {\it variety} is a reduced and irreducible scheme.

Every locally free coherent sheaf $B$ of rank $n$ on a scheme $X$ of
finite type over $\CC $ is isomorphic to the sheaf of sections of an
algebraic vector bundle of rank $n$ over $X$.  For notational
purposes, we identify the vector bundle with the sheaf $B$ and write
$B(P)$ for the fiber of $B$ at a closed point $P\in X$.  In
sheaf-theoretic terms, the fiber is given by $B(P) = B\otimes_{\Ocal
_{X}} (\Ocal _{X,P}/P)$.

A scheme $X$ is {\it Cohen-Macaulay} or {\it Gorenstein},
respectively, if its local ring $\Ocal _{X,P}$ at every point is a
Cohen-Macaulay or Gorenstein local ring.  For either condition it
suffices that it hold at closed points $P$.  At the end of the
section, in \S \ref{ss:hilbert/duality}, we provide a brief summary of
the facts we need from duality theory and the theory of Cohen-Macaulay
and Gorenstein schemes.


\subsection{The schemes $H_{n}$ and $X_{n}$}
\label{ss:hilbert/Hn-Xn}

Let $R = \CC [x,y]$ be the coordinate ring of the affine plane $\CC
^{2}$.  By definition, closed subschemes $S\subseteq \CC ^{2}$ are in
one-to-one correspondence with ideals $I\subseteq R$.  The subscheme
$S=V(I)$ is finite if and only if $R/I$ has Krull dimension zero, or
finite dimension as a vector space over $\CC $.  In this case, the
{\it length} of $S$ is defined to be $\dim _{\CC }R/I$.

The Hilbert scheme $H_{n} = \Hilb ^{n}(\CC ^{2})$ parametrizes finite
closed subschemes $S\subseteq \CC ^{2}$ of length $n$.  The scheme
structure of $H_{n}$ and the precise sense in which it parametrizes
the subschemes $S$ are defined by a universal property, which
characterizes $H_{n}$ up to unique isomorphism.  The universal
property is actually a property of $H_{n}$ together with a closed
subscheme $F\subseteq H_{n}\times \CC ^{2}$, called the {\it universal
family}.

\begin{prop}\label{p:hilbert/universal}
There exist schemes $H_{n} = \Hilb ^{n}(\CC ^{2})$ and $F\subseteq
H_{n}\times \CC ^{2}$ enjoying the following properties, which
characterize them up to unique isomorphism:
\begin{itemize}
\item [(1)] $F$ is flat and finite of degree $n$ over $H_{n}$, and
\item [(2)] if $Y\subseteq T\times \CC ^{2}$ is a closed subscheme,
flat and finite of degree $n$ over a scheme $T$, then there is a
unique morphism $\phi \colon T\rightarrow H_{n}$ giving a commutative
diagram
\begin{equation*}		
\begin{CD}
Y &	@>>>&	T\times \CC ^{2}&	@>>>&	T \\
@VVV &	&	@VVV &			&	@V{\phi }VV\\
F&	@>>>&	H_{n}\times \CC ^{2}&	@>>>& 	H_{n},
\end{CD}
\end{equation*}
that is, the flat family $Y$ over $T$ is the pullback through $\phi $
of the universal family $F$.
\end{itemize}
\end{prop}


\begin{proof}
The Hilbert scheme $\hat{H} = \Hilb ^{n}(\PP ^{2})$ of points in the
projective plane exists as a special case of Grothendieck's
construction in \cite{Gro61}, with a universal family $\hat{F}$ having
the analogous universal property.  We identify $\CC ^{2}$ as usual
with an open subset of $\PP ^{2}$, the complement of the projective
line $Z$ ``at infinity.''

The projection of $\hat{F}\cap (\hat{H}\times Z)$ onto $\hat{H}$ is a
closed subset of $\hat{H}$.  Its complement $H_{n}\subseteq \hat{H}$
is clearly the largest subset such that the restriction $F$ of
$\hat{F}$ to $H_{n}$ is contained in $H_{n}\times \CC ^{2}$.  The
required universal property of $H_{n}$ and $F$ now follows immediately
from that of $\hat{H}$ and $\hat{F}$.
\end{proof}

To see how $H_{n}$ parametrizes finite closed subschemes $S\subseteq
\CC ^{2}$ of length $n$, note that the latter are exactly the families
$Y$ in Proposition \ref{p:hilbert/universal} for $T = \Spec \CC $.  By
the universal property they correspond one-to-one with the closed
points of $H_{n}$.  For notational purposes we will identify the
closed points of $H_{n}$ with ideals $I \subseteq R$ satisfying $\dim
_{\CC }R/I = n$, rather than with the corresponding subschemes
$S=V(I)$.


We have the following fundamental theorem of Fogarty \cite{Fog68}.

\begin{prop}\label{p:hilbert/fogarty}
The Hilbert scheme $H_{n}$ is a nonsingular, irreducible variety over
$\CC $ of dimension $2n$.
\end{prop}

The generic examples of finite closed subschemes $S\subseteq \CC ^{2}$
of length $n$ are the reduced subschemes consisting of $n$ distinct
points.  They form an open subset of $H_{n}$, and the irreducibility
aspect of Fogarty's theorem means that this open set is dense.

The most special closed subschemes in a certain sense are those
defined by monomial ideals.  If $I\subseteq R$ is a monomial ideal
then the {\it standard monomials} $x^{p}y^{q}\not \in I$ form a basis
of $R/I$.  If $\dim _{\CC }R/I = n$, the exponents $(p,q)$ of the
standard monomials form the diagram $D(\mu )$ of a partition $\mu $ of
$n$, and conversely.  We use the partition $\mu $ to index the
corresponding monomial ideal, denoting it by $I_{\mu }$.  Note that
$\rad{I_{\mu }} = (x,y)$ for all $\mu $, so the subscheme $V(I_{\mu
})$ is concentrated at the origin $(0,0)\in \CC ^{2}$, and owes its
length entirely to its non-reduced scheme structure.

The algebraic torus
\begin{equation}\label{e:hilbert/torus}
\TT ^{2} = (\CC ^{*})^{2}
\end{equation}
acts on $\CC ^{2}$ as the group of invertible diagonal $2\times 2$
matrices.  The monomial ideals $I_{\mu }$ are the torus invariant
ideals, and thus they are the fixed points of the induced action of
$\TT ^{2}$ on the Hilbert scheme.  Every ideal $I\in H_{n}$ has a
monomial ideal in the closure of its $\TT ^{2}$ orbit (\cite{Hai98} ,
Lemma 2.3).


Let $\CC [\xx ,\yy ] = \CC [x_{1},y_{1},\ldots,x_{n},y_{n}]$ be the
polynomial ring in $2n$ variables, so $\Spec \CC [\xx ,\yy ] = (\CC
^{2})^{n}$, where $x_{i}$, $y_{i}$ are the coordinates on the $i$-th
factor.  The symmetric group $S_{n}$ acts on $(\CC ^{2})^{n}$ by
permuting the Cartesian factors.  In coordinates, this corresponds to
the action of $S_{n}$ on $\CC [\xx ,\yy ]$ given in
\eqref{e:n!/Sn-action}.  We can identify $\Spec \CC [\xx ,\yy
]^{S_{n}}$ with the variety
\begin{equation}\label{e:hilbert/SnC2}
S^{n}\CC ^{2} = (\CC ^{2})^{n}/S_{n}
\end{equation}
of unordered $n$-tuples, or $n$-element multisets, of points in $\CC
^{2}$.

\begin{prop}[\cite{Hai98}, Proposition 2.2]\label{p:hilbert/chow}
For $I\in H_{n}$, let $\sigma (I)$ be the multiset of points of
$V(I)$, counting each point $P$ with multiplicity equal to the length
of the local ring $(R/I)_{P}$.  Then the map $\sigma \colon
H_{n}\rightarrow S^{n}\CC ^{2}$ is a projective morphism (called the
{\rm Chow morphism}).
\end{prop}

\begin{defn}\label{d:hilbert/Xn}
The {\it isospectral Hilbert scheme} $X_{n}$ is the reduced fiber
product
\begin{equation}\label{e:hilbert/Xn-diagram}
\begin{CD}
X_{n}&	@>{f}>>&	(\CC ^{2})^{n}\\
@V{\rho }VV&	&	@VVV\\
H_{n}&	@>{\sigma }>>&	S^{n}\CC ^{2},
\end{CD}
\end{equation}
that is, the reduced closed subscheme of $H_{n}\times (\CC ^{2})^{n}$
whose closed points are the tuples $(I,P_{1},\ldots,P_{n})$ satisfying
$\sigma (I) = \llbracket P_{1},\ldots,P_{n} \rrbracket$.
\end{defn}

We will continue to refer to the morphisms $\rho $, $\sigma $ and $f$
in diagram \eqref{e:hilbert/Xn-diagram} by those names in what
follows.


For each $I\in H_{n}$, the operators $\overline{x}$, $\overline{y}$ of
multiplication by $x,y$ are commuting endomorphisms of the
$n$-dimensional vector space $R/I$.  As such, they have a well-defined
joint spectrum, a multiset of pairs of eigenvalues
$(x_{1},y_{1}),\ldots,(x_{n},y_{n})$ determined by the identity
\begin{equation}\label{e:hilbert/spectrum}
\det \nolimits _{R/I}(1+\alpha \overline{x}+\beta \overline{y}) =
\prod _{i=1}^{n}(1+\alpha x_{i}+\beta y_{i}).
\end{equation}
On the local ring $(R/I)_{P}$ at a point $P = (x_{0},y_{0})$, the
operators $\overline{x}$, $\overline{y}$ have the sole joint
eigenvalue $(x_{0},y_{0})$, with multiplicity equal to the length of
$(R/I)_{P}$.  Hence $\sigma (I)$ is equal as a multiset to the joint
spectrum of $\overline{x}$ and $\overline{y}$.  This is the motivation
for the term {\it isospectral}.

The action of $S_{n}$ on $(\CC ^{2})^{n}$ induces a compatible action
of $S_{n}$ on $X_{n}$ by automorphisms of $X_{n}$ as a scheme over
$H_{n}$.  Explicitly, for $w\in S_{n}$ we have
$w(I,P_{1},\ldots,P_{n}) = (I,P_{w^{-1}(1)},\ldots,P_{w^{-1}(n)})$.

We caution the reader that the {\it scheme-theoretic} fiber product in
\eqref{e:hilbert/Xn-diagram} is not reduced, even for $n=2$.  For
every invariant polynomial $g\in \CC [\xx ,\yy ]^{S_{n}}$,
the global regular function
\begin{equation}\label{e:hilbert/not-ideal}
g(x_{1},y_{1},\ldots,x_{n},y_{n}) - \sigma ^{\sharp }g
\end{equation}
on $H_{n}\times (\CC ^{2})^{n}$ vanishes on $X_{n}$.  By definition
these equations generate the ideal sheaf of the scheme-theoretic fiber
product.  They cut out $X_{n}$ set-theoretically, but not as a reduced
subscheme.  The full ideal sheaf defining $X_{n}$ as a reduced scheme
must necessarily have a complicated local description, since it is a
consequence of Theorem~\ref{t:main} and Proposition
\ref{p:hilbert/equiv} (to follow) that generators for all the ideals
$J_{\mu }$ in \S \ref{ss:hilbert/n!}, eq.~\eqref{e:hilbert/Jmu} are
implicit in the local ideals of $X_{n}$ at the distinguished points
$Q_{\mu }$ lying over the fixed points $I_{\mu }\in H_{n}$.


\subsection{Elementary properties of $X_{n}$}
\label{ss:hilbert/simple}

We now develop some elementary facts about the isospectral Hilbert
scheme $X_{n}$.  The first of these is its product structure, which
reduces local questions on $X_{n}$ to the case of smaller instances
$X_{k}$ $(k<n)$ at any point whose corresponding multiset $\llbracket
P_{1},\ldots, P_{n} \rrbracket$ is not of the form $\llbracket n\cdot
P \rrbracket$.

\begin{lem}\label{l:hilbert/product}
Let $k$ and $l$ be positive integers with $k+l = n$.  Suppose
$U\subseteq (\CC ^{2})^{n}$ is an open set consisting of points
$(P_{1},\ldots,P_{k},Q_{1},\ldots,Q_{l})$ where no $P_{i}$ coincides
with any $Q_{i}$.  Then, identifying $(\CC ^{2})^{n}$ with $(\CC
^{2})^{k}\times (\CC ^{2})^{l}$, the preimage $f^{-1}(U)$ of $U$ in
$X_{n}$ is isomorphic as a scheme over $(\CC ^{2})^{n}$ to the
preimage $(f_{k}\times f_{l})^{-1}(U)$ of $U$ in $X_{k}\times X_{l}$.
\end{lem}

\begin{proof}
Let $Y = (\rho \times 1)^{-1}(F)\subseteq X_{n}\times \CC ^{2}$ be the
universal family over $X_{n}$.  The fiber $V(I)$ of $Y$ over a point
$(I,P_{1}, \ldots, P_{k}, Q_{1}, \ldots, Q_{l})\in f^{-1}(U)$ is the
disjoint union of closed subschemes $V(I_{k})$ and $V(I_{l})$ in $\CC
^{2}$ of lengths $k$ and $l$, respectively, with $\sigma (I_{k}) =
\llbracket P_{1},\ldots,P_{k} \rrbracket$ and $\sigma (I_{l}) =
\llbracket Q_{1},\ldots,Q_{l} \rrbracket$.  Hence over $f^{-1}(U)$,
$Y$ is the disjoint union of flat families $Y_{k}$, $Y_{l}$ of degrees
$k$ and $l$.  By the universal property, we get induced morphisms
$\phi _{k}\colon f^{-1}(U)\rightarrow H_{k}$, $\phi _{l}\colon
f^{-1}(U)\rightarrow H_{l}$ and $\phi _{k}\times \phi _{l}\colon
f^{-1}(U)\rightarrow H_{k}\times H_{l}$.  The equations $\sigma
(I_{k}) = \llbracket P_{1},\ldots,P_{k} \rrbracket$, $\sigma (I_{l}) =
\llbracket Q_{1},\ldots,Q_{l} \rrbracket$ imply that $\phi _{k}\times
\phi _{l}$ factors through a morphism $\alpha \colon
f^{-1}(U)\rightarrow X_{k}\times X_{l}$ of schemes over $(\CC
^{2})^{n}$.


Conversely, on $(f_{k}\times f_{l})^{-1}(U)\subseteq X_{k}\times
X_{l}$, the pullbacks of the universal families from $X_{k}$ and
$X_{l}$ are disjoint and their union is a flat family of degree $n$.
By the universal property there is an induced morphism $\psi \colon
(f_{k}\times f_{l})^{-1}(U)\rightarrow H_{n}$, which factors through a
morphism $\beta \colon (f_{k}\times f_{l})^{-1}(U)\rightarrow X_{n}$
of schemes over $(\CC ^{2})^{n}$.

By construction, the universal families on $f^{-1}(U)$ and
$(f_{k}\times f_{l})^{-1}(U)$ pull back to themselves via $\beta \circ
\alpha $ and $\alpha \circ \beta $, respectively.  This implies that
$\beta \circ \alpha $ is a morphism of schemes over $H_{n}$ and
$\alpha \circ \beta $ is a morphism of schemes over $H_{k}\times
H_{l}$.  Since they are also morphisms of schemes over $(\CC
^{2})^{n}$, we have $\beta \circ \alpha = 1_{f^{-1}(U)}$ and $\alpha
\circ \beta = 1_{(f_{k}\times f_{l})^{-1}(U)}$.  Hence $\alpha $ and
$\beta$ induce mutually inverse isomorphisms $f^{-1}(U)\cong
(f_{k}\times f_{l})^{-1}(U)$.
\end{proof}


\begin{prop}\label{p:hilbert/irreducible}
The isospectral Hilbert scheme $X_{n}$ is irreducible, of dimension $2n$.
\end{prop}

\begin{proof}
Let $U$ be the preimage in $X_{n}$ of the open set $W\subseteq (\CC
^{2})^{n}$ of points $(P_{1},\ldots,P_{n})$ where the $P_{i}$ are all
distinct.  It follows from Lemma \ref{l:hilbert/product} that $f$
restricts to an isomorphism $f\colon U\rightarrow W$, so $U$ is
irreducible.  We are to show that $U$ is dense in $X_{n}$.

Let $Q$ be a closed point of $X_{n}$, which we want to show belongs to
the closure $\overline{U}$ of $U$.  If $f(Q) = (P_{1},\ldots,P_{n})$
with the $P_{i}$ not all equal, then by Lemma \ref{l:hilbert/product}
there is a neighborhood of $Q$ in $X_{n}$ isomorphic to an open set in
$X_{k}\times X_{l}$ for some $k, l < n$.  The result then follows by
induction, since we may assume $X_{k}$ and $X_{l}$ irreducible.  If
all the $P_{i}$ are equal, then $Q$ is the unique point of $X_{n}$
lying over $I = \rho (Q)\in H_{n}$.  Since $\rho $ is finite, $\rho
(\overline{U})\subseteq H_{n}$ is closed.  But $\rho (U)$ is dense in
$H_{n}$, so $\rho (\overline{U}) = H_{n}$.  Therefore $\overline{U}$
contains a point lying over $I$, which must be $Q$.
\end{proof}

\begin{prop}\label{p:hilbert/V(y)-in-X}
The closed subset $V(y_{1},\ldots,y_{n})$ in $X_{n}$ has dimension $n$.
\end{prop}

\begin{proof}
It follows from the cell decomposition of Ellingsrud and Str{\o}mme
\cite{ElSt87} that the closed locus $Z$ in $H_{n}$ consisting of
points $I$ with $V(I)$ supported on the $x$-axis $V(y)$ in $\CC ^{2}$
is the union of locally closed affine cells of dimension $n$.  The
subset $V(\yy )\subseteq X_{n}$ is equal to $\rho ^{-1}(Z)$ and $\rho
$ is finite.
\end{proof}


The product structure of $X_{n}$ is inherited in a certain sense by
$H_{n}$, but its description in terms of $X_{n}$ is more transparent.
As a consequence, passage to $X_{n}$ is sometimes handy for proving
results purely about $H_{n}$.  The following lemma is an example of
this.  We remark that one can show by a more careful analysis that
locus $G_{r}$ in the lemma is in fact irreducible.

\begin{lem}\label{l:hilbert/Gr}
Let $G_{r}$ be the closed subset of $H_{n}$ consisting of ideals $I$
for which $\sigma (I)$ contains some point with multiplicity at least
$r$.  Then $G_{r}$ has codimension $r-1$, and has a unique irreducible
component of maximal dimension.
\end{lem}

\begin{proof}
By symmetry among the points $P_{i}$ of $\sigma (I)$ we see that
$G_{r} = \rho (V_{r})$, where $V_{r}$ is the locus in $X_{n}$ defined
by the equations $P_{1}=\cdots =P_{r}$.  It follows from Lemma
\ref{l:hilbert/product} that $V_{r}\setminus \rho ^{-1}(G_{r+1})$ is
isomorphic to an open set in $W_{r}\times X_{n-r}$, where $W_{r}$ is
the closed subset $P_{1}=\cdots =P_{r}$ in $X_{r}$.  As a reduced
subscheme of $X_{r}$, the latter is isomorphic to $\CC ^{2}\times
Z_{r}$, where $Z_{r} = \sigma ^{-1}(\{0 \})$ is the {\it zero fiber}
in $H_{r}$,  the factor $\CC ^{2}$ accounting for the choice of $P =
P_{1}=\cdots =P_{r}$.

By a theorem of Brian\c{c}on \cite{Bri77}, $Z_{r}$ is irreducible of
dimension $r-1$, so $V_{r}\setminus \rho ^{-1}(G_{r+1})$ is
irreducible of dimension $2(n-r)+r+1 = 2n-(r-1)$.  Since
$G_{r}\setminus G_{r+1} = \rho (V_{r}\setminus \rho ^{-1}(G_{r+1}))$
and $\rho $ is finite, the result follows by descending induction on
$r$, starting with $G_{n+1} = \emptyset $.
\end{proof}


\subsection{Blowup construction of $H_{n}$ and $X_{n}$}
\label{ss:hilbert/blowup}

Let $A = \CC [\xx ,\yy ]^{\epsilon }$ be the space of
$S_{n}$-alternating elements, that is, polynomials $g$ such that $wg =
\epsilon (w)g$ for all $w\in S_{n}$, where $\epsilon $ is the sign
character.  To describe $A$ more precisely, we note that $A$ is the
image of the alternation operator
\begin{equation}\label{e:hilbert/alt}
\Theta ^{\epsilon }g = \sum _{w\in S_{n}} \epsilon (w) w g.
\end{equation}
If $D = \{ (p_{1},q_{1}),\ldots,(p_{n},q_{n})\}$ is an $n$-element
subset of $\NN \times \NN $, then the determinant $\Delta _{D}$
defined in \eqref{e:n!/Delta-D} can also be written
\begin{equation}\label{e:hilbert/Delta-D}
\Delta _{D} = \Theta ^{\epsilon }(\xx ^{p}\yy ^{q}),
\end{equation}
where $\xx ^{p}\yy ^{q} = x_{1}^{p_{1}}y_{1}^{q_{1}}\cdots
x_{n}^{p_{n}}y_{n}^{q_{n}}$.  For a monomial $\xx ^{p} \yy ^{q}$ whose
exponent pairs $(p_{i},q_{i})$ are not all distinct, we have $\Theta
^{\epsilon }(\xx ^{p}\yy ^{q}) = 0$.  From this it is easy to see that
the set of all elements $\Delta _{D}$ is a basis of $A$.  Another way
to see this is to identify $A$ with the $n$-th exterior power $\ext
^{n}\CC [x,y]$ of the polynomial ring in two variables $x,y$.  Then
the basis elements $\Delta _{D}$ are identified with the wedge
products of monomials in $\CC [x,y]$.

For $d>0$, let $A^{d}$ be the space spanned by all products of $d$
elements of $A$.  We set $A^{0} = \CC [\xx ,\yy ]^{S_{n}}$.  Note that
$A$ and hence every $A^{d}$ is a $\CC [\xx ,\yy ]^{S_{n}}$-submodule
of $\CC [\xx ,\yy ]$, so we have $A^{i}A^{j} = A^{i+j}$ for all $i$, $j$,
including $i=0$ or $j=0$.


\begin{prop}[\cite{Hai98}, Proposition 2.6]\label{p:hilbert/blowup-H}
The Hilbert scheme $H_{n}$ is isomorphic as a scheme projective over
$S^{n}\CC ^{2}$ to $\Proj T$, where $T$ is the graded $\CC [\xx ,\yy
]^{S_{n}}$-algebra $T = \bigoplus _{d\geq 0} A^{d}$.
\end{prop}

Now set $S = \CC [\xx ,\yy ]$ and let $J = SA$ be the ideal generated
by $A$.

\begin{prop}\label{p:hilbert/blowup}
The isospectral Hilbert scheme $X_{n}$ is isomorphic as a scheme over
$(\CC ^{2})^{n}$ to the blowup of $(\CC ^{2})^{n}$ at the ideal $J =
SA$ defined above.
\end{prop}

\begin{proof}
By definition the blowup of $(\CC ^{2})^{n}$ at $J$ is $Z = \Proj
S[tJ]$, where $S[tJ]\cong \bigoplus _{d\geq 0} J^{d}$ is the Rees
algebra.  The ring $T$ is a homogeneous subring of $S[tJ]$ in an
obvious way, and since $A^{d}$ generates $J^{d}$ as a $\CC [\xx ,\yy
]$-module, we have $\CC [\xx ,\yy ]T = S[tJ]$, that is, $S[tJ]\cong
(\CC [\xx ,\yy ]\otimes _{A^{0}}T)/I$ for some homogeneous ideal $I$.
In geometric terms, using Proposition \ref{p:hilbert/blowup-H} and the
fact that $S^{n}\CC ^{2} = \Spec A^{0}$, this says that $Z$ is a
closed subscheme of the scheme-theoretic fiber product $H_{n}\times
_{S^{n}\CC ^{2}}(\CC ^{2})^{n}$.  Since $Z$ is reduced, it follows that
$Z$ is a closed subscheme of $X_{n}$.  By Proposition
\ref{p:hilbert/irreducible}, $X_{n}$ is irreducible, and since both
$Z$ and $X_{n}$ have dimension $2n$, it follows that $Z=X_{n}$.
\end{proof}

In the context of either $H_{n}$ or $X_{n}$ we will always write
$\Ocal (k)$ for the $k$-th tensor power of the ample line bundle
$\Ocal (1)$ induced by the representation of $H_{n}$ as $\Proj T$ or
$X_{n}$ as $\Proj S[tJ]$.  It is immediate from the proof of
Proposition \ref{p:hilbert/blowup} that $\Ocal _{X_{n}}(k) = \rho
^{*}\Ocal _{H_{n}}(k)$.


In full analogy to the situation for the Pl\"ucker embedding of a
Grassmann variety, there is an intrinsic description of $\Ocal (1)$ as
the highest exterior power of the {\it tautological vector bundle}
whose fiber at a point $I\in H_{n}$ is $R/I$.  Let
\begin{equation}\label{e:hilbert/pi:F->H}
\pi \colon F\rightarrow H_{n}
\end{equation}
be the projection of the universal family on the Hilbert scheme.
Since $\pi $ is an affine morphism, we have $F = \Spec B$, where $B$
is the sheaf of $\Ocal _{H_{n}}$-algebras
\begin{equation}\label{e:hilbert/B}
B = \pi _{*}\Ocal _{F}.
\end{equation}
The fact that $F$ is flat and finite of degree $n$ over $H_{n}$ means
that $B$ is a locally free sheaf of $\Ocal _{H_{n}}$-modules of rank
$n$.  Its associated vector bundle is the tautological bundle.

\begin{prop}[\cite{Hai98}, Proposition 2.12]\label{p:hilbert/O(1)}
We have an isomorphism $\ext ^{n} B \cong \Ocal (1)$ of line bundles
on $H_{n}$.
\end{prop}


\subsection{Nested Hilbert schemes}
\label{ss:hilbert/nested}

The proof of our main theorem will be by induction on $n$.
For the inductive step we interpolate between $H_{n-1}$ and $H_{n}$ using
the {\it nested} Hilbert scheme.

\begin{defn}\label{d:hilbert/nested}
The {\it nested Hilbert scheme} $H_{n-1,n}$ is the reduced closed subscheme
\begin{equation}\label{e:hilbert/nested}
H_{n-1,n} = \{(I_{n-1},I_{n}):I_{n}\subseteq I_{n-1} \}\subseteq
H_{n-1}\times H_{n}.
\end{equation}
\end{defn}

The analog of Fogarty's theorem (Proposition \ref{p:hilbert/fogarty})
for the nested Hilbert scheme is the following result of Tikhomirov,
whose proof can be found in \cite{Chea98}.

\begin{prop}\label{p:hilbert/tikh-cheah}
The nested Hilbert scheme $H_{n-1,n}$ is nonsingular and
irreducible, of dimension $2n$.
\end{prop}


As with $H_{n}$, the nested Hilbert scheme is an open set in a
projective nested Hilbert scheme $\Hilb ^{n-1,n}(\PP ^{2})$.  Clearly,
$H_{n-1,n}$ is the preimage of $H_{n}$ under the projection $\Hilb
^{n-1,n}(\PP ^{2})\rightarrow \Hilb ^{n}(\PP ^{2})$.  Hence the
projection $H_{n-1,n}\rightarrow H_{n}$ is a projective morphism.

If $(I_{n-1},I_{n})$ is a point of $H_{n-1,n}$ then $\sigma (I_{n-1})$
is an $n-1$ element sub-multiset of $\sigma (I_{n})$.  In symbols, if
$\sigma (I_{n-1}) = \llbracket P_{1},\ldots,P_{n-1} \rrbracket $, then
$\sigma (I_{n}) = \llbracket P_{1},\ldots,P_{n-1},P_{n} \rrbracket $
for a distinguished last point $P_{n}$.  The $S_{n-1}$-invariant
polynomials in the coordinates $x_{1},y_{1},\ldots,x_{n-1},y_{n-1}$ of
the points $P_{1},\ldots,P_{n-1}$ are global regular functions on
$H_{n-1,n}$, pulled back via the projection on $H_{n-1}$.  Similarly,
the $S_{n}$-invariant polynomials in $x_{1},y_{1},\ldots,x_{n},y_{n}$
are regular functions pulled back from $H_{n}$.  It follows that the
coordinates of the distinguished point $P_{n}$ are regular functions,
since $x_{n} = (x_{1}+\cdots +x_{n})-(x_{1}+\cdots +x_{n-1})$, and
similarly for $y_{n}$.  Thus we have a morphism
\begin{equation}\label{e:hilbert/nested-sigma}
\sigma \colon H_{n-1,n}\rightarrow S^{n-1}\CC ^{2}\times \CC ^{2} =
(\CC ^{2})^{n}/S_{n-1}, 
\end{equation}
such that both the maps $H_{n-1,n}\rightarrow S^{n-1}\CC ^{2}$ and
$H_{n-1,n}\rightarrow S^{n}\CC ^{2}$ induced by the Chow morphisms
composed with the projections on $H_{n-1}$ and $H_{n}$ factor through
$\sigma $.

The distinguished point $P_{n}$ belongs to $V(I_{n})$, and can be any
point of $V(I_{n})$, so the image of the morphism
\begin{equation}\label{e:hilbert/nested->F}
\alpha \colon H_{n-1,n}\rightarrow H_{n}\times \CC ^{2}
\end{equation}
sending $(I_{n-1},I_{n})$ to $(I_{n},P_{n})$ is precisely the
universal family $F$ over $H_{n}$.

The following proposition, in conjunction with Lemma
\ref{l:hilbert/Gr}, provides dimension estimates needed for the
calculation of the canonical line bundle on $H_{n-1,n}$ in \S
\ref{ss:hilbert/omega} and the proof of the main theorem in \S
\ref{ss:hilbert/main}.


\begin{prop}\label{p:hilbert/fibers}
Let $d$ be the dimension of the fiber of the morphism $\alpha $ in
\eqref{e:hilbert/nested->F} over a point $(I,P)\in F$, and let $r$ be
the multiplicity of $P$ in $\sigma (I)$.  Then $d$ and $r$ satisfy the
inequality
\begin{equation}\label{e:hilbert/fiber-inequality}
r\geq \binom{d+2}{2}.
\end{equation}
\end{prop}

\begin{proof}
The possible ideals $I_{n-1}$ for the given $(I_{n},P_{n}) = (I,P)$
are the length $1$ ideals in the local ring $(R/I)_{P}$, or
equivalently, the one-dimensional subspaces of its socle.  The fiber
of $\alpha $ is therefore the projective space $\PP (\soc (R/I)_{P})$,
and we have $d+1 = \dim \soc (R/I)_{P}$.

First consider the maximum possible dimension of any fiber of $\alpha
$.  Since both $H_{n-1,n}$ and $F$ are projective over $H_{n}$, the
morphism $\alpha $ is projective and its fiber dimension is upper
semicontinuous.  Since every point of $H_{n}$ has a monomial ideal
$I_{\mu }$ in the closure of its $\TT ^{2}$-orbit, and since $F$ is
finite over $H_{n}$, every point of $F$ must have a pair $(I_{\mu
},0)\in F$ in the closure of its orbit.  The fiber dimension is
therefore maximized at some such point.  The socle of $R/I_{\mu }$ has
dimension equal to the number of corners of the diagram of $\mu $.  If
this number is $s$, we clearly have $n\geq \binom{s+1}{2}$.  This
implies that for every Artin local ring $R/I$ generated
over $\CC $ by two elements, the socle dimension $s$ and the length
$n$ of $R/I$ satisfy $n\geq \binom{s+1}{2}$.

Returning to the original problem, $(R/I)_{P}$ is an Artin local
ring of length $r$ generated by two elements, with socle dimension
$d+1$, so \eqref{e:hilbert/fiber-inequality} follows.
\end{proof}


We now introduce the nested version of the isospectral Hilbert scheme.
It literally plays a pivotal role in the proof of the main theorem by
induction on $n$: we transfer the Gorenstein property from $X_{n-1}$
to the nested scheme $X_{n-1,n}$ by pulling back, and from there to
$X_{n}$ by pushing forward.

\begin{defn}\label{d:hilbert/nested-X}
The {\it nested isospectral Hilbert scheme} $X_{n-1,n}$ is the reduced
fiber product $H_{n-1,n}\times _{H_{n-1}} X_{n-1}$.
\end{defn}

There is an alternative formulation of the definition, which is useful
to keep in mind for the next two results.  Namely, $X_{n-1,n}$ can be
identified with the reduced fiber product in the diagram
\begin{equation}\label{e:hilbert/Xn-1,n-diagram}
\begin{CD}
X_{n-1,n}&	@>>>&	(\CC ^{2})^{n}\\
@VVV&		&	@VVV\\
H_{n-1,n}&	@>{\sigma }>>&	S^{n-1}\CC ^{2}\times \CC ^{2},
\end{CD}
\end{equation}
that is, the reduced closed subscheme of $H_{n-1,n}\times (\CC
^{2})^{n}$ consisting of tuples $(I_{n-1},I_{n},P_{1},\ldots,P_{n})$
such that $\sigma (I_{n}) = \llbracket P_{1},\ldots,P_{n} \rrbracket $
and $P_{n}$ is the distinguished point.  To see that this agrees with
the definition, note that a point of $H_{n-1,n}\times _{H_{n-1}}
X_{n-1}$ is given by the data $(I_{n-1},I_{n},P_{1},\ldots,P_{n-1})$,
and that these data determine the distinguished point $P_{n}$.  We
obtain the alternative description by identifying $X_{n-1,n}$ with the
graph in $X_{n-1,n}\times \CC ^{2}$ of the morphism
$X_{n-1,n}\rightarrow \CC ^{2}$ sending
$(I_{n-1},I_{n},P_{1},\ldots,P_{n-1})$ to $P_{n}$.


We have the following nested analogs of Lemma \ref{l:hilbert/product}
and Proposition \ref{p:hilbert/V(y)-in-X}.  The analog of Proposition
\ref{p:hilbert/irreducible} also holds, {\it i.e.}, $X_{n-1,n}$ is
irreducible.  We do not prove this here, as it will follow
automatically as part of our induction: see the observations following
diagram \eqref{e:hilbert/generic-square} in \S \ref{ss:hilbert/main}.

\begin{lem}\label{l:hilbert/product-nested}
Let $k+l = n$ and $U\subseteq (\CC ^{2})^{n}$ be as in Lemma
\ref{l:hilbert/product}.  Then the preimage of $U$ in $X_{n-1,n}$ is
isomorphic as a scheme over $(\CC ^{2})^{n}$ to the preimage of $U$ in
$X_{k}\times X_{l-1,l}$.
\end{lem}

\begin{proof}
Lemma \ref{l:hilbert/product} gives us isomorphisms on the preimage of
$U$ between $X_{n}$ and $X_{k}\times X_{l}$, and between $X_{n-1}$ and
$X_{k}\times X_{l-1}$.

We can identify $X_{n-1,n}$ with the closed subset of $X_{n-1}\times
X_{n}$ consisting of points where $P_{1},\ldots,P_{n-1}$ are the same
in both factors, and $I_{n-1}$ contains $I_{n}$.  On the preimage of
$U$, under the isomorphisms above, this corresponds to the closed
subset of $(X_{k}\times X_{l-1})\times (X_{k}\times X_{l})$ where
$P_{1},\ldots,P_{k}$, $Q_{1},\ldots,Q_{l-1}$ and $I_{k}$ are the same
in both factors and $I_{l-1}$ contains $I_{l}$.  The latter can be
identified with $X_{k}\times X_{l-1,l}$.
\end{proof}


\begin{prop}\label{p:hilbert/V(y)-in-nest-X}
The closed subset $V(y_{1},\ldots,y_{n})$ in $X_{n-1,n}$ has dimension
$n$.
\end{prop}

\begin{proof}
We have the corresponding result for $X_{n}$ in Proposition
\ref{p:hilbert/V(y)-in-X}.  We can assume by induction that the result
for the nested scheme holds for smaller values of $n$ (for the base
case note that $X_{0,1} \cong X_{1} \cong \CC ^{2})$.  Locally on a
neighborhood of any point where $P_{1},\ldots,P_{n}$ are not all
equal, the result then follows from Lemma
\ref{l:hilbert/product-nested}.

The locus where all the $P_{i}$ are equal is isomorphic to $\CC ^{1}
\times Z$ where $Z = \sigma ^{-1}(\{0 \})\subseteq H_{n-1,n}$ is the
zero fiber in the nested Hilbert scheme, the factor $\CC ^{1}$
accounting for the choice of the common point $P=P_{1}=\cdots =P_{n}$
on the $x$-axis $V(y)\subseteq \CC ^{2}$.  By a theorem of Cheah
(\cite{Chea98}, Theorem 3.3.3, part (5)) we have $\dim Z = n-1$.
\end{proof}


\subsection{Calculation of canonical line bundles}
\label{ss:hilbert/omega}

We will need to know the canonical sheaves $\omega $ on the smooth
schemes $H_{n}$ and $H_{n-1,n}$.  To compute them we make use of the
fact that invertible sheaves on a normal variety are isomorphic if
they have isomorphic restrictions to an open set whose complement has
codimension at least $2$.

\begin{defn}\label{e:hilbert/Ux}
Let $z = a x + b y$ be a linear form in the variables $x,y$.  We
denote by $U_{z}$ the open subset of $H_{n}$ consisting of ideals $I$
for which $z$ generates $R/I$ as an algebra over $\CC $.  We also
denote by $U_{z}$ the preimage of $U_{z}$ under the projection
$H_{n-1,n}\rightarrow H_{n}$.
\end{defn}

Note that $z$ generates $R/I$ if and only if the set
$\{1,z,\ldots,z^{n-1} \}$ is linearly independent, and thus a basis, in
$R/I$.  That given sections of a vector bundle determine linearly
independent elements in its fiber at a point $I$ is an open condition
on $I$, so $U_{z}$ is indeed an open subset of $H_{n}$.  The ring
$R/I$ can be generated by a single linear form $z$ if and only if the
scheme $V(I)$ is a subscheme of a smooth curve in $\CC ^{2}$.  For
this reason the union of all the open sets $U_{z}$ is called the {\it
curvilinear locus}.


\begin{lem}\label{l:hilbert/UxUy}
The complement of $U_{x}\cup U_{y}$ in has codimension $2$, both in
$H_{n}$ and in $H_{n-1,n}$.
\end{lem}

\begin{proof}
Let $Z = H_{n}\setminus (U_{x}\cup U_{y})$.  Let $W$ be the generic
locus, that is, the open set of ideals $I\in H_{n}$ for which $\sigma
(I) = \{P_{1},\ldots,P_{n} \}$ is a set of $n$ distinct points.  The
Chow morphism $\sigma $ induces an isomorphism of $W$ onto its image
in $S^{n}\CC ^{2}$, and $Z\cap W$ is the locus where some two of the
$P_{i}$ have the same $x$-coordinate, and another two have the same
$y$-coordinate.  This locus has codimension $2$.  The complement of
$W$ is the closed subset $G_{2}$ in Lemma \ref{l:hilbert/Gr}, which
has one irreducible component of dimension $2n-1$.  An open set in
this component consists of those $I$ for which $\sigma (I)$ has one
point $P_{i}$ of multiplicity $2$ and the rest are distinct.  This
open set is not contained in $Z$, so $Z\cap G_{2}$ has codimension at
least $2$.  This takes care of $H_{n}$.

If $I$ is curvilinear, then the socle of $(R/I)_{P}$ has length $1$
for all $P\in V(I)$.  Hence the morphism $\alpha \colon
H_{n-1,n}\rightarrow F$ in \eqref{e:hilbert/nested->F} restricts to a
bijection on the curvilinear locus.  If the complement of $U_{x}\cup
U_{y}$ in $H_{n-1,n}$ had a codimension $1$ component, it would
therefore have to be contained in the complement of the curvilinear
locus, by the result for $H_{n}$.

By Proposition \ref{p:hilbert/fibers}, $\alpha $ has fibers of
dimension $d$ only over $G_{r}$ for $r\geq \binom{d+2}{2}$, and it
follows from Lemma \ref{l:hilbert/Gr} that the union of these fibers
has codimension at least $\binom{d+2}{2}-1-d = \binom{d+1}{2}$.  For
$d>1$ this exceeds $2$.  The fibers of dimension $1$ over $G_{3}$
occur only over non-curvilinear points.  However, the codimension $2$
component of $G_{3}$ contains all $I$ such that $\sigma (I)$ has one
point of multiplicity $3$ and the rest distinct, and such an $I$ can
be curvilinear.  Hence the non-curvilinear locus in $G_{3}$ has
codimension at least $3$ and its preimage in $H_{n-1,n}$ has
codimension at least $2$.
\end{proof}


The following proposition is well-known; it holds for the Hilbert
scheme of points on any smooth surface with trivial canonical sheaf.
We give an elementary proof for $H_{n}$, since we need it as a
starting point for the proof of the corresponding result for
$H_{n-1,n}$.

\begin{prop}\label{p:hilbert/omega-Hn}
The canonical sheaf $\omega _{H_{n}}$ on the Hilbert scheme is
trivial, {\it i.e.}, $\omega _{H_{n}}\cong \Ocal _{H_{n}}$.
\end{prop}

\begin{proof}
The $2n$-form $d\xx \, d\yy  = dx_{1}\wedge \cdots \wedge dx_{n}\wedge
dy_{1}\wedge \cdots \wedge dy_{n}$ is $S_{n}$-invariant and thus defines
a $2n$-form on the smooth locus in $S^{n}\CC ^{2}$ and therefore a
rational $2n$-form on $H_{n}$.

If $I$ is a point of $U_{x}$, then $I$ is generated as an ideal in $R$
by two polynomials,
\begin{equation}\label{e:hilbert/x-generator}
x^{n}-e_{1}x^{n-1}+e_{2}x^{n-2}-\cdots +(-1)^{n}e_{n}
\end{equation}
and
\begin{equation}\label{e:hilbert/y-generator}
y - (a_{n-1}x^{n-1}+a_{n-2}x^{n-2}+\cdots +a_{0}),
\end{equation}
where the complex numbers $e_{1},\ldots,e_{n},a_{0},\ldots,a_{n-1}$
are regular functions of $I$.  This is so because the tautological
sheaf $B$ is free with basis $1,x,\ldots,x^{n-1}$ on $U_{x}$, and the
sections $x^{n}$ and $y$ must be unique linear combinations of the
basis sections with regular coefficients.  Conversely, for any choice
of the parameters ${\bf e}$, ${\aa}$, the polynomials in
\eqref{e:hilbert/x-generator} and \eqref{e:hilbert/y-generator}
generate an ideal $I\in U_{x}$.  Hence $U_{x}$ is an affine $2n$-cell
with coordinates ${\bf e},\aa $.


On the open set where $\sigma (I)$ consists of points $P_{i}$ with
distinct $x$ coordinates, the polynomial in
\eqref{e:hilbert/x-generator} is  $\prod _{i = 1}^{n} (x-x_{i})$.
This implies that $e_{k}$ coincides as a rational function (and
as a global regular function) with the $k$-th elementary symmetric
function $e_{k}(\xx )$.  Likewise, $a_{k}$ is given as a rational
function of $\xx $ and $\yy $ by the coefficient of $x^{k}$ in the
unique interpolating polynomial $\phi _{a}(x)$ of degree $n-1$
satisfying $\phi_{a} (x_{i}) = y_{i}$ for $i = 1,\ldots,n$.

The equations $\phi_{a} (x_{i}) = y_{i}$ can be expressed as a matrix
identity $(y_{1},\ldots,y_{n}) = (a_{0},\ldots,a_{n-1})M$, where $M$
is the Vandermonde matrix in the $\xx $ variables.  Modulo terms
involving the $dx_{i}$, this yields the identity of rational $n$-forms
on $H_{n}$,
\begin{equation}\label{e:hilbert/da}
d\aa  = \Delta (\xx )^{-1}\, d\yy,
\end{equation}
where $\Delta (\xx )$ is the Vandermonde determinant.
For the elementary symmetric functions $e_{k} = e_{k}(\xx )$ we have
the well-known identity
\begin{equation}\label{e:hilbert/de}
d{\bf e} = \Delta (\xx )\, d\xx .
\end{equation}
Together these show that $d\xx \, d\yy = d\aa \, d{\bf e}$ is a nowhere
vanishing regular section of $\omega $ on $U_{x}$.  By symmetry, the
same holds on $U_{y}$.  This shows that we have $\omega \cong \Ocal $
on $U_{x}\cup U_{y}$ and hence everywhere, by Lemma
\ref{l:hilbert/UxUy}.
\end{proof}


On $H_{n-1,n}$ we have two groups of twisting sheaves, $\Ocal
_{n-1}(k)$ and $\Ocal _{n}(l)$, pulled back from $H_{n-1}$ and $H_{n}$
respectively.  We abbreviate $\Ocal (k,l) = \Ocal _{n-1}(k)\otimes
\Ocal _{n}(l)$.

\begin{prop}\label{p:hilbert/omega-nested}
The canonical sheaf $\omega _{H_{n-1,n}}$ on the nested Hilbert scheme
$H_{n-1,n}$ is given by $\Ocal (1,-1)$ in the above notation.
\end{prop}

\begin{proof}
We have tautological sheaves $B_{n-1}$ and $B_{n}$ pulled back from
$H_{n-1}$ and $H_{n}$.  The kernel $L$ of the canonical surjection
$B_{n}\rightarrow B_{n-1}$ is the line bundle with fiber
$I_{n-1}/I_{n}$ at the point $(I_{n-1},I_{n})$.  From Proposition
\ref{p:hilbert/O(1)} we have $L = \Ocal (-1,1)$.  On the generic
locus, the fiber $I_{n-1}/I_{n}$ can be identified with the
one-dimensional space of functions on $V(I_{n})$ that vanish except at
$P_{n}$.  Thus the ratio of two sections of $L$ is determined by
evaluation at $x = x_{n}$, $y=y_{n}$.

Regarding the polynomials in \eqref{e:hilbert/x-generator} and
\eqref{e:hilbert/y-generator} as regular functions on $U_{x}\times \CC
^{2}$, they are the defining equations of the universal family $F_{x}
= \pi ^{-1}(U_{x})$ over $U_{x}\subseteq H_{n}$, as a closed subscheme
of the affine scheme $U_{x}\times \CC ^{2}$.  We can use these
defining equations to eliminate $e_{n}$ and $y$, showing that $F_{x}$
is an affine cell with coordinates $x, e_{1}, \ldots, e_{n-1}, a_{0},
\ldots, a_{n-1}$.


Over the curvilinear locus, the morphism $\alpha \colon
H_{n-1,n}\rightarrow F$ in \eqref{e:hilbert/nested->F} restricts to a
bijective morphism of smooth schemes, hence an isomorphism.  Under
this isomorphism $x$ corresponds to the $x$-coordinate $x_{n}$ of the
distinguished point, and modulo $x_{n}$ we can replace the elementary
symmetric functions $e_{k}(\xx )$ with $e'_{k} =
e_{k}(x_{1},\ldots,x_{n-1})$, for $k = 1,\ldots,n-1$.  As in the proof
of Proposition \ref{p:hilbert/omega-Hn}, we now calculate that a
nowhere vanishing regular section of $\omega $ on $U_{x}\subseteq
H_{n-1,n}$ is given by
\begin{equation}\label{e:hilbert/dxdeda}
t_{x} = dx_{n}\, d{\bf e}'\, d\aa = \frac{1}{\prod
_{i=1}^{n-1}(x_{n}-x_{i})}d\xx \, d\yy .
\end{equation}
By symmetry, $t_{y} = \left(1/\prod _{i=1}^{n-1}(y_{n}-y_{i}) \right)
d\xx \, d\yy $ is a nowhere vanishing regular section of $\omega $ on
$U_{y}$.

Now, at every point of $U_{x}$, the ideal $I_{n-1}$ is generated
modulo $I_{n}$ by
\begin{equation}\label{e:hilbert/x'-generator}
x^{n-1}-e'_{1}x^{n-2}+\cdots +(-1)^{n-1}e'_{n-1} = \prod
_{i=1}^{n-1}(x-x_{i}),
\end{equation}
so this expression represents a nowhere vanishing section $s_{x}$ of
$L$ on $U_{x}$.  Similarly, $\prod _{i=1}^{n-1}(y-y_{i})$ represents a
nowhere vanishing section $s_{y}$ of $L$ on $U_{y}$.  By the
observations in the first paragraph of the proof, the ratio
$s_{x}/s_{y}$ is the rational function $\prod
_{i=1}^{n-1}(x_{n}-x_{i})/\prod _{i=1}^{n-1}(y_{n}-y_{i})$ on
$H_{n-1,n}$.  Since we have nowhere vanishing sections $t_{x}$,
$t_{y}$ of $\omega $ on $U_{x}$ and $U_{y}$ with $t_{y}/t_{x} =
s_{x}/s_{y}$ it follows that we have $\omega \cong L^{-1} = \Ocal
(1,-1)$ on $U_{x}\cap U_{y}$ and hence everywhere, by Lemma
\ref{l:hilbert/UxUy}.
\end{proof}


\subsection{Geometry of $X_{n}$ and the $n!$ conjecture}
\label{ss:hilbert/n!}

Recall that the $n!$ conjecture, Conjecture \ref{conj:n!/n!}, concerns
the space $D_{\mu }$ spanned by all derivatives of the polynomial
$\Delta _{\mu }$ in \S \ref{ss:n!/n!}, eq.~\eqref{e:n!/Delta-mu}.  Let
\begin{equation}\label{e:hilbert/Jmu}
J_{\mu } = \{p\in \CC [\xx ,\yy ]:p(\partial \xx ,\partial \yy )\Delta
_{\mu } = 0 \}
\end{equation}
be the ideal of polynomials whose associated partial differential
operators annihilate $\Delta _{\mu }$, and set
\begin{equation}\label{e:hilbert/Rmu}
R_{\mu } = \CC [\xx ,\yy ]/J_{\mu }.
\end{equation}
Clearly $D_{\mu }$ and $R_{\mu }$ have the same dimension as vector
spaces.  The ideal $J_{\mu }$ is doubly homogeneous and
$S_{n}$-invariant, so $R_{\mu }$ is a doubly graded ring with an
action of $S_{n}$ which respects the grading.


We have the following useful characterization of the ideal $J_{\mu }$.

\begin{prop}[\cite{Hai99}, Proposition 3.3]
\label{p:hilbert/Jmu-characterization}%
The ideal $J_{\mu }$ in \eqref{e:hilbert/Jmu} is equal to the set of
polynomials $p\in \CC [\xx ,\yy ]$ such that the coefficient of
$\Delta _{\mu }$ in $\Theta ^{\epsilon }(gp)$ is zero for all $g\in
\CC [\xx ,\yy ]$.
\end{prop}

Note that it makes sense to speak of the coefficient of $\Delta _{\mu
}$ in an alternating polynomial, since the polynomials $\Delta _{D}$
form a basis of $ \CC [\xx ,\yy ]^{\epsilon }$.


It seems infeasible to describe explicitly the ideal sheaf of $X_{n}$
as a closed subscheme of $H_{n}\times (\CC ^{2})^{n}$, but we can give
an implicit description of this ideal sheaf by regarding $X_{n}$ as a
closed subscheme of $F^{n}$.  Here $F^{n}$ denotes the relative
product of $n$ copies of the universal family $F$, as a scheme over
$H_{n}$.  Like $X_{n}$, $F^{n}$ is a closed subscheme of $H_{n}\times
(\CC ^{2})^{n}$.  Its fiber over a point $I=I(S)$ is $S^{n}$, so the
closed points of $F^{n}$ are the tuples $(I,P_{1},\ldots,P_{n})$
satisfying $P_{i}\in V(I)$ for all $i$.  In particular $F^{n}$
contains $X_{n}$.  Being a closed subscheme of $F^{n}$, $X_{n}$ can be
defined as a scheme over $H_{n}$ by
\begin{equation}\label{e:hilbert/Xn-in-Fn}
X_{n} = \Spec B^{\otimes n}/\Jcal 
\end{equation}
for some sheaf of ideals $\Jcal $ in the sheaf of $\Ocal _{H_{n}}$-algebras
$B^{\otimes n}$.

\begin{prop}\label{p:hilbert/J-characterization}
Let
\begin{equation}\label{e:hilbert/Bn.Bn->wedgeB}
B^{\otimes n}\otimes B^{\otimes n}\rightarrow B^{\otimes n}\rightarrow
\ext ^{n}B
\end{equation}
be the map of $\Ocal _{H_{n}}$-module sheaves given by multiplication,
followed by the alternation operator $\Theta ^{\epsilon }$.  Then the
ideal sheaf $\Jcal $ of $X_{n}$ as a subscheme of $F^{n}$ is the
kernel of the map
\begin{equation}\label{e:hilbert/Bn->Bn wedgeB}
\phi \colon B^{\otimes n}\rightarrow (B^{\otimes n})^{*}\otimes \ext
^{n}B
\end{equation}
induced by \eqref{e:hilbert/Bn.Bn->wedgeB}.
\end{prop}


\begin{proof}
Let $U\subseteq H_{n}$ be the generic locus, the open set consisting
of ideals $I\in H_{n}$ for which $S = V(I)$ consists of $n$ distinct,
reduced points.  Note that $F^{n}$ is clearly reduced over $U$, and
since $F^{n}$ is flat over $H_{n}$ and $U$ is dense, $F^{n}$ is
reduced everywhere.  (For this argument we do not have to assume that
$F^{n}$ is irreducible, and indeed for $n>1$ it is not: $X_{n}$ is one
of its irreducible components.)  Sections of $B^{\otimes n}$ can be
identified with regular functions on suitable open subsets of $F^{n}$.
Since $X_{n}$ is reduced and irreducible (Proposition
\ref{p:hilbert/irreducible}), the open set $W = \rho ^{-1}(U)$ is
dense in $X_{n}$, and $\Jcal $ consists of those sections of
$B^{\otimes n}$ whose restrictions to $U$ define regular functions
vanishing on $W$.

Let $s$ be a section of $\Jcal $.  For any section $g$ of $B^{\otimes
n}$, the section $\Theta ^{\epsilon }(gs)$ also belongs to $\Jcal $.
Since it is alternating, $\Theta ^{\epsilon }(gs)$ vanishes at every
point $(I,P_{1},\ldots,P_{n})\in F^{n}$ for which two of the $P_{i}$
coincide.  In particular $\Theta ^{\epsilon }(gs)$ vanishes on
$F^{n}\setminus X_{n}$ and hence on $F^{n}$, that is, $\Theta
^{\epsilon }(gs) = 0$.  This is precisely the condition for $s$ to
belong to the kernel of $\phi $.

Conversely, if $s$ does not vanish on $W$ there is a point $Q =
(I,P_{1},\ldots,P_{n})$, with all $P_{i}$ distinct, where the regular
function represented by $s$ is non-zero.  Multiplying $s$ by a
suitable $g$, we can arrange that $gs$ vanishes at every point in the
$S_{n}$ orbit of $Q$, except for $Q$.  Then $\Theta ^{\epsilon
}(gs)\not =0$, so $s$ is not in the kernel of $\phi $.
\end{proof}


\begin{prop}\label{p:hilbert/equiv}
Let $Q_{\mu }$ be the unique point of $X_{n}$ lying over $I_{\mu }\in
H_{n}$.  The following are equivalent:
\begin{itemize}
\item [(1)] $X_{n}$ is locally Cohen-Macaulay and Gorenstein at
$Q_{\mu }$;
\item [(2)] the $n!$ conjecture holds for the partition $\mu $.
\end{itemize}
When these conditions hold, moreover, the ideal of the
scheme-theoretic fiber $\rho ^{-1}(I_{\mu })\subseteq X_{\mu }$
as a closed subscheme of $(\CC ^{2})^{n}$ coincides with the ideal
$J_{\mu }$ in \eqref{e:hilbert/Jmu}.
\end{prop}

\begin{proof}
Since $B^{\otimes n}$ and $(B^{\otimes n})^{*}\otimes \ext ^{n} B$ are
locally free sheaves, the sheaf homomorphism $\phi $ in
\eqref{e:hilbert/Bn->Bn wedgeB} can be identified with a linear
homomorphism of vector bundles over $H_{n}$.  Let
\begin{equation}\label{e:hilbert/phi(I)}
\phi (I)\colon B^{\otimes n}(I)\rightarrow B^{\otimes n}(I)^{*}\otimes
\ext ^{n} B(I)
\end{equation}
denote the induced map on the fiber at $I$.  The rank $\rk \phi (I)$
of the fiber map is a lower semicontinuous function, that is, the set
$\{I\colon \rk \phi (I)\geq r \}$ is open for all $r$.  If $\rk \phi
(I)$ is constant on an open set $U$, then the cokernel of $\phi $ is
locally free on $U$, and conversely.  When this holds, $\im \phi $ is
also locally free, of rank equal to the constant value of $\rk \phi
(I)$.  By Proposition \ref{p:hilbert/J-characterization}, we have $\im
\phi = \rho _{*}\Ocal _{X_{n}}$.  The fiber of $X_{n}$ over a point
$I$ in the generic locus consists of $n!$ reduced points, so the
generic rank of $\phi $ is $n!$.


The fiber $B(I)$ of the tautological bundle at $I$ is $R/I$.
Identifying $R^{\otimes n}$ with $\CC [\xx ,\yy ]$, we have a
linear map
\begin{equation}\label{e:hilbert/eta}
\eta \colon \CC [\xx ,\yy ]\rightarrow B^{\otimes n}(I)^{*}\otimes
\ext ^{n} B(I)
\end{equation}
given by composing $\phi (I)$ with the canonical map $\CC [\xx ,\yy ]
\rightarrow (R/I)^{\otimes n}$.  It follows from the definition of
$\phi $ that $\eta (p) = 0$ if and only if $\lambda \Theta ^{\epsilon
}(gp) = 0$ for all $g\in \CC [\xx ,\yy ]$, where
\begin{equation}\label{e:hilbert/lambda}
\lambda \colon A\rightarrow \ext ^{n} (R/I)\cong \CC 
\end{equation}
is the restriction of the canonical map $\CC [\xx ,\yy ] \rightarrow
(R/I)^{\otimes n}$ to $S_{n}$-alternating elements.  In particular,
for $I = I_{\mu }$, we have $\lambda (\Delta _{D}) = 0$ for all $D\not
=D(\mu )$, and $\lambda (\Delta _{\mu })$ spans $\ext ^{n}(R/I)$.
Hence, by Proposition \ref{p:hilbert/Jmu-characterization}, the kernel
of $\eta $ is exactly the ideal $J_{\mu }$ in this case.

Suppose the $n!$ conjecture holds for $\mu $.  Then since $\eta $ and
$\phi (I)$ have the same image, we have $\rk \phi (I_{\mu }) = n!$.
By lower semicontinuity, since $n!$ is also the generic rank of $\phi
$, $\rk \phi (I)$ is locally constant and equal to $n!$ on a
neighborhood of $I$.  Hence $\im \phi  = \rho _{*}\Ocal _{X_{n}}$ is
locally free, $\rho $ is flat of degree $n!$, and $X_{n}$ is locally
Cohen-Macaulay at $Q_{\mu }$.


Let $M$ be the the maximal ideal of the regular local ring $\Ocal
_{H_{n},I_{\mu }}$.  Since $X_{n}$ is finite over $H_{n}$, the ideal
$N = M\Ocal _{X_{n},Q_{\mu }}$ is a parameter ideal.  Assuming the
$n!$ conjecture holds for $\mu $, the Cohen-Macaulay ring $\Ocal
_{X_{n},Q_{\mu }}$ is Gorenstein if and only if $\Ocal _{X_{n},Q_{\mu
}}/N$ is Gorenstein.  We have $\Ocal _{X_{n},Q_{\mu }}/N\cong (\im
\phi )\otimes _{\Ocal _{H_{n}}}\Ocal _{H_{n},I_{\mu }}/M$, so the fiber
map $\phi (I_{\mu })$ factors as
\begin{equation}\label{e:hilbert/factored-phi}
\phi (I_{\mu })\colon B^{\otimes n}(I_{\mu }) \rightarrow \Ocal
_{X_{n},Q_{\mu }}/N \rightarrow B^{\otimes n}(I_{\mu })^{*} \otimes
\ext ^{n} B(I_{\mu }),
\end{equation}
and since $\coker \phi $ is locally free, the second homomorphism
above is injective.  Its image is $\im \phi (I_{\mu }) = \im \eta =
\CC [\xx ,\yy ]/J_{\mu }$, so we have $\Ocal _{X_{n},Q_{\mu }}/N \cong
\CC [\xx ,\yy ]/J_{\mu }$, and the latter ring is Gorenstein by
\cite{Ems78}, Proposition 4.


Conversely, suppose $X_{n}$ is locally Gorenstein at $Q_{\mu }$.  Then
$\Ocal _{X_{n},Q_{\mu }}/N$ is a Gorenstein Artin local ring
isomorphic to $\CC [\xx ,\yy ]/J$ for some ideal $J$.  Since $\rho
_{*}\Ocal _{X_{n}}$ is locally free on a neighborhood of $I_{\mu }$,
necessarily of rank $n!$, we have $\dim _{\CC }\CC [\xx ,\yy ]/J =
n!$.  The locally free sheaf $\rho _{*}\Ocal _{X_{n}}$ is the sheaf of
sections of a vector bundle, which is actually a bundle of $S_{n}$
modules, since $S_{n}$ acts on $X_{n}$ as scheme over $H_{n}$.  The
isotypic components of such a bundle are direct summands of it and
hence locally free themselves, so the character of $S_{n}$ on the
fibers in constant.  In our case, the generic fibers are the
coordinate rings of the $S_{n}$ orbits of points
$(P_{1},\ldots,P_{n})\in (\CC ^{2})^{n}$ with all $P_{i}$ distinct.
Therefore every fiber affords the regular representation of $S_{n}$.

The socle of $\CC [\xx ,\yy ]/J$ is a one-dimensional
$S_{n}$-invariant subspace.
Since $\CC [\xx ,\yy ]/J$ affords the regular representation its only
such subspaces are $(\CC [\xx ,\yy ]/J)^{S_{n}}$, which consists of
the constants, and $(\CC [\xx ,\yy ]/J)^{\epsilon }$.  The
socle must therefore be the latter space (unless $n=1$, in which case
the Proposition is trivial).

The factorization of $\phi (I_{\mu })$ in
\eqref{e:hilbert/factored-phi} implies that $J\subseteq J_{\mu } =
\ker \eta $.  If $J_{\mu }/J\not =0$ then we must have $\soc (\CC [\xx
,\yy ]/J)\subseteq J_{\mu }/J$, as the socle is contained in every
non-zero ideal.  But this would imply $(C[\xx ,\yy ]/J_{\mu
})^{\epsilon } = 0$ and hence $\Delta _{\mu }\in J_{\mu }$, which is
absurd.
\end{proof}


\begin{prop}\label{p:hilbert/n!->main}
If the $n!$ conjecture holds for all partitions $\mu $ of a given $n$,
then $X_{n}$ is Gorenstein with canonical line bundle $\omega _{X_{n}}
= \Ocal (-1)$.
\end{prop}

\begin{proof}
The set $U$ of points $I\in H_{n}$ such that $\rk \phi (I) = n!$ is
open and $\TT ^{2}$-invariant.  From the proof of Proposition
\ref{p:hilbert/equiv} we see that the $n!$ conjecture implies that $U$
contains all the monomial ideals $I_{\mu }$.  Since every $I\in H_{n}$
has a monomial ideal in the closure of its orbit, this implies $U =
H_{n}$, so $\rho \colon X_{n}\rightarrow H_{n}$ is flat, and $X_{n}$
is Cohen-Macaulay.

Let $P = \rho _{*}\Ocal _{X_{n}}$, a locally free sheaf of rank $n!$.
By Proposition \ref{p:hilbert/J-characterization}, the map in
\eqref{e:hilbert/Bn.Bn->wedgeB} induces a pairing
\begin{equation}\label{e:hilbert/P P->wedgeB}
P\otimes P\rightarrow \ext ^{n}B = \Ocal (1)
\end{equation}
and $\phi $ factors through the induced homomorphism
\begin{equation}\label{e:hilbert/P->P* O(1)}
\tilde{\phi }\colon P\rightarrow P^{*}\otimes \Ocal (1)\cong \Hom
(P,\Ocal (1)). 
\end{equation}
Note that the multiplication by a section $s$ of $P$ in $\Hom (P,\Ocal
(1))$ as a sheaf of $P$-modules is given by $(s\lambda )(h) = \lambda
(sh)$.  By the definition of $\phi $ we have $\tilde{\phi }(sg)(h) =
\Theta ^{\epsilon }(sgh) = \tilde{\phi }(g)(sh)$, so $\tilde{\phi }$
is a homomorphism of sheaves of $P$-modules.  Since $\rk \phi (I)$ is
constant and equal to $n!$, which is the rank of both $P$ and
$P^{*}\otimes \Ocal (1)$, $\tilde{\phi }$ is an isomorphism.

Now, $X_{n} = \Spec P$, so by the duality theorem, $\omega _{X_{n}}$
is the sheaf of $\Ocal _{X_{n}}$-modules associated to the sheaf of
$P$-modules $\omega _{H_{n}}\otimes P^{*}$.  By Proposition
\ref{p:hilbert/omega-Hn} we have $\omega _{H_{n}} \cong \Ocal
_{H_{n}}$, and we have just shown $P^{*}\cong P\otimes \Ocal (-1)$.
Together these imply $\omega _{X_{n}} = \Ocal_{X_{n}} (-1)$.
\end{proof}


\subsection{Main theorem}
\label{ss:hilbert/main}

\begin{thm}\label{t:main}
The isospectral Hilbert scheme $X_{n}$ is normal, Cohen-Macaulay, and
Gorenstein, with canonical sheaf $\omega _{X_{n}} \cong  \Ocal (-1)$.
\end{thm}

The proof of this theorem will occupy us for the rest of this
subsection.  In principle, to show that $X_{n}$ is Cohen-Macaulay (at
a point $Q$, say), we would like to exhibit a local regular sequence
of length $2n = \dim X_{n}$.  In practice, we are unable to do this,
but we can show that the $\yy $ coordinates form a regular sequence of
length $n$ wherever they vanish.  As it turns out, showing this much
is half of the battle.  A geometric induction argument takes care of
the other half.
 
The key geometric property of $X_{n}$, which implies that the $\yy $
coordinates form a regular sequence, is given by the following pair of
results, the first of which is proven in \S \ref{ss:poly/proof}.

\begin{prop}\label{p:hilbert/J-is-free}
Let $J = \CC [\xx ,\yy ] A$ be the ideal generated by the space of
alternating polynomials $A = \CC [\xx ,\yy ]^{\epsilon }$.  Then
$J^{d}$ is a free $\CC [\yy ]$-module for all $d$.
\end{prop}

\begin{cor}\label{c:hilbert/X-flat-over-y}
The projection $X_{n}\rightarrow \CC ^{n} = \Spec \CC [\yy ]$ of
$X_{n}$ on the $\yy $ coordinates is flat.
\end{cor}

\begin{proof}
Let $S = \CC [\xx ,\yy ]$.  By Proposition \ref{p:hilbert/blowup}, we
have $X_{n} = \Proj S[tJ]$, and Proposition \ref{p:hilbert/J-is-free}
implies that $S[tJ]$ is a free $\CC [\yy ]$-module.
\end{proof}


An alternating polynomial $g\in A$ must vanish at every point
$(P_{1},\ldots,P_{n})\in (\CC ^{2})^{n}$ where two of the $P_{i}$
coincide.  Hence we have
\begin{equation}\label{e:hilbert/J-vs-Jij}
J\subseteq \bigcap _{i<j} (x_{i}-x_{j},y_{i}-y_{j}),
\end{equation}
and it is natural to conjecture by analogy to the univariate case that
equality holds here.  In \cite{Hai99}, Proposition 6.2 we proved that
the following more general identity holds once we know that $J^{d}$ is
a free $\CC [\yy ]$-module for all $n$ and $d$.

\begin{cor}\label{c:hilbert/free->what-J-is}
We have
\begin{equation}\label{e:hilbert/Jd-vs-Jijd}
J^{d} = \bigcap _{i<j} (x_{i}-x_{j},y_{i}-y_{j})^{d}
\end{equation}
for all $n$ and $d$.
\end{cor}

Remarkably, even though this seems like it should be an elementary
result, we know of no proof not using Proposition \ref{p:hilbert/J-is-free},
even for $d=1$.  Using Corollary \ref{c:hilbert/free->what-J-is}, we
can deal with the normality question.

\begin{prop}\label{e:hilbert/normal}
The isospectral Hilbert scheme $X_{n}$ is arithmetically normal in its
projective embedding over $(\CC ^{2})^{n}$ as the blowup $X_{n} =
\Proj S[tJ]$.  In particular $X_{n}$ is normal.
\end{prop}

\begin{proof}
By definition {\it arithmetically normal} means that $S[tJ]$ is a
normal domain.  Since $S$ itself is a normal domain, this is
equivalent to the ideals $J^{d}\subseteq S$ being integrally closed
ideals for all $d$.  The powers of an ideal generated by a regular
sequence are integrally closed, as is an intersection of integrally
closed ideals, so $J^{d}$ is integrally closed by Corollary
\ref{c:hilbert/free->what-J-is}.
\end{proof}


For the Cohen-Macaulay and Gorenstein properties of $X_{n}$ we use an
inductive argument involving the nested Hilbert scheme, duality, and
the following lemma.

\begin{lem}\label{l:hilbert/Rg*}
Let $g\colon Y\rightarrow X$ be a proper morphism.  Let
$z_{1},\ldots,z_{m}\in \Ocal _{X}(X)$ be global regular functions on
$X$ (and, via $g$, on $Y$).  Let $Z\subseteq X$ be the closed subset
$Z = V(z_{1},\ldots,z_{m})$ and let $U = X\setminus Z$ be its
complement.  Suppose the following conditions hold.
\begin{itemize}
\item [(1)] The $z_{i}$ form a regular sequence in the local ring
$\Ocal _{X,P}$ for all $P\in Z$.
\item [(2)] The $z_{i}$ form a regular sequence in the local ring
$\Ocal _{Y,Q}$ for all $Q\in g^{-1}(Z)$.
\item [(3)] Every fiber of $g$ has dimension less than $m-1$.
\item [(4)] The canonical homomorphism $\Ocal _{X}\rightarrow
Rg_{*}\Ocal _{Y}$ restricts to an isomorphism on $U$.
\end{itemize}
Then $Rg_{*}\Ocal _{Y} = \Ocal _{X}$, {\it i.e.}, the canonical
homomorphism is an isomorphism.
\end{lem}


\begin{proof}
The question is local on $X$, so without loss of generality we can
assume $X = \Spec S$ is affine.  Then we are to show that
$H^{i}(Y,\Ocal _{Y}) = 0$ for $i>0$ and that $S\rightarrow
H^{0}(Y,\Ocal _{Y})$ is an isomorphism.

Condition (3) implies that $H^{i}(Y,\Ocal _{Y}) = 0$ for $i\geq m-1$.
Let $Z' = g^{-1}(Z)$ and $U' = g^{-1}(U)$.  Conditions (1) and (2)
imply that $\depth _{Z}\Ocal _{X}$ and $\depth _{Z'}\Ocal _{Y}$ are
both at least $m$.  Hence the local cohomology modules
$H^{i}_{Z}(\Ocal _{X})$ and $H^{i}_{Z'}(\Ocal _{Y})$ vanish for $i\leq
m-1$.  By the exact sequence of local cohomology \cite{Har67} we
therefore have
\begin{equation}\label{e:hilbert/local-Y}
H^{i}(Y,\Ocal _{Y})\cong H^{i}(U',\Ocal _{Y})
\end{equation}  
and
\begin{equation}\label{e:hilbert/local-X}
H^{i}(X,\Ocal _{X})\cong H^{i}(U,\Ocal _{X})
\end{equation}
for all $i<m-1$.  Condition (4) yields  $H^{i}(U',\Ocal _{Y})\cong
H^{i}(U,\Ocal _{X})$.  Thus we have $H^{i}(Y,\Ocal _{Y})\cong
H^{i}(X,\Ocal _{X})$ for all $i<m-1$.

Since $X$ is affine, this shows $H^{i}(Y,\Ocal _{Y}) = 0$ for all
$i>0$, and $S\cong H^{0}(Y,\Ocal _{Y})$.  The isomorphism $S\cong
H^{0}(Y,\Ocal _{Y})$ is the canonical homomorphism, since it is
determined by its restriction to $U$.
\end{proof}


We now prove Theorem~\ref{t:main} by induction on $n$.  In the proof
of the inductive step we will assume for technical reasons that $n >
3$.  Thus we first need to dispose of the cases $n=1,2,3$.  By
Proposition \ref{p:hilbert/n!->main}, this reduces to verifying the
$n!$ conjecture for $n\leq 3$, which is a simple matter of
calculation.

Note that for $n=1$, we have $X_{1} = \CC ^{2}$ trivially, and for
$n=2$, $X_{2}$ is the blowup of $(\CC ^{2})^{2}$ along the diagonal.
Thus $X_{1}$ and $X_{2}$ are actually nonsingular.  To this we may add
that, up to automorphisms induced by translations in $\CC ^{2}$, the
scheme $X_{3}$ has an essentially isolated singularity at the point
$Q_{(2,1)}$ lying over $I_{(2,1)}\in H_{n}$.  Observe that $I_{(2,1)}
= (x,y)^{2}$ is the ideal of a ``fat point'' of length $3$ at the
origin in $\CC ^{2}$, and that such fat points are the smallest
non-curvilinear subschemes of the plane.  It can be shown in general
that the singular locus in $X_{n}$ coincides with the non-curvilinear
locus.

Assume now by induction that $X_{n-1}$ is Cohen-Macaulay and
Gorenstein with $\omega _{X_{n-1}} = \Ocal _{X_{n-1}}(-1)$.  Then
$\rho _{n-1}\colon X_{n-1}\rightarrow H_{n-1}$ is flat, so in the
scheme-theoretic fiber square
\begin{equation}\label{e:hilbert/inductive-square}
\begin{CD}
Y&	@>{\rho '}>>&	H_{n-1,n}\\
@VVV&	&		@VVV\\
X_{n-1}&@>{\rho }>>&	H_{n-1},
\end{CD}
\end{equation}
the morphism $\rho '$ is also flat.  On the generic locus, where
$P_{1},\ldots,P_{n}$ are all distinct, the above diagram coincides
locally with the fiber square
\begin{equation}\label{e:hilbert/generic-square}
\begin{CD}
Y&	@>{\rho '}>>&	S^{n-1}\CC ^{2}\times \CC ^{2}\\
@VVV&	&		@VVV\\
(\CC ^{2})^{n-1}&@>{\rho }>>&	S^{n-1}\CC ^{2}.
\end{CD}
\end{equation}
This shows that $Y$ is generically reduced, hence reduced, as well as
irreducible and birational to $(\CC ^{2})^{n}$.  The reduced fiber
product in \eqref{e:hilbert/inductive-square} is $X_{n-1,n}$ by
definition, so we have $Y = X_{n-1,n}$.  Since $\rho '$ is flat and
finite, and $H_{n-1,n}$ is nonsingular, $X_{n-1,n}$ is
Cohen-Macaulay.


Furthermore, the relative canonical sheaf of $X_{n-1,n}$ over
$H_{n-1,n}$ is the pullback of that of $X_{n-1}$ over $H_{n-1}$.  By
Proposition \ref{p:hilbert/omega-Hn} and the induction hypothesis, the
latter is $\Ocal _{X_{n-1}}(-1)$ and its pullback to $X_{n-1,n}$ is
$\Ocal (-1,0)$.  By Proposition \ref{p:hilbert/omega-nested} it
follows that the canonical sheaf on $X_{n-1,n}$ is $\omega
_{X_{n-1,n}} = \Ocal (-1,0)\otimes \Ocal (1,-1) = \Ocal (0,-1)$.  In
particular, $X_{n-1,n}$ is Gorenstein.

Now consider the projection
\begin{equation}\label{e:hilbert/g}
g\colon Y = X_{n-1,n}\rightarrow X_{n}.
\end{equation}
We claim that $Rg_{*}\Ocal _{Y} = \Ocal _{X_{n}}$.  By the projection
formula, since $\Ocal (0,-1) = g^{*}\Ocal _{X_{n}}(-1)$ is pulled back
from $X_{n}$, this implies also $Rg_{*}\Ocal (0,-1) = \Ocal
_{X_{n}}(-1)$.  Now $\Ocal (0,-1)[2n] = \omega _{X_{n-1,n}}[2n]$ is
the dualizing complex on $X_{n-1,n}$, so by the duality theorem it
follows that $\Ocal (-1)[2n]$ is the dualizing complex on $X_{n}$.  In
other words, $X_{n}$ is Gorenstein, with canonical sheaf $\omega
_{X_{n}} = \Ocal (-1)$, which is what we wanted to prove.

For the claim, we verify the conditions of Lemma \ref{l:hilbert/Rg*}
with $Y = X_{n-1,n}$, $X = X_{n}$, and $z_{1},\ldots,z_{n-1}$ equal to
$y_{1}-y_{2},\ldots,y_{n-1}-y_{n}$.  For the $z_{i}$ to form a regular
sequence it suffices that $y_{1},\ldots,y_{n}$ is a regular
sequence.  On $Y = X_{n-1,n}$ this follows from the Cohen-Macaulay
property of $X_{n-1,n}$, together with Proposition
\ref{p:hilbert/V(y)-in-nest-X}, which says that $V(\yy )$ is a
complete intersection in $Y$.  On $X = X_{n}$ (this is the crucial
step!), the regular sequence condition follows from Corollary
\ref{c:hilbert/X-flat-over-y}.


On the open set $U = X_{n}\setminus V(\zz )$, the coordinates $y_{i}$
are nowhere all equal, by definition.  It follows from Lemmas
\ref{l:hilbert/product} and \ref{l:hilbert/product-nested} that $U$
can be covered by open sets on which the projection $g\colon
X_{n-1,n}\rightarrow X_{n}$ is locally isomorphic to $1\times
g_{l}\colon X_{k}\times X_{l-1,l}\rightarrow X_{k}\times X_{l}$ for
some $k+l = n$ with $l<n$.  We can assume as part of the induction
that $R(g_{l})_{*}\Ocal _{X_{l-1,l}} = \Ocal _{X_{l}}$.  Hence we
have $Rg_{*}\Ocal _{X_{n-1,n}}|_{g^{-1}(U)} = \Ocal _{X_{n}}|_{U}$.

For the fiber dimension condition, note that the fiber of $g$ over a
point $(I,P_{1},\ldots,P_{n})$ of $X_{n}$ is the same as the fiber of
the morphism $\alpha $ in \eqref{e:hilbert/nested->F} over
$(I,P_{n})$.  By Proposition \ref{p:hilbert/fibers}, its dimension $d$
satisfies the inequality $\binom{d+2}{2}\leq n$.  Since we are assuming
$n>3$, this inequality easily implies $d<n-2$, as required for Lemma
\ref{l:hilbert/Rg*} to hold.


\subsection{Proof of the graded character conjecture}
\label{ss:hilbert/character}

We now review in outline the proof from \cite{Hai99} that the $n!$
conjecture and the Cohen-Macaulay property of $X_{n}$ imply Conjecture
\ref{conj:n!/GH}.  The proof involves some technical manipulations
with Frobenius series and Macdonald polynomials which it would take us
too far afield to repeat in full here.  Conceptually, however, the
argument is straightforward.  The main point is the connection between
$S_{n}$ characters and symmetric functions given by the {\it Frobenius
characteristic}
\begin{equation}\label{e:hilbert/frobenius}
\Psi (\chi ) = \frac{1}{n!} \sum _{w\in S_{n}} \chi (w) p_{\tau (w)}(x).
\end{equation}
Here $p_{\tau }(x)$ denotes a power-sum symmetric function and $\tau
(w)$ is the partition of $n$ given by the cycle lengths in the
expression for $w$ as a product of disjoint cycles.  The Frobenius
characteristics of the irreducible characters are the Schur functions
(\cite{Mac95}, Chapter I, eq.~(7.5))
\begin{equation}\label{e:hilbert/frobenius-chi-lamda}
\Psi (\chi ^{\lambda }) = s_{\lambda }(x).
\end{equation}

Let $J_{\mu }$ be be the ideal of operators annihilating $\Delta _{\mu
}$, and set $R_{\mu } = \CC [\xx ,\yy ]/J_{\mu }$, as in
\eqref{e:hilbert/Jmu}--\eqref{e:hilbert/Rmu}.  It follows from
\cite{Hai99}, Proposition 3.4, that $\ch (R_{\mu })_{r,s} = \ch
(D_{\mu })_{r,s}$ for all degrees $r,s$, so we may replace $D_{\mu }$
with $R_{\mu }$ in the statement of Conjecture \ref{conj:n!/GH}.

We define the {\it Frobenius series} of $R_{\mu }$ to be
\begin{equation}\label{e:hilbert/frob-series}
\Fcal _{R_{\mu }}(x;q,t) = \sum _{r,s} t^{r}q^{s}\Psi (\ch (R_{\mu
})_{r,s}).
\end{equation}
The Frobenius series is a kind of doubly graded Hilbert series that
keeps track of characters instead of just dimensions.  In this
notation Conjecture \ref{conj:n!/GH} takes the form of an identity
\begin{equation}\label{e:hilbert/F=H}
\Fcal _{R_{\mu }}(x;q,t) = \tilde{H}_{\mu }(x;q,t).
\end{equation}


There is a well-defined formal extension of the notion of Frobenius
series to certain local rings with $S_{n}$ and $\TT ^{2}$ actions, as
explained in \cite{Hai99}, Section 5.  In particular, we can define
$\Fcal _{S_{\mu }}(x;q,t)$, where $S_{\mu }=\Ocal _{X_{n},Q_{\mu }}$
is the local ring of $X_{n}$ at the distinguished point $Q_{\mu }$
lying over $I_{\mu }$.

By Theorem~\ref{t:main} and Proposition \ref{p:hilbert/equiv}, the
ring $R_{\mu }$ is the coordinate ring of the scheme-theoretic fiber
$\rho ^{-1}(I_{\mu })\subseteq X_{n}$.  This and the flatness of $\rho
$ imply (\cite{Hai99}, eq.~(5.3)) that the Frobenius series of $S_{\mu
}$ and $R_{\mu }$ are related by a scalar factor,
\begin{equation}\label{e:hilbert/Smu-vs-Rmu}
\Fcal _{S_{\mu }}(x;q,t) = \Hcal (q,t)\Fcal _{R_{\mu }}(x;q,t),
\end{equation}
where $\Hcal (q,t)$ is the formal Hilbert series of the local ring
$\Ocal _{H_{n},I_{\mu }}$, as defined in \cite{Hai99}, a rational
function of $q$ and $t$.  In fact, it follows from the determination
in \cite{Hai98} of explicit regular local parameters on $H_{n}$ at
$I_{\mu }$ that $\Hcal (q,t)$ is the reciprocal of the
polynomial 
\begin{equation}\label{e:hilbert/Hmu(q,t)-inverse}
\prod_{x\in D(\mu )} (1-q^{-a(x)}t^{1+l(x)}) \;\cdot 
\prod_{x\in D(\mu )} (1-q^{1+a(x)}t^{-l(x)}),
\end{equation}
where $a(x)$, $l(x)$ denote the lengths of the {\it arm} and {\it leg}
of the cell $x$ in the diagram of $\mu $.  It is interesting to note
that the two products above are essentially the normalizing factors
introduced in \cite{Mac95}, Chapter VI, eqs.~(8.1--8.1$'$) to define
the integral form Macdonald polynomials $J_{\mu }$.  This coincidence
is typical of the numerological parallels that already pointed to a
link between Macdonald polynomials and the Hilbert scheme, before the
theory presented here was fully developed.


Now, the global regular functions $y_{1},\ldots,y_{n}$ form a regular
sequence in the local ring $S_{\mu }$, as follows from the proof of
Theorem~\ref{t:main}.  By \cite{Hai99}, Proposition 5.3, part (3),
this implies that
\begin{equation}\label{e:hilbert/FSmu/y}
\Fcal _{S_{\mu }/(\yy )}(x;q,t) = \Fcal _{S_{\mu }}[X(1-q);q,t].
\end{equation}
Furthermore, \cite{Hai99}, Proposition 5.3, part (1) implies that if
the irreducible character $\chi ^{\lambda }$ has multiplicity zero in
$\ch R_{\mu }/(\yy )$, then the coefficient of the Schur function
$s_{\lambda }$ in $\Fcal _{S_{\mu }/(\yy )}(x;q,t)$ is equal to zero.


The ring $R_{\mu }/(\yy )$ is of course nothing other than the
$y$-degree zero part of $R_{\mu }$, whose Frobenius series is given by
the $q=0$ specialization $\Fcal _{R_{\mu }}(x;0,t)$.  Thus $R_{\mu
}/(\yy )$ coincides with a well-understood algebraic object: by the
results of \cite{BeGa92,GaPr92}, it is isomorphic to the cohomology
ring of the Springer fiber over a unipotent element of $GL(n)$, as
defined in \cite{Spr78}, and its Frobenius series is the transformed
Hall-Littlewood polynomial
\begin{equation}\label{e:hilbert/FRmu(q=0)}
\Fcal _{R_{\mu }}(x;0,t) = t^{n(\mu )}Q_{\mu }[X/(1-t^{-1});t^{-1}],
\end{equation} 
where $Q_\mu (x;t)$ is defined as in \cite{Mac95}, Chapter III,
eq.~(2.11).  We remark that \eqref{e:hilbert/FRmu(q=0)} is just the
$q=0$ specialization of \eqref{e:hilbert/F=H}.

In particular, the only characters $\chi ^{\lambda }$ that occur with
non-zero multiplicity in $\ch R_{\mu }/(\yy )$ are those with $\lambda
\geq \mu $.  Hence by \eqref{e:hilbert/Smu-vs-Rmu} and
\eqref{e:hilbert/FSmu/y} we have
\begin{equation}\label{e:hilbert/triangularity-1}
\Fcal _{R_{\mu }}[X(1-q);q,t]\in \QQ (q,t)\{s_{\lambda }(x)\colon
\lambda \geq \mu \}.
\end{equation}
This shows that $\Fcal _{R_{\mu }}(x;q,t)$ satisfies condition (1) of
the characterization of $\tilde{H}_{\mu }(x;q,t)$ in Proposition
\ref{p:n!/H-tilde-triang}.  By the symmetry between $\xx $ and $\yy $,
$\Fcal _{R_{\mu }}$ also satisfies condition (2).  Finally, condition
(3) for $\Fcal _{R_{\mu }}$ says that the trivial character occurs
only in degree zero, with multiplicity one.  This is obvious, since
the only $S_{n}$-invariants in $R_{\mu }$ are the constants.  Hence
\eqref{e:hilbert/F=H} holds.

To summarize, we have the following consequence of
Theorem~\ref{t:main} and Proposition \ref{p:hilbert/equiv}.

\begin{thm}\label{t:MPK}
Conjectures \ref{conj:n!/MPK}, \ref{conj:n!/n!}, and \ref{conj:n!/GH},
that is, the Macdonald positivity conjecture, the $n!$ conjecture, and
the graded character interpretation of $\tilde{K}_{\lambda \mu
}(q,t)$, are all true.
\end{thm}


\subsection{Appendix: Cohen-Macaulay and Gorenstein schemes and
duality theory}
\label{ss:hilbert/duality}

We review here facts from duality theory used in this section.  The
formulation below is valid for quasiprojective schemes over $\CC $,
which is sufficiently general for our purposes.  Proofs can be found
in \cite{Har66}.

Let $D(X)$ denote the derived category of complexes of sheaves of
$\Ocal _{X}$-modules with bounded, coherent cohomology.  If $A$ is a
coherent sheaf, we denote by $A[d]$ the complex which is $A$ in degree
$-d$ and zero in every other degree.  If $A^{\bullet }$ is any complex
in $D(X)$ whose only non-zero homology sheaf is $H^{-d}(A^{\bullet })
\cong A$, then we have $A^{\bullet }\cong A[d]$.  We say that such a
complex $A^{\bullet }$ {\it reduces to a sheaf}.

If a complex $\omega ^{\bullet } \in D(X)$ has finite injective
dimension, then $\Dcal _{\omega ^{\bullet }} = R\Hom (- ,\omega
^{\bullet })$ is a functor from $D(X)$ into itself.  A {\it dualizing
complex} on $X$ is a complex $\omega ^{\bullet }$ of finite injective
dimension such that the canonical natural transformation
$1_{D(X)}\rightarrow \Dcal_{\omega ^{\bullet}} \circ \Dcal_{\omega
^{\bullet}} $ is an isomorphism.

There exists a dualizing complex $\omega ^{\bullet }_{X}$ and
corresponding dualizing functor $\Dcal _{X} = \Dcal _{\omega ^{\bullet
}_{X}}$ on every quasiprojective scheme $X$ over $\CC $, with the
following properties.

(1) {\it Duality theorem:} if $f\colon Y\rightarrow X$ is a proper
morphism, then there is a canonical natural isomorphism $\Dcal
_{X}\circ Rf_{*}\cong Rf_{*}\circ \Dcal _{Y}$.

(2) If $f\colon Y\rightarrow X$ is smooth of relative dimension $d$
then $\omega ^{\bullet }_{Y} = \omega _{Y/X}[d]\otimes f^{*}\omega
^{\bullet }_{X}$, where $\omega _{Y/X}$ is the relative canonical line
bundle, that is, the sheaf of relative exterior $d$-forms $\Omega^{d}
_{Y/X}$.


A scheme $X$ is {\it Cohen-Macaulay} (respectively, {\it Gorenstein})
if its local ring $\Ocal _{X,P}$ at every point is a Cohen-Macaulay
(or Gorenstein) local ring.  If $\rho \colon X\rightarrow H$ is a
finite morphism of equidimensional schemes of the same dimension, with
$H$ smooth, then $X$ is Cohen-Macaulay if and only if $\rho $ is flat.

For $X$ quasiprojective over $\CC $ it follows from duality theory
that $X$ is Cohen-Macaulay if and only if the dualizing complex
$\omega ^{\bullet }_{X}$ reduces to a sheaf on each connected
component of $X$.  If $X$ is Cohen-Macaulay and equidimensional of
dimension $d$, then the dualizing complex is $\omega _{X}[d]$, where
$\omega _{X}$ is the {\it canonical sheaf} on $X$, whose stalk at each
point $P$ is the canonical module of $\Ocal _{X,P}$.  In particular,
$X$ is Gorenstein if and only if $\omega ^{\bullet }_{X}$ reduces to a
line bundle ({\it i.e.}, a locally free sheaf of rank $1$) on each
connected component of $X$.


\section{Polygraphs}
\label{s:polygraphs}

\subsection{First definitions}
\label{ss:poly/first-defs}

Let $E = \AA ^{2}(k)$ be the affine plane over a field $k$ of
characteristic zero (the restriction on $k$ is not really
necessary---see \S \ref{ss:poly/ground-ring}).  We are going to study
certain unions of linear subspaces, or {\it subspace arrangements}, in
$E^{n}\times E^{l}$.  We call these arrangements {\it polygraphs}
because their constituent subspaces are the graphs of linear maps from
$E^{n}$ to $E^{l}$.

Let $[n]$ denote the set of integers $\{1,\ldots,n \}$.  Given a
function $f\colon [l]\rightarrow [n]$, there is a linear morphism
\begin{equation}\label{e:poly/pi_f}
\pi _{f}\colon E^{n}\rightarrow E^{l}
\end{equation}
defined by
\begin{equation}\label{e:poly/pi_f(P)}
\pi _{f}(P_{1},\ldots,P_{n}) = (P_{f(1)},\ldots,P_{f(l)}).
\end{equation}
Let 
\begin{equation}\label{e:poly/W_f}
W_{f}\subseteq E^{n}\times E^{l}
\end{equation}
be the graph of $\pi _{f}$.  We denote the coordinates on $E^{n}\times
E^{l}$ by
\begin{equation}\label{e:poly/xyab}
\xx ,\yy ,\aa ,\bb
=x_{1},y_{1},\ldots,x_{n},y_{n},\; a_{1},b_{1},\ldots,a_{l},b_{l},
\end{equation}
where $x_{j},y_{j}$ are the coordinates on the $j$-th factor in
$E^{n}$ and $a_{i},b_{i}$ are the coordinates on the $i$-th factor in
$E^{l}$.  In coordinates, $W_{f}$ is then defined by the equations
\begin{equation}\label{e:poly/I_f}
W_{f} = V(I_{f}),\quad \text{where}\quad I_{f} = \sum _{i\in [l]}
(a_{i}-x_{f(i)},b_{i}-y_{f(i)}).
\end{equation}


\begin{defn}\label{d:poly/Z(n,l)}
The {\it polygraph} $Z(n,l)\subseteq E^{n}\times E^{l}$ is the subspace
arrangement
\begin{equation}\label{e:poly/def-of-Z(n,l)}
Z(n,l)=\bigcup W_{f},\quad \text{over all $f\colon [l]\rightarrow [n]$}.
\end{equation}
\end{defn}

The geometric points of $Z(n,l)$ are the points $(\xx ,\yy ,\aa ,\bb)
= (P_{1},\ldots,P_{n},Q_{1},\ldots,Q_{l})$ satisfying the following
condition: for all $i\in [l]$ there is a $j\in [n]$ such that
$Q_{i}=(a_{i},b_{i})$ is equal to $P_{j}=(x_{j},y_{j})$.  Note that
for $l=0$, the polygraph $Z(n,0)$ makes sense and is equal to $E^{n}$.
The index set $[l] = [0]$ is empty in this case, and the unique
function $f\colon [0]\rightarrow [n]$ is the empty function $f =
\emptyset $, with $I_{\emptyset } = 0$ and $W_{\emptyset } = E^{n}$.

A word on terminology is in order here.  For clarity it is is often
useful to describe a scheme or a morphism in terms of geometric
points, that is, points defined by values of the coordinates in some
algebraically closed extension $K$ of $k$.  In the preceding paragraph
and in \eqref{e:poly/pi_f(P)}, we have given geometric descriptions.
In all other contexts, however, we will use the term {\it point} in
the scheme-theoretic sense.  When discussing the local geometry of
$Z(n,l)$ at a point $P$, for instance, we mean that $P\in \Spec
R(n,l)$ is a prime ideal in the coordinate ring $R(n,l)$ of $Z(n,l)$.
Since all schemes under consideration are closed subschemes of
$E^{n}\times E^{l}$, we can also identify $P$ with a prime ideal of
$k[\xx ,\yy ,\aa ,\bb ]$.  Ultimately, everything we do reduces to
commutative algebra and ideal theory in the polynomial ring $k[\xx
,\yy ,\aa ,\bb ]$.  Any geometric descriptions we give are best
understood merely as guides to formal definitions in terms of ideals
and ring homomorphisms.  Thus the polygraph $Z(n,l)$ is correctly
defined as the subscheme of $E^{n}\times E^{l}$ whose ideal is the
intersection of the ideals $I_{f}$ defined in \eqref{e:poly/I_f}.

Our purpose in this section is to prove the following theorem.

\begin{thm}\label{t:polygraphs-are-free}
The coordinate ring $R(n,l) =\Ocal (Z(n,l))$ of the polygraph
$Z(n,l)$ is a free $k[\yy ]$-module.
\end{thm}


\subsection{Examples}
\label{ss:poly/examples}

As motivation for Theorem~\ref{t:polygraphs-are-free}, it may be helpful to
consider the case of polygraphs in one set of variables, that is, the
subspace arrangements defined as in \ref{d:poly/Z(n,l)}, but with
$E=\AA ^{1}$.

Polygraphs as we have defined them are in {\it two sets of variables},
namely the $\xx ,\aa $ and the $\yy ,\bb $.  Had we begun with
$E=\AA ^{d}$ instead of $E=\AA ^{2}$, we would have $d$ sets of
variables (and notational headaches galore).  In one set of variables
the coordinates are just $\xx ,\aa $, and the ideal of $Z(n,l)$ is
simply
\begin{equation}\label{e:poly/univariate-I}
I = \sum _{i\in [l]}\left( \vphantom{\prod } \right. \prod _{j\in
[n]}(a_{i}-x_{j})\left. \vphantom{\prod } \right).
\end{equation}
Indeed, $I$ clearly defines $Z(n,l)$ set-theoretically, and it is a
complete intersection ideal, since $Z(n,l)$ has codimension $l$, while
$I$ has $l$ generators.  It is easy to see that $V(I)$ is generically
reduced: if $P$ is a point where the $x_{j}$ are all distinct, then
only one $W_{f}$ passes through $P$, and each factor $(a_{i}-x_{j})$
in \eqref{e:poly/univariate-I} with $j\not =f(i)$ is a unit in the
local ring $k[\xx ,\aa ]_{P}$, so $I_{P}$ coincides with the local
ideal $(I_{f})_{P}$.  A generically reduced complete intersection is
reduced, so $I = I(Z(n,l))$.

Now $R(n,l) = \Ocal (Z(n,l))$ is a complete intersection ring, hence
Cohen-Macaulay, and since $Z(n,l)$ is finite over $E^{n}$, the
variables $\xx $ form a homogeneous system of parameters.  This
implies that $R(n,l)$ is a free $k[\xx ]$-module.  In fact, it is free
with basis consisting of monomials $\aa ^{e} = a_{1}^{e_{1}}\cdots
a_{l}^{e_{l}}$, $0\leq e_{i}<n$, since these monomials span modulo $I$
and their number is equal to $n^{l}$, which is the number of $W_{f}$'s
and thus the degree of the finite flat morphism $Z(n,l)\rightarrow
E^{n}$.


Returning to polygraphs in two sets of variables, the analog of
\eqref{e:poly/univariate-I}, namely
\begin{equation}\label{e:poly/bivariate-I}
I=\sum _{i\in [l]}\prod _{j\in [n]}(a_{i}-x_{j},b_{i}-y_{j}),
\end{equation}
clearly defines $Z(n,l)$ set-theoretically, but now it is not a
complete intersection ideal, and does {\it not} define $Z(n,l)$ as a
reduced scheme, that is, $I\not =\rad{I}$.  At present, we do not have
a good conjecture as to a set of generators for the full ideal
$I(Z(n,l))$ in general.  To get a feeling for the possibilities, the
reader might look ahead to \S \ref{ss:poly/n=2},
eq.~\eqref{e:poly/I(Z(2,l))-gens}, where generators for the ideal are
given in the case $n=2$.

In the bivariate situation, $Z(n,l)$ is not Cohen-Macaulay.  In fact,
$Z(2,1)$ is isomorphic to $E\times Y$, where $Y$ is the union of two
linear $2$-spaces in $\AA ^{4}$ that meet only at the origin.  This
$Y$ is the classic simplest example of an equidimensional affine
algebraic set whose coordinate ring is not Cohen-Macaulay (and also
not a complete intersection, even set-theoretically).

The property that does extend from the univariate case, according to
Theorem~\ref{t:polygraphs-are-free}, is that $R(n,l)$ is free over the
ring of polynomials in {\it one} set of the $\xx ,\yy $ variables,
which by symmetry, we have taken without loss of generality to be $\yy
$.  We expect the analog of Theorem~\ref{t:polygraphs-are-free} to
hold in $d$ sets of variables, with freeness over one set of the
variables (Conjecture \ref{conj:other/d-variate}).  The proof given
here, however, is specific to the bivariate case.


As a further example, let us consider the case $l=1$.  We will give a
simple proof of Theorem~\ref{t:polygraphs-are-free} for $Z(n,1)$,
which also works in $d$ sets of variables after some obvious
modifications.

For $l=1$, we write $W_{j}$ instead of $W_{f}$, where $j = f(1)$.  The
subspaces $W_{j} = V(a_{1}-x_{j},b_{1}-y_{j})$ meet transversely,
since a change of variables to $x_{j}' = x_{j}-a_{1}$, $y_{j}' =
y_{j}-b_{1}$ makes them into coordinate subspaces.  It follows that
the ideal $I(Z(n,1))$ of their union is the product of the ideals
$I_{j}$.  In other words, for $l=1$, the ideal $I(Z(n,1))$ is equal to
the ideal in \eqref{e:poly/bivariate-I}.

Now $W_{n}$ projects isomorphically on $E^{n}$, the projection being
given in coordinates by the substitutions $a_{1}\mapsto x_{n}$,
$b_{1}\mapsto y_{n}$.  These substitutions carry the ideal of
$W_{1}\cup \cdots \cup W_{n-1}$ to 
\begin{equation}\label{e:poly/Wn-intersect-rest}
\prod _{j=1}^{n-1} (x_{n}-x_{j},y_{n}-y_{j}).
\end{equation}
This shows that the scheme-theoretic intersection of $W_{n}$ with
$W_{1}\cup \cdots \cup W_{n-1}$ is isomorphic to $Z(n-1,1)$, the
coordinates $a_{1},b_{1}$ on $Z(n-1,1)$ being identified with
$x_{n},y_{n}$.  Since $W_{1}\cup \cdots \cup W_{n-1} = Z(n-1,1)\times
E$, we have an exact sequence of $k[\xx ,\yy ,a_{1}, b_{1}]$-modules
\begin{equation}\label{e:poly/l=1-exact-seq}
0\rightarrow R(n,1)\rightarrow k[\xx ,\yy ]\oplus (R(n-1,1)\otimes
k[x_{n},y_{n}])\rightarrow R(n-1,1) \rightarrow 0.
\end{equation}
The middle term here is the direct sum of the coordinate rings of
$W_{n}$ and $W_{1}\cup \cdots \cup W_{n-1}$; the outer terms are the
coordinate rings of their union and their scheme-theoretic
intersection.  By induction on $n$ (the case $n=1$ is trivial), the
middle term is a free $k[\yy ]$-module and the last term is a free
$k[y_{1},\ldots,y_{n-1}]$-module.  All terms are graded $k[\yy
]$-modules, finitely generated in each $x$-degree (see \S
\ref{ss:poly/generic-hilbert} for more explanation).  It follows that
the first term is a free $k[\yy ]$-module.


\subsection{Overview}
\label{ss:poly/overview}

Before embarking on the proof of Theorem~\ref{t:polygraphs-are-free},
let us outline our method.  We will actually prove a stronger result,
by induction on $n$ and $l$, namely, that $R=R(n,l)$ is free with a
basis $B$ which is a common basis for a certain family of ideals in
$R$.  By this we mean that each ideal in the family is itself a free
module, spanned by a subset of the overall basis $B$.  We remark that
this is only possible because the family of ideals we consider is
rather special: inside the lattice of all ideals of $R(n,l)$, with
intersection and sum as lattice operations, our family generates a
distributive sublattice $\Lcal $, all of whose members are radical
ideals, and such that $R(n,l)/I$ is a free $k[\yy ]$-module for all
$I\in \Lcal $.

The ideals $I$ in our family will be the ideals of certain subspace
arrangements contained as closed subsets within the polygraph
$Z(n,l)$.  Before proceeding further, we define these new subspace
arrangements, and state the theorem to be proven in its full strength.


\begin{defn}\label{d:poly/Y(m,r,k)}
Let $Z(n,l)$ be a polygraph.  Given integers $r\in [n]\cup \{0 \}$,
$k\in [l]\cup \{0 \}$, and $m$, we denote by $Y(m,r,k)$ the subspace
arrangement
\begin{equation}\label{e:poly/def-of-Y(m,r,k)}
Y(m,r,k) = \bigcup _{f,T} V(x_{j}:j\in T)\cap W_{f},
\end{equation}
where $f$ ranges over functions $f\colon [l]\rightarrow [n]$, as in
Definition \ref{d:poly/Z(n,l)}, and $T$ ranges over subsets of $[n]$
such that
\begin{equation}\label{e:poly/T-in-Y(m,r,k)-def}
|T \cap  [r]\setminus f([k])|\geq m.
\end{equation}
We denote by $I(m,r,k)$ the ideal of $Y(m,r,k)$ as a
reduced closed subscheme of $Z(n,l)$.
\end{defn}

Note that in some cases $Y(m,r,k)$ is trivial---either empty or equal
to $Z(n,l)$.  Specifically, if $m\leq 0$ then $Y(m,r,k)$ is the whole
of $Z(n,l)$, and $I(m,r,k)=0$.  If $m>r$, or if $m=r=n$ and $k>0$,
then $Y(m,r,k)$ is empty, and $I(m,r,k)=(1)$.  Of course we could have
simply ruled out the trivial cases by definition, but it will simplify
notation later on to admit them.

Roughly speaking, $Y(m,r,k)$ consists of points at which $x_{j}$
vanishes for at least $m$ indices $j\in [r]$ such that $j\not =f(i)$
for $i\leq k$.  Stated this way, the criterion for membership in
$Y(m,r,k)$ is ambiguous for points lying on more than one $W_{f}$.  A
precise formulation will be given in \S \ref{ss:poly/Y(m,r,k)}.


Now we state the full theorem to be proven by induction.

\begin{thm}\label{t:common-basis}
The coordinate ring $R(n,l)$ of the polygraph $Z(n,l)$ is a free
$k[\yy ]$-module with a basis $B$ such that every ideal $I(m,r,k)$ is
spanned as a $k[\yy ]$-module by a subset of $B$.
\end{thm}

The proof of Theorem~\ref{t:common-basis} will occupy us for most of
the rest of this section.  The basic strategy is to construct the
common ideal basis $B$ explicitly by induction, using bases that we may
assume are already given in $R(n-1,l)$ and $R(n,l-1)$.

In order to show that the set we construct is a basis, we rely on a
simple but crucial algebraic device.  Because we only consider
arrangements of subspaces on which the $\yy $ coordinates are
independent, their coordinate rings are torsion-free $k[\yy
]$-modules.  To prove that a subset of a torsion-free $k[\yy ]$-module
is a free module basis, it suffices to verify it locally on an open
locus $\hat{U}_{2}\subseteq \Spec k[\yy ]$ whose complement has
codimension two.  In our situation there is a natural choice of the
open set $\hat{U}_{2}$, namely, the set of points where at most two of
the $\yy $ coordinates coincide.  We shall see that the local geometry
of $Z(n,l)$ over $\hat{U}_{2}$ essentially reduces to the case $n=2$.
In this case we are able to verify everything explicitly.


There are several sources of difficulty that make the proof more
complicated than the program just outlined would suggest.  The first
difficulty is that we have to use facts we can only verify over
$\hat{U}_{2}$, not only to prove that our elements form a basis, but
to construct them in the first place.  This forces us to employ a
delicate general method.  From a common ideal basis constructed in a
smaller ring by induction, we first produce a less restrictive kind of
basis.  We then define elements of our new common ideal basis in terms
of the less restrictive basis, allowing the coefficients to be
rational functions of $\yy$.  Finally we use information from the case
$n=2$ to show that the coefficients are regular on $\hat{U}_{2}$,
and hence they are polynomials.

Even in the construction of the less restrictive basis, the role of
the inductively constructed common ideal basis is quite subtle.  Our
ability to make use of it depends on fortunate scheme-theoretic
relationships between the subspace arrangements $Y(m,r,k)$ and certain
other special arrangements.  Most of the work in \S\S
\ref{ss:poly/lifting-II} and \ref{ss:poly/basis} goes into establishing
and carefully exploiting these relationships.

The easier, preliminary portions of the proof are in \S \S
\ref{ss:poly/local-geometry}--\ref{ss:poly/lifting-I}.  Some of this
preliminary material is conceptually fundamental, especially the local
reduction in \S \ref{ss:poly/local-geometry}, the description of the
basis and the ideals $I(m,r,k)$ in the case $n=2$ in \S
\ref{ss:poly/n=2}, and the general basis lifting method in \S
\ref{ss:poly/lifting-I}.  The preliminary material also includes a
number of lemmas giving assorted details about the $n=2$ picture and
the corresponding local picture over $\hat{U}_{2}$ for use later on.
At a first reading, it may be wise to skip over the proofs of some of
these lemmas.


As a guide to the reader we summarize here the contents of the
remaining subsections.

\begin{itemize}

\item [\S \ref{ss:poly/local-geometry}] Definition of the open sets
$U_{1}$ and $U_{2}$ and reduction of the local geometry of $Z(n,l)$ to
the case $n=2$.

\item [\S \ref{ss:poly/generic-hilbert}] Double grading of
$R(n,l)$ and computation of $x$-degree Hilbert series for generic $\yy $.

\item [\S \ref{ss:poly/n=2}] Full working out of the case $n=2$ and
some of its consequences.

\item [\S \ref{ss:poly/Y(m,r,k)}] Further information about the
arrangements $Y(m,r,k)$.

\item [\S \ref{ss:poly/lifting-I}] General lifting principle for
extending common ideal bases to schemes with an extra coordinate.

\item [\S \ref{ss:poly/lifting-II}] Application of the lifting
principle in the case of special arrangements.

\item [\S \ref{ss:poly/basis}] The three stages of the basis
construction procedure.

\item [\S \ref{ss:poly/proof}] Proofs of
Theorems~\ref{t:polygraphs-are-free} and \ref{t:common-basis} and
Proposition \ref{p:hilbert/J-is-free}.

\item [\S \ref{ss:poly/ground-ring}] Extension of
Theorem~\ref{t:polygraphs-are-free} to arbitrary ground rings.

\end{itemize}
The geometric results in Section \ref{s:hilbert} only depend on the
contents of \S \ref{ss:poly/first-defs} and \S \S
\ref{ss:poly/overview}--\ref{ss:poly/proof}.  The material in \S
\ref{ss:poly/examples} and \S \ref{ss:poly/ground-ring} has been
included for the sake of illustration and completeness.


We close this section with a brief discussion of the base cases for
the induction.  Theorem~\ref{t:common-basis} is essentially trivial
for $n=1$.  This case can therefore serve as the base of induction on
$n$, although we still need to work out the case $n=2$ in full detail
since it is used in the induction step for $n>1$.  The base of
induction on $l$ will be the case $l=0$.  This case is not altogether
trivial, and already usefully illustrates both
Theorem~\ref{t:common-basis} and the definition of the arrangements
$Y(m,r,k)$.  For this reason we give it here rather than later.

\begin{lem}\label{l:poly/l=0}
Theorem~\ref{t:common-basis} holds for $Z(n,0)$.
\end{lem}

\begin{proof}
We have $Z(n,0) = E^{n}$ and $R(n,0) = k[\xx ,\yy ]$.  Obviously
$R(n,0)$ is a free $k[\yy ]$-module, but there are also the ideals
$I(m,r,k)$ to consider.  Necessarily we must have $k=0$.  Then
$Y(m,r,0)$ is the union of subspaces in $E^{n}$ defined by the
vanishing of at least $m$ of the coordinates $x_{1},\ldots,x_{r}$.
Its ideal $I(m,r,0)$ is generated by all square-free monomials $\prod
_{j\in T}x_{j}$, where $T$ is a subset of $[r]$ of size $|T| = r-m+1$.
The set $B$ of all monomials in the $\xx $ coordinates is a free
$k[\yy ]$-module basis of $R(n,0)$, with subsets spanning every ideal
generated by monomials in $\xx $.  In particular, each ideal
$I(m,r,0)$ is spanned by a subset of $B$.
\end{proof}


\subsection{Local geometry of $Z(n,l)$}
\label{ss:poly/local-geometry}

Nearly everything in our proof of Theorem~\ref{t:common-basis}
ultimately depends on a process of local geometric reduction over
certain open sets $\hat{U}_{1}$ and $\hat{U}_{2}$ to the cases $n=1$
(which is trivial) and $n=2$ (which we will examine in detail).  Here
we define the relevant open sets, and set up the required algebraic,
geometric, and notational machinery.

\begin{defn}\label{d:poly/U_k}
The set $\hat{U}_{k}$ is the open locus in $\Spec k[\yy ]$ where the
coordinates $y_{1},\ldots,y_{n}$ assume at least $n+1-k$ distinct
values.  In other words, $\hat{U}_{k}$ is the complement of the union
of all linear subspaces $V = V(y_{i_{1}}-y_{j_{1}}, \ldots,
y_{i_{k}}-y_{j_{k}})$ defined by $k$ independent forms
$(y_{i}-y_{j})$.  In particular, $\hat{U}_{1}$ is the locus where the
$y_{i}$ are all distinct, and $\hat{U}_{2}$ is the locus where there
is at most one coincidence $y_{p}=y_{q}$.  For any scheme $\pi \colon
Z\rightarrow \Spec k[\yy ]$ over $\Spec k[\yy ]$, we define
$U_{k}\subseteq Z$ to be the open set $\pi ^{-1}(\hat{U}_{k})$.
\end{defn}

The definition of $U_{k}$ involves an abuse of notation, since we
might, for instance, have $U_{k}$ defined as a subset of $E^{n}\times
E^{l}$ in one place and as a subset of $Z(n,l)$ in another.  In
practice it will be clear from the context what is meant.  Note that
the definitions are consistent in the sense that the subset
$U_{k}\subseteq Z(n,l)$ is the intersection of $Z(n,l)$ with the
subset $U_{k}$ defined in $E^{n}\times E^{l}$.


To treat $R(n,l)$ as a $k[\yy ]$-module, we will want to localize with
respect to prime ideals in $k[\yy ]$, that is, at points $Q\in
\hat{U}_{k}\subseteq \Spec k[\yy ]$.  To extract local geometric
information about $Z(n,l)$ as a subscheme of $E^{n}\times E^{l}$, by
contrast, we want to localize at points $P\in U_{k}\subseteq
E^{n}\times E^{l}$.  A simple technical lemma relates these two types
of localization, as follows.

\begin{lem}\label{l:poly/hat(U)-vs-U}
Let $R$ be a $k[\yy ]$-algebra, let $\pi \colon \Spec R\rightarrow
\Spec k[\yy ]$ be the projection on the $\yy $ coordinates, let
$\hat{U}$ be an open subset of $\Spec k[\yy ]$, and let $U = \pi
^{-1}(\hat{U})$, as in Definition \ref{d:poly/U_k}.  If $I,J\subseteq
R$ are ideals such that $I_{P}=J_{P}$ locally for all $P\in U$
(localized as $R$-modules), then $I_{Q}=J_{Q}$ for all $Q\in \hat{U}$
(localized as $k[\yy ]$-modules).
\end{lem}

\begin{proof}
The points of $\Spec R_{Q}$ are exactly the ideals $P_{Q}$, where
$P\in \Spec R$ is such that $\pi (P)\subseteq Q$.  In particular,
every such $P$ belongs to $U$, for $Q\in \hat{U}$.  Since
$(I_{Q})_{P_{Q}} = I_{P} = J_{P} =(J_{Q})_{P_{Q}}$ for all such $P$,
we have $I_{Q}=J_{Q}$.
\end{proof}


We will really only be interested in the open sets $U_{1}$ and
$U_{2}$, where we can fully understand the local geometry of $Z(n,l)$.
Since $Z(n,l)$ and every other subspace arrangement we consider
consists of subspaces on which the $\yy $ coordinates are independent,
$U_{1}$ and $U_{2}$ are dense, and have complements of codimension $1$
and $2$ respectively.  As we shall see, it is possible to extract
significant global geometric information from the careful use of local
information on $U_{1}$ or $U_{2}$.  The following lemma shows that the
local geometry of $Z(n,l)$ on $U_{1}$ is essentially trivial.

\begin{lem}\label{l:poly/U_1-local-pix}
The components $W_{f}$ of $Z(n,l)$ have disjoint intersections with
$U_{1}$.  Thus for $P\in U_{1}\cap Z(n,l)$, there is a unique $W_{f}$
containing $P$, and $Z(n,l)$ coincides locally with $W_{f}$, that is,
$I(Z(n,l))_{P}=(I_{f})_{P}$.
\end{lem}

\begin{proof}
Let $P$ be a point of $W_{f}\cap W_{g}$, where $f(i)\not =g(i)$ for
some $i$.  On $W_{f}$ we have $b_{i}=y_{f(i)}$, while on $W_{g}$ we
have $b_{i}=y_{g(i)}$.  Hence $P\in V(y_{f(i)}-y_{g(i)})$, so $P\not
\in U_{1}$.
\end{proof}


Here is one easy but useful consequence of the local picture on
$U_{1}$.  Recall that the {\it lattice of ideals} in a ring $R$ is the
set of all ideals in $R$, with meet and join operations given by
intersection and sum.

\begin{lem}\label{l:poly/Y-reduced-U1}
Let $\Lcal $ be the sublattice of the lattice of ideals in $R(n,l)$
generated by the ideals of all subspaces of the form
\begin{equation}\label{e:poly/Y-type-space}
V(x_{j}:j\in T)\cap W_{f}.
\end{equation}
Then for every $I\in \Lcal $, $V(I)\cap U_{1}$ is reduced, that is,
$I_{P} = \rad{I}_{P}$ for all $P\in U_{1}$.
\end{lem}

\begin{proof}
At a point $P\in U_{1}$, $Z(n,l)$ coincides locally with $W_{g}\cong
E^{n}$ for a unique $g$.  Identifying $\Ocal (W_{g})$ with $\Ocal
(E^{n}) = k[\xx ,\yy ]$, the ideal of the subspace in
\eqref{e:poly/Y-type-space} is locally either $(x_{j}:j\in T)$, if
$g=f$, or $(1)$, otherwise.  It follows that every $I\in \Lcal $
coincides locally with an ideal in $k[\xx ,\yy ]$ generated by
square-free monomials in the variables $\xx $, and for such an ideal
we have $I=\rad{I}$.
\end{proof}

\begin{cor}\label{c:poly/Y-reduced-U1}
If $I$ belongs to the lattice generated by the ideals $I(m,r,k)$ in
$R(n,l)$, then $V(I)\cap U_{1}$ is reduced.
\end{cor}


Next we want to give the analog of Lemma \ref{l:poly/U_1-local-pix} for
$U_{2}$, showing that the local geometry of $Z(n,l)$ at a point of
$U_{2}$ essentially reduces to the case $n=2$.

To state the next lemma precisely, we will need to consider the
following type of situation.  Let $Z\subseteq Z(n,l)$ be the union of
those components $W_{f}$ of $Z(n,l)$ for which $f(i)$ takes some
assigned value $h(i)$, for all $i$ in a subset $I$ of the index set
$[l]$.  On $Z$ we have identically $a_{i}=x_{h(i)}$, $b_{i}=y_{h(i)}$
for $i\in I$.  Using these equations to eliminate the coordinates
$a_{i},b_{i}$ for $i\in I$ we see that $Z$ is isomorphic to a
polygraph $Z(n,l-j)\subseteq E^{n}\times E^{l-j}$, where $j=|I|$.

This situation creates a problem of notation, as the natural index set
for the coordinates on $E^{l-j}$ here is not $[l-j]$ but $L =
[l]\setminus I$.  Similar problems arise involving the index set
$[n]$, for example if $Z$ is the union of those components $W_{f}$ for
which $f$ takes values in a subset $N\subseteq [n]$.  As such
situations will arise repeatedly in what follows, we adopt the
following notational convention to deal with them.

\begin{conv}\label{conv:N,L}
Let $N$ and $L$ be finite sets of positive integers of sizes $|N|=n$,
$|L|=l$.  To every construct that we will define in terms of $n$ and
$l$, there is a corresponding construct in which the roles of $[n]$
and $[l]$ are played by $N$ and $L$, respectively.  We will refer to
the $N$, $L$ version as the construct {\it in indices $N$, $L$}.
\end{conv}

A few examples should suffice to make the meaning of this convention
clear.  In indices $N$, $L$, the ambient space $E^{n}\times E^{l}$ is
replaced by the space $E^{N}\times E^{L}$ with coordinates
\begin{equation}\label{e:poly/N,L-coordinates}
\xx _{N},\yy _{N},\aa _{L},\bb _{L} =
x_{j_{1}},y_{j_{1}},\ldots,x_{j_{n}},y_{j_{n}},\;
a_{i_{1}},b_{i_{1}},\ldots,a_{i_{l}},b_{i_{l}},
\end{equation}
where $N = \{j_{1},\ldots,j_{n} \}$ and $L = \{i_{1},\ldots,i_{l} \}$.
Functions $f:[l]\rightarrow [n]$ become functions $f:L\rightarrow N$,
the subspaces $W_{f}$ in indices $N$, $L$ have the obvious meaning,
and their union is the polygraph $Z(N,L)$ in indices $N$, $L$.  Its
coordinate ring is $R(N,L)$.  The arrangements $Y_{N,L}(m,r,k)$ and
their ideals $I_{N,L}(m,r,k)$ are defined as in \ref{d:poly/Y(m,r,k)}
but with $[r]$ and $[k]$ referring to the smallest $r$ elements of $N$
and the smallest $k$ elements of $L$, respectively.


Using the above convention we can readily describe the local geometry of
$Z(n,l)$ at points of $U_{2}$.

\begin{lem}\label{l:poly/U_2-local-pix}
Let $P$ be a point of $U_{2}\setminus U_{1}$ and let $\{p,q \}$ be the
unique pair of indices such that $P\in V(y_{p}-y_{q})$.  Let ${}\sim
{}$ be the equivalence relation on functions $f\colon [l]\rightarrow
[n]$ defined by $f\sim g$ if and only if for all $i\in [l]$, $y_{f(i)}
- y_{g(i)}$ vanishes at $P$, that is, $f(i)=g(i)$ or $\{f(i),g(i)
\}=\{p,q \}$.
\begin{itemize}
\item [(i)] We have $P\in W_{f}$ only for $f$ in a unique
${}\sim$-equivalence class $F$, so $Z(n,l)$ coincides locally at $P$
with
\begin{equation*}		
Z=\bigcup _{f\in F} W_{f}.
\end{equation*}
\item [(ii)] Let $h$ be a member of $F$, let $N=\{p,q \}$ and let
$L=h^{-1}(N)$ (note that $L$ depends only on $F$).  The projection of
$Z$ on the coordinates $\xx$, $\yy $, $\aa _{L}$, $\bb _{L}$ is an
isomorphism
\begin{equation*}		
Z\cong E^{[n]\setminus N}\times Z(N,L),
\end{equation*}
where $Z(N,L)$ is the polygraph in indices $N$ and $L$.
\end{itemize}
\end{lem}

\begin{proof}
If $P\in W_{f}\cap W_{g}$ then for all $i$ we have $P\in
V(y_{f(i)}-y_{g(i)})$, just as in the proof of Lemma
\ref{l:poly/U_1-local-pix}.  This implies (i).  On $Z$ we have
identically $a_{i} = x_{h(i)}$, $b_{i} = y_{h(i)}$ for $i\not \in L$,
so the coordinate ring of $Z$ is generated by the remaining variables,
namely $\xx $, $\yy $, $\aa _{L}$, $\bb _{L}$.  This implies that the
projection on these coordinates is an isomorphism of $Z$ onto its
image, which is clearly $E^{[n]\setminus N}\times Z(N,L)$.
\end{proof}


Reasoning as in the proof of part (ii) of the preceding lemma, we also
obtain the following ideal-theoretic result, which we record for
future reference.

\begin{lem}\label{l:poly/elimination}
Let $I\subseteq [l]$ and $h:I\rightarrow [n]$ be given, and set
\begin{equation}\label{e:poly/Z-elim}
Z = \bigcup W_{f},\quad \text{over $f:[l]\rightarrow [n]$ such that
$f|_{I} = h$},
\end{equation}
a subarrangement of $Z(n,l)$.  Setting $L = [l]\setminus I$, the
projection of $Z$ on the coordinates $\xx $, $\yy $, $\aa _{L}$, $\bb
_{L}$ is an isomorphism $Z\cong Z([n],L)$, and the ideal of $Z$ as a
closed subscheme of $Z(n,l)$ is given by
\begin{equation}\label{e:poly/I-elim}
I(Z) = \sum _{i\in I}(a_{i}-x_{h(i)},b_{i}-y_{h(i)}).
\end{equation}
\end{lem}

\begin{proof}
Let $\pi $ be the projection on the coordinates $\xx ,\yy ,\aa _{L},
\bb _{L}$.  On $Z$ we have identically $a_{i} = x_{h(i)}$, $b_{i} =
y_{h(i)}$ for $i\not \in L$.  As in the proof of the preceding lemma,
this implies that $\pi $ induces an isomorphism of $Z$ onto its image,
which in this case is $Z([n],L)$.  Note that $Z([n],L)$ is also the
image under $\pi $ of the whole of $Z(n,l)$.

Let $J$ be the ideal on the right-hand side in \eqref{e:poly/I-elim},
and let $\tilde{Z} = V(J)$.  Then $\tilde{Z}$ is a closed subscheme of
$Z(n,l)$, conceivably non-reduced.  Modulo $J$ we again have $a_{i} =
x_{h(i)}$, $b_{i} = y_{h(i)}$ for $i\not \in L$, so $\pi $ induces an
isomorphism of $\tilde{Z}$ onto its (scheme-theoretic) image.  But we
have $Z\subseteq \tilde{Z}\subseteq Z(n,l)$, and $Z([n],L)$ is the
image of both $Z$ and $Z(n,l)$, so it is also the scheme-theoretic
image $\pi (\tilde{Z})$.  The isomorphism $\pi \colon Z\cong Z([n],L)$
factors as the isomorphism $\pi \colon \tilde{Z}\cong Z([n],L)$
composed with the closed embedding $Z\hookrightarrow \tilde{Z}$, so
this implies $Z = \tilde{Z}$ or equivalently $I(Z) = J$.
\end{proof}


We come now to a technical lemma which despite its simplicity is
really the motor of the whole machine.  It allows us to verify that a
purported free $k[\yy ]$-module basis (of $R(n,l)$, for instance)
really is one by restricting our attention to $\hat{U}_{2}$, where we
have good control of the local geometry.

\begin{lem}\label{l:poly/basis-on-U2}
Let $M$ be a torsion-free $k[\yy ]$-module, and let $B$ be a subset of
$M$.  Suppose that for every $Q\in \hat{U}_{2}$, the localization $M_{Q}$ is
a free $k[\yy ]_{Q}$-module with basis $B$.  Then $M$ is a free $k[\yy
]$-module with basis $B$.
\end{lem}

\begin{proof}
We are to show that every $x\in M$ can be uniquely expressed as
\begin{equation}\label{e:poly/x=sum-p_a b_a-2}
x = \sum _{\alpha }p_{\alpha }b_{\alpha }
\end{equation}
with $p_{\alpha }\in k[\yy ]$ and $b_{\alpha }\in B$.  By hypothesis,
this is true for the image of $x$ in $M_{Q}$, with $p_{\alpha }\in
k[\yy ]_{Q}$, for all $Q\in \hat{U}_{2}$.  The local ring $k[\yy
]_{Q}$ is a subring of $k[\yy ]_{(0)} = k(\yy )$, and the unique
coefficients $p_{\alpha }\in k[\yy ]_{Q}\subseteq k(\yy )$ satisfying
\eqref{e:poly/x=sum-p_a b_a-2} for any $Q$ also satisfy
\eqref{e:poly/x=sum-p_a b_a-2} for $Q=0$.  Hence they do not depend on
$Q$.  Since the complement of $\hat{U}_{2}$ has codimension $2$, every
rational function regular on $\hat{U}_{2}$ is regular everywhere.
Thus the $p_{\alpha }$ belong to $k[\yy ]$.  Since $M$ is torsion-free
and \eqref{e:poly/x=sum-p_a b_a-2} holds locally on the dense open set
$\hat{U}_{2}$, \eqref{e:poly/x=sum-p_a b_a-2} holds identically.
\end{proof}

\begin{cor}\label{c:poly/2-free-submodules}
Let $I$ and $J$ be free submodules of a torsion-free $k[\yy ]$-module
$M$ and suppose that $I_{Q} = J_{Q}$ for all $Q\in \hat{U_{2}}$.  Then $I=J$.
\end{cor}

\begin{proof}
Lemma \ref{l:poly/basis-on-U2} implies that any free $k[\yy ]$-module
basis of $I$ is also a basis of $I+J$, so $I=I+J$.  Similarly $J=I+J$.
\end{proof}


We should stress that the role of the explicit basis $B$ in Lemma
\ref{l:poly/basis-on-U2} is crucial.  In general, a torsion-free
$k[\yy ]$-module $M$ which is locally free on $\hat{U}_{2}$ certainly
need not be free.  As an example, take $M$ to be an ideal in $k[\yy ]$
with $V(M)$ non-empty and disjoint from $\hat{U}_{2}$.  Then
$M_{Q}=(1)_{Q}$ for all $Q\in \hat{U}_{2}$, so $M$ is locally free
with basis $\{1 \}$ on $\hat{U}_{2}$.  But no element {\it contained
in $M$} generates $M$ locally on $\hat{U}_{2}$, so we cannot conclude
that $M$ is free.  Indeed if $M$ were free, Corollary
\ref{c:poly/2-free-submodules} would then imply $M=(1)$.

The following companion to Lemma \ref{l:poly/basis-on-U2} is useful to
establish that a known basis of a free $k[\yy ]$-module (of $R(n,l)$,
for instance) is a common basis for a submodule or submodules (such as
the ideals $I(m,r,k)$).

\begin{lem}\label{l:poly/basis-of-submodule}
Let $B$ be a basis of a free $k[\yy ]$-module $M$.  Let $J$ be a
submodule of $M$, and suppose that $B_{1} = B\cap J$ spans $k(\yy
)\otimes J$.  Then $J=k[\yy ]B_{1}$.
\end{lem}

\begin{proof}
Let $x$ be an element of $J$.  Since $M$ is free with basis $B$ we can
write
\begin{equation}\label{e:poly/x=sum-p_a b_a-1}
x = \sum _{\alpha }p_{\alpha }b_{\alpha },
\end{equation}
with $p_{\alpha }\in k[\yy ]$ and $b_{\alpha }\in B$.  Of course this
is also the unique expression for $x$ in terms of the basis $B$ of the
$k(\yy )$-vector space $k(\yy )\otimes M$.  Hence $p_{\alpha } = 0$
for $b_{\alpha }\not \in B_{1}$, so we have $x\in k[\yy ]B_{1}$.
\end{proof}

Note that the condition that $B_{1}$ spans $k(\yy )\otimes J$ can be
checked locally on $\hat{U}_{1}$, since tensoring with $k(\yy )$ is
the same as localizing at $Q=0$, and the zero ideal
belongs to every non-empty open subset of $\Spec k[\yy ]$.  The
philosophy governing the application of Lemmas
\ref{l:poly/basis-on-U2} and \ref{l:poly/basis-of-submodule} is
therefore as follows.  To show that a candidate $B$ is a free module
basis, we can check it locally on $\hat{U}_{2}$; then to show that $B$
is a common ideal basis, we can check it locally on $\hat{U}_{1}$.


\subsection{Generic Hilbert series}
\label{ss:poly/generic-hilbert}

Our next task is to work our the case $n=2$ in detail.  This will be
done in \S \ref{ss:poly/n=2}, using what amounts to a Gr\"obner basis
argument in disguise.  To make this argument work, we need some
enumerative information in advance about the Hilbert series of
$R(n,l)$.  Here we gather the required information.  We also note the
general fact that $R(n,l)$ is doubly graded and finite over $k[\xx,
\yy]$, a fact that will later be used implicitly in several places.

The coordinate ring
\begin{equation}\label{e:poly/O(En x El)=k[x,y,a,b]}
\Ocal (E^{n}\times E^{l})=k[\xx ,\yy ,\aa ,\bb ]
\end{equation}
of $E^{n}\times E^{l}$ is doubly graded, by degree in the $\xx ,\aa $
variables (or {\it $x$-degree}) and the $\yy ,\bb $ variables ({\it
$y$-degree}) respectively.  The ideals $I_{f}$ are obviously doubly
homogeneous, and since the defining ideal of $Z(n,l)$ is their
intersection, the coordinate ring
\begin{equation}\label{e:poly/def-of-R(n,l)}
R(n,l) = k[\xx ,\yy ,\aa ,\bb ]/I(Z(n,l))
\end{equation}
is doubly graded.  All ideals considered throughout will be doubly
homogeneous, and all coordinate rings doubly graded.

By construction $Z(n,l)$ is finite over $E^{n}$, so $R(n,l)$ is a
finitely generated $k[\xx ,\yy ]$-module.  Hence if 
\begin{equation}\label{e:poly/R(n,l)-x-grading}
R(n,l)=\bigoplus _{d} R(n,l)_{d}
\end{equation}
is the grading of $R(n,l)$ by $x$-degree, then each homogeneous
component $R(n,l)_{d}$ is a finitely generated graded $k[\yy ]$-module
(graded by $y$-degree).  We have the following well-known graded
version of Nakayama's lemma.

\begin{lem}\label{l:poly/numerical-criterion-of-freeness}
Let $M$ be a finitely generated graded $k[\yy ]$-module.  If $B$ is a
set of homogeneous elements of $M$ that spans $M/\yy M$ as a
$k$-vector space, then $B$ generates $M$.  If $B$ further satisfies
$|B| = \dim _{k(\yy )}(k(\yy )\otimes M)$, then $M$ is a free $k[\yy
]$-module with basis $B$.
\end{lem}


For the coordinate ring $R$ of any union of subspaces of the form
$V(x_{j}:j\in T)\cap W_{f}$, including $Z(n,l)$ and $Y(m,r,k)$, we can
readily determine the Hilbert series of $k(\yy )\otimes R$ as a graded
$k(\yy )$-algebra (graded by $x$-degree).

\begin{lem}\label{l:poly/torsion-free}
Let $R$ be the coordinate ring of a union of subspaces
\begin{equation}\label{e:poly/arbitrary-union}
\bigcup_{C} V(x_{j}:j\in T)\cap W_{f},
\end{equation}
over some collection $C$ of pairs $T\subseteq [n]$, $f\colon
[l]\rightarrow [n]$.  Then $R$ is torsion-free as a $k[\yy ]$-module,
and the dimension of the $x$-degree homogeneous component $k(\yy
)\otimes R_{d}$ is equal to the number of pairs
\begin{equation}\label{e:poly/basis-pairs}
e\in \NN ^{n},\quad f\colon [l]\rightarrow [n]
\end{equation}
such that $e_{1}+\cdots +e_{n} = d$ and there is some $(T,f)\in C$ for
this $f$ with $e_{j} = 0$ for all $j\in T$.
\end{lem}


\begin{proof}
Let $W_{T,f}$ denote the subspace $V(x_{j}:j\in T)\cap W_{f}$, and let
$I_{T,f} = I_{f}+(x_{j}:j\in T)$ be its ideal and $R_{T,f} = k[\xx
,\yy ,\aa ,\bb ]/I_{T,f}$ its coordinate ring.  By definition, $R =
k[\xx ,\yy ,\aa ,\bb ]/I$, where $I$ is the intersection of the ideals
$I_{T,f}$ for all $(T,f)\in C$.  Since the $\yy $ coordinates are
independent on each $W_{T,f}$, the coordinate rings $R_{T,f}$ are
free, and hence torsion-free, $k[\yy ]$-modules.  The ring $R$ is
isomorphic to a subring of the direct sum $\bigoplus _{C}R_{T,f}$, so
$R$ is also a torsion-free $k[\yy ]$-module.

Let $C_{f}$ be the set of pairs $(T,f)\in C$ with a given $f$, and let
$R_{f}$ be the coordinate ring of the partial union $Z_{f} =
\bigcup_{C_{f}} W_{T,f}$ We have an injective ring homomorphism
\begin{equation}\label{e:poly/R-in-Rf-sum}
R\hookrightarrow \bigoplus _{f} R_{f}.
\end{equation}
By Lemma \ref{l:poly/U_1-local-pix}, the partial unions $Z_{f}$ have
disjoint restrictions to $U_{1}$.  By Lemma \ref{l:poly/hat(U)-vs-U},
this implies that \eqref{e:poly/R-in-Rf-sum} localizes to an
isomorphism at each point of $\hat{U}_{1}$, and in particular, upon
tensoring with $k(\yy )$.

Now the projection of $W_{f}$ on $E^{n}$ is an isomorphism, so $Z_{f}$
projects isomorphically on $V\times \Spec k[\yy ]$, where $V$ is the
union of coordinate subspaces $\bigcup _{C_{f}} V(x_{j}:j\in T)$ in
$\Spec k[\xx ]$.  The coordinate ring of $V$, say $k[\xx ]/J$, is the
face ring of a simplicial complex.  It has a homogeneous vector space
basis consisting of all monomials $\xx ^{e}$ such that there is some
$(T,f)\in C_{f}$ with $e_{j} = 0$ for all $j\in T$.  The ring $R_{f}$
in turn is a free $k[\yy ]$ module with this same basis.  Since $k(\yy
)\otimes R\cong \bigoplus _{f}k(\yy )\otimes R_{f}$, the result
follows.
\end{proof}

\begin{cor}\label{c:poly/R(n,l)-generic}
The Hilbert series of $k(\yy )\otimes R(n,l)$ as a $k(\yy )$-algebra
graded by $x$-degree is given by
\begin{equation}\label{e:poly/R(n,l)-x-hilb-series}
\sum _{d} t^{d} \dim _{k(\yy )}(k(\yy )\otimes R(n,l)_{d}) =
\frac{n^{l}}{(1-t)^{n}}.
\end{equation}
\end{cor}


In \S \ref{ss:poly/n=2} we will make use of Lemma
\ref{l:poly/numerical-criterion-of-freeness} and Corollary
\ref{c:poly/R(n,l)-generic} in the following guise.

\begin{cor}\label{c:poly/hilbert-series}
Let $B$ be a set of doubly homogeneous polynomials whose images in
$R(n,l)$ span $R(n,l)/(\yy )$ as a $k$-vector space.  Denoting the
$x$-degree of $p\in B$ by $d(p)$, suppose the degree enumerator of $B$
satisfies
\begin{equation}\label{e:poly/degree-enumerator}
\sum _{p\in B} t^{d(p)} = \frac{n^{l}}{(1-t)^{n}}.
\end{equation}
Then $R(n,l)$ is a free $k[\yy ]$-module with basis $B$.
\end{cor}

We also have a version of Corollary \ref{c:poly/R(n,l)-generic} for
the arrangements $Y(m,r,k)$.

\begin{cor}\label{c:poly/Y(m,r,k)-generic}
The Hilbert series
\begin{equation}\label{e:poly/O(Y(m,r,k))-x-hilb-series}
\sum _{d} t^{d} \dim _{k(\yy )}\left(k(\yy )\otimes \Ocal (Y(m,r,k))_{d} \right)
\end{equation}
of $k(\yy )\otimes \Ocal (Y(m,r,k))$ as a graded $k(\yy
)$-algebra is equal to the enumerator
\begin{equation}\label{e:poly/Y(m,r,k)-basis-enumerator}
\sum _{e,f}t^{|e|},\quad \text{over $e\in \NN ^{n}$, $f\colon
[l]\rightarrow [n]$ satisfying $|[r]\setminus S_{k}(e,f)|\geq m$},
\end{equation}
where $|e| = e_{1}+\cdots +e_{n}$ and $S_{k}(e,f) = \{j:e_{j}>0 \}\cup
f([k])$.
\end{cor}

\begin{proof}
The requirement that $e_{j} = 0$ for all $j\in T$, for a given $f$ and
some $T$ satisfying \eqref{e:poly/T-in-Y(m,r,k)-def}, is equivalent to
$|[r]\setminus S_{k}(e,f)|\geq m$.
\end{proof}


\subsection{The case $n=2$}
\label{ss:poly/n=2}

At this point we are ready to analyze the case $n=2$ in full detail.
Apart from our later need for the results, we hope the reader may find
that working out the case $n=2$ usefully illustrates the concepts
introduced so far.  We begin by writing down explicit polynomials that
will form the common ideal basis required by the conclusion of
Theorem~\ref{t:common-basis}.  For this discussion we fix $l$, and of
course we fix $n=2$.

To each pair $(e,f)$, for $e\in \NN ^{2}$ and $f\colon [l]\rightarrow
[2]$, we will associate a basis element $p[e,f]$, homogeneous of
$x$-degree $|e| = e_{1}+e_{2}$.


For $e=(0,0)$ we set
\begin{equation}\label{e:poly/p[e,f]-case-1}
p[e,f] = \prod _{\substack{j>1\\ f(j)\not =f(1)}} (b_{j}-b_{1}) \cdot
\begin{cases}
(b_{1}-y_{2})&	\text{if $f(1)=1$}\\
1&	\text{otherwise}.
\end{cases}
\end{equation}

For $e=(0,h)$ with $h>0$, let $f^{-1}(\{1 \}) = S\cup T$, where $S$
and $T$ are disjoint and $S$ consists of the smallest $h$ elements of
$f^{-1}(\{1 \})$, or the whole set if $h\geq |f^{-1}(\{1 \})|$.  Then
we set
\begin{equation}\label{e:poly/p[e,f]-case-x2}
p[e,f] = x_{2}^{h-|S|}\prod _{i\in S}(a_{i}-x_{1}-x_{2})\prod _{j\in
T}(b_{j}-y_{2}).
\end{equation}

For $e = (h,0)$ with $h>0$ we set
\begin{equation}\label{e:poly/p[e,f]-case-x1}
p[e,f] = x_{1}\theta p[(0,h-1),\theta f],
\end{equation}
where $\theta $ denotes the transposition $(1\;2)$, acting on $f$ in the
obvious way, and on the polynomial ring $k[\xx ,\yy ,\aa ,\bb
]$ by exchanging $x_{1}$ with $x_{2}$ and $y_{1}$ with $y_{2}$, while
fixing the coordinates $\aa $, $\bb $.

Finally, for $e=(h_{1},h_{2})$ with both $h_{1},h_{2}>0$ we set $h =
\min (h_{1},h_{2})$ and
\begin{equation}\label{e:poly/p[e,f]-case-x1x2}
p[e,f] = (x_{1}x_{2})^{h}p[e-(h,h),f].
\end{equation}


The complicated definition of the elements $p[e,f]$ is forced on us by
the requirement that they should form a common basis for the ideals
$I(m,r,k)$, with the rule for membership in $I(m,r,k)$ being given by
Lemma \ref{l:poly/n=2-common-basis}, below.  To help orient the reader
let us consider the simplest example.  For $e = (0,0)$ and $f$ equal
to the constant function $f(i) = 2$, we have $p[e,f]=1$.  For these
$e,f$, the rule in Lemma \ref{l:poly/n=2-common-basis} places $p[e,f]$
only in the ideals $I(m,r,k)$ with $m>r$ and $I(2,2,k)$ with $k>0$,
which are trivially equal to $(1)$.

\begin{lem}\label{l:poly/n=2-freeness}
For $n=2$, the coordinate ring $R(2,l)$ of $Z(2,l)$ is a free $k[\yy
]$ module with basis $B$ the set of elements $p[e,f]$ defined in
\eqref{e:poly/p[e,f]-case-1}--\eqref{e:poly/p[e,f]-case-x1x2}.
\end{lem}

\begin{proof}
We will show that $B$ spans $R(2,l)/(\yy )$ as a $k$-vector space.
Since the $x$-degree enumerator of $B$ is clearly $2^{l}/(1-t)^{2}$,
this implies the result by Corollary \ref{c:poly/hilbert-series}.

Let $B_{0} = \{p[(0,h),f] \}$ be the subset of $B$ containing only the
elements defined in \eqref{e:poly/p[e,f]-case-1} and
\eqref{e:poly/p[e,f]-case-x2}.  By \eqref{e:poly/p[e,f]-case-x1} and
\eqref{e:poly/p[e,f]-case-x1x2} we have
\begin{equation}\label{e:poly/n=2-B-from-B_0}
B = (B_{0}\cup x_{1}\theta B_{0})\cdot
\{1,x_{1}x_{2},(x_{1}x_{2})^{2},\ldots \}.
\end{equation}
It suffices to show that $B_{0}\cup x_{1}\theta B_{0}$ spans $S =
R(2,l)/((x_{1}x_{2})+(\yy ))$.  For this it suffices in turn to show
that $B_{0}$ spans $S/(x_{1})$ and $x_{1}\theta B_{0}$ spans $x_{1}S$.
Since $x_{1}x_{2}=0$ in $S$, multiplication by $x_{1}$ gives a
well-defined surjective homomorphism $S/(x_{2})\rightarrow x_{1}S$.
If $B_{0}$ spans $S/(x_{1})$, then $\theta B_{0}$ spans $S/(x_{2})$,
and therefore $x_{1}\theta B_{0}$ spans $x_{1}S$.  Hence we need only
show that $B_{0}$ spans $S/(x_{1})$.


The ideal of $Z(2,l)$ contains the ideal
\begin{equation}\label{e:poly/I(Z(2,l))-gens}
\sum _{i\in [l]}(a_{i}-x_{1},b_{i}-y_{1})(a_{i}-x_{2},b_{i}-y_{2}) +
\sum _{i,j\in [l]} \bigl( \det
\begin{bmatrix}
	a_{i}&	b_{i}&	1\\
	a_{j}&	b_{j}&	1\\
	x_{1}&	y_{1}&	1
\end{bmatrix} \bigr).
\end{equation}
It is a consequence of the argument below that $I(Z(2,l))$ is actually
equal to the ideal in \eqref{e:poly/I(Z(2,l))-gens}, but we do not
need this result.  We only verify that the generators displayed in
\eqref{e:poly/I(Z(2,l))-gens} do indeed vanish on $Z(2,l)$.  This is
clear for the generators in the first term, which is the special case
of \eqref{e:poly/bivariate-I} for $n=2$.  The determinants in the
second term vanish because on every $W_{f}$, either $f(i)=f(j)$,
making the first two rows equal, or $f(i)\not =f(j)$, making one of
the first two rows equal to the last.  From
\eqref{e:poly/I(Z(2,l))-gens} it is easy to see that the ideal
\begin{equation}\label{e:poly/I=I(Z(2,l))+(x1,y)}
I = I(Z(2,l))+(x_{1})+(\yy )
\end{equation}
contains $a_{i}^{2}-a_{i}x_{2}$, $a_{i}b_{i}$, $b_{i}^{2}$, and
$x_{2}b_{i}$ for all $i$, and $a_{j}b_{i}-a_{i}b_{j}$ for all $i<j$.

With respect to a suitable term ordering, the initial ideal of $I$
contains $a_{i}^{2}$, $a_{i}b_{i}$, $b_{i}^{2}$, and $x_{2}b_{i}$ for
all $i$, and $a_{j}b_{i}$ for all $i<j$.  Hence $S/(x_{1})=k[\xx ,\yy
,\aa ,\bb ]/I$ is spanned by monomials in $k[x_{2},\aa ,\bb ]$ not
divisible by any of these.  In other words, $S/(x_{1})$ is spanned by
monomials
\begin{equation}\label{e:poly/n=2-spanning-monoms}
x_{2}^{k}\prod _{i\in S}a_{i}\prod _{j\in T}b_{j},
\end{equation}
where every element of $S$ is less than every element of $T$, and
$k=0$ if $T\not =\emptyset $.  Let us order the monomials in
\eqref{e:poly/n=2-spanning-monoms} so that those with smaller values
of $k$ precede those with larger values, and for $k=0$, those that
don't contain $b_{1}$ as a factor precede those that do.  Then it is
easy to see that each monomial in \eqref{e:poly/n=2-spanning-monoms}
occurs as the leading term in the reduction of an element of $B_{0}$
modulo $I$.  This implies that $B_{0}$ spans $S/(x_{1})$, as desired.
\end{proof}


We now obtain the special case of Theorem~\ref{t:common-basis} for
$n=2$.

\begin{lem}\label{l:poly/n=2-common-basis}
For $n=2$, each ideal $I(m,r,k)\subseteq R(2,l)$ is spanned as a
$k[\yy ]$-module by the set of elements $p[e,f]\in B$ indexed by $e$,
$f$ satisfying
\begin{equation}\label{e:poly/n=2-I(m,r,k)-membership}
|[r]\setminus S_{k}(e,f)|<m,
\end{equation}
where $S_{k}(e,f) = \{j:e_{j}>0 \}\cup f([k])$.
\end{lem}

\begin{proof}
First we verify that the specified elements do belong to $I(m,r,k)$.
Observe (glancing ahead) that for each $m$, $r$, $k$, the ideal $I$
displayed on the right-hand side in \eqref{e:poly/n=2-I(m,r,k)-gens}
is generated by polynomials which vanish on $Y(m,r,k)$, so we have
$I\subseteq I(m,r,k)$.  It is mostly routine now to check case-by-case
that the relevant elements $p[e,f]$ belong to $I$.  The only
tricky case is to show that $p[e,f]$ belongs to $I(1,2,k)$ for $e =
(1,0)$, $f(1)=2$, and $k>0$.  In this case $p[e,f] = x_{1}\theta
p[(0,0),\theta f]$ contains a factor $x_{1}(b_{1}-y_{1})$, which
is not so obviously in $I$.  However, in this case $I$ contains
$a_{1}-x_{1}-x_{2}$, and in $R(2,l)$ we have
$(a_{1}-x_{2})(b_{1}-y_{1})=0$, so $x_{1}(b_{1}-y_{1}) =
-(a_{1}-x_{1}-x_{2})(b_{1}-y_{1})$.

Using Lemma \ref{l:poly/basis-of-submodule}, to complete the proof it
is enough to show that the specified elements $p[e,f]$ span $k(\yy
)\otimes I(m,r,k)$.  Since $B$ is a homogeneous basis of $k(\yy
)\otimes R(2,l)$, it is equivalent to show that pairs $e,f$ {\it not}
satisfying \eqref{e:poly/n=2-I(m,r,k)-membership}, counted according
to the $x$-degree $|e|$ of $p[e,f]$, are enumerated by the Hilbert
series of $k(\yy )\otimes \Ocal (Y(m,r,k))$.  This is true by
Corollary \ref{c:poly/Y(m,r,k)-generic}.
\end{proof}


\begin{cor}\label{c:poly/n=2-I(m,r,k)-gens}
In the case $n=2$, the (non-trivial) ideals $I(m,r,k)$ are generated
as ideals in $R(2,l)$ as follows:
\begin{equation}\label{e:poly/n=2-I(m,r,k)-gens}
\begin{aligned}
I(2,2,0)	& = (\xx ,\aa )\\
I(1,2,k)	& = (x_{1}x_{2})+\sum _{i\in
[k]}(a_{i}-x_{1}-x_{2},b_{i}-b_{1})\\ 
I(1,1,k)	& = (x_{1})+\sum _{i\in [k]}(a_{i}-x_{2},b_{i}-y_{2})
\end{aligned}
\end{equation}
\end{cor}

\begin{proof}
In each case the ideal listed on the right-hand side is clearly
contained in $I(m,r,k)$.  Conversely, in the proof of Lemma
\ref{l:poly/n=2-common-basis} we showed that $I(m,r,k)$ is generated
(and even spanned as a $k[\yy ]$-module) by elements belonging to the
ideal on the right-hand side.
\end{proof}


Knowing the ideals $I(m,r,k)$ for $n=2$, we now want to draw
conclusions for general $n$.  The lemmas expressing our conclusions
will come in pairs.  In each pair we first establish a fact about the
case $n=2$, then deduce its analog on $U_{2}$ for general $n$ by local
reduction.  Our first pair of lemmas extends the local reducedness
result in Corollary \ref{c:poly/Y-reduced-U1} from $U_{1}$ to $U_{2}$.
In this instance we will need the first lemma in the pair again later,
for the proof of Lemma \ref{l:poly/I_m=I(m,r,k)-on-U2} (looking ahead,
the reader may notice that Lemmas \ref{l:poly/I_m=I(m,r,k)-case-n=2}
and \ref{l:poly/I_m=I(m,r,k)-on-U2} form another pair of the same
type).

\begin{lem}\label{l:poly/V(xT) x Y(n=2)-reduced}
Fix a two-element subset $N\subseteq [n]$.  Let $R_{1}$ be the
coordinate ring of $E^{[n]\setminus N}$, so $R = R_{1}\otimes
_{k}R(N,l)$ is the coordinate ring of $E^{[n]\setminus N}\times
Z(N,l)$.  Let $\Lcal $ be the lattice of ideals in $R$ generated by
the ideals $(x_{j})$ for $j\not \in N$ and $R_{1}\otimes I_{N,l}(m,r,k)$ for
all $m$, $r$, $k$.  Then $I = \rad{I}$ for all $I\in \Lcal $.
\end{lem}

\begin{proof}
Let $\xx ',\yy '$ be the coordinates in indices $[n]\setminus N$, so
$R_{1} = k[\xx ',\yy ']$.  Obviously $R_{1}$ is a free $k[\yy
']$-module with basis $B_{1}$ the set of all monomials in the
variables $\xx '$, and each ideal $(x_{j})\cap R_{1}$ is spanned by a
subset of $B_{1}$.  By Lemma \ref{l:poly/n=2-common-basis}, $R_{2} =
R(N,l)$ is a free $k[\yy _{N}]$-module with a basis $B_{2}$ which is a
common basis for the ideals $I_{N,l}(m,r,k)$.  It follows that $R$ is
a free $k[\yy ]$-module with basis $B_{1}\otimes B_{2}$ and that every
ideal in $\Lcal $ is spanned by a subset of $B_{1}\otimes B_{2}$.
This implies that $\Lcal $ is a distributive lattice, and hence every
ideal in $\Lcal $ is an intersection of ideals of the form $J =
(x_{j}:j\in T)\otimes R_{2}+ R_{1}\otimes I$, where $I$ is a sum of
ideals $I_{N,l}(m,r,k)$.  Now $V(J) = V(x_{j}:j\in T)\times V(I)$, so
we only have to show that $V(I)$ is reduced.

Lemma \ref{l:poly/n=2-common-basis} implies that $R_{2}/I$ is a free
$k[\yy _{N}]$ module, and hence a subring of its localization
$(R_{2}/I)_{Q}$ at any point $Q\in \Spec k[\yy _{N}]$.  By Corollary
\ref{c:poly/Y-reduced-U1} the latter ring is reduced for $Q\in
\hat{U}_{1}$ (in particular, $Q=0$ suffices).
\end{proof}


\begin{lem}\label{l:poly/Y-reduced-U2}
If $I$ belongs to the lattice generated by the ideals $I(m,r,k)$ in
$R(n,l)$, then $V(I)\cap U_{2}$ is reduced.
\end{lem}

\begin{proof} We already have the result on $U_{1}$ by Corollary
\ref{c:poly/Y-reduced-U1}.  Let $P$ be a point of $U_{2}\setminus
U_{1}$.  By Lemma \ref{l:poly/U_2-local-pix}, $Z(n,l)$ coincides
locally at $P$ with a subspace arrangement isomorphic to
$E^{[n]\setminus N}\times Z(N,L)$, for a set $N$ with two elements.
Unraveling the definitions of $Y(m,r,k)$ and the local isomorphism,
we find that $Y(m,r,k)$ coincides locally, as a reduced closed
subscheme of $E^{[n]\setminus N}\times Z(N,L)$, with a subspace
arrangement
\begin{equation}\label{e:poly/U_2-local-Y(m,r,k)}
\bigcup _{T_{1},m'} V(x_{j}:j\in T_{1})\times Y_{N,L}(m',r',k').
\end{equation}

Here we set $r' = |[r]\cap N|$ and $k' = |[k]\cap L|$.  The union
ranges over $T_{1}\subseteq [n]\setminus N$ and $m'$ satisfying
$m'+|T_{1}\cap ([r]\setminus N)\setminus h([k])| \geq m$, where $h$ is
chosen so that $P\in W_{h}$.  Note that the set $([r]\setminus
N)\setminus h([k])$ does not depend on the choice of $h$.  It follows
from Lemma \ref{l:poly/V(xT) x Y(n=2)-reduced} that the ideals of the
subspace arrangements in \eqref{e:poly/U_2-local-Y(m,r,k)}, for $N$,
$L$ fixed and $m$, $r$, $k$ arbitrary, belong to a lattice consisting
entirely of radical ideals.
\end{proof}

We remark that once Theorem~\ref{t:common-basis} is established, we
can conclude {\it a posteriori} that the lattice $\Lcal $ generated by
the ideals $I(m,r,k)$ is distributive, and that the ring $R(n,l)/I$ is
a free $k[\yy ]$-module for every $I\in \Lcal $.  In fact, it can be
shown that the existence of a common free module basis is equivalent
to these conditions, plus the freeness of $R(n,l)$ itself.  Combined
with Corollary \ref{c:poly/Y-reduced-U1}, the freeness of $R(n,l)/I$
for $I\in \Lcal $ implies that $V(I)$ is reduced everywhere, not just
on $U_{2}$, or in other words, $I = \rad{I}$.  However, we do not see
a way to obtain these stronger results without appealing to
Theorem~\ref{t:common-basis}.  We must therefore content ourselves for
now with the local result on $U_{2}$, and we will ultimately prove
Theorem~\ref{t:common-basis} for general $n$ by means that require
only this local information.


We conclude our study of the case $n=2$ with another pair of lemmas
giving a reducedness property for $n=2$ and the corresponding local
property on $U_{2}$ for all $n$.  For this (and again later) we need a
preliminary lemma.  By construction, the subspaces $W_{f}$ are
invariant with respect to simultaneous translation of the $\xx $ and
$\aa $ coordinates by a common quantity.  This fact has the following
consequence.

\begin{lem}\label{l:poly/translation}
Let $R$ be the coordinate ring of any union of the subspaces
$W_{f}\subseteq E^{n}\times E^{l}$.  If $x$ is one of the coordinate
variables $\xx $, $\aa $, we have $(x) = \rad{(x)}$ in $R$.
\end{lem}

\begin{proof}
A manifestation of translation invariance is that the endomorphism
$\tau _{x}$ of $k[\xx ,\yy ,\aa ,\bb ]$ defined by the substitutions
\begin{equation}\label{e:poly/trans-sub}
x_{j}\mapsto x_{j}-x,\quad a_{i}\mapsto a_{i}-x,\quad \text{for all
$j\in [n]$, $i\in [l]$}
\end{equation}
carries the ideals $I_{f}$ into themselves, as is clear from their
definition in \eqref{e:poly/I_f}.  We therefore have $\tau
_{x}I\subseteq I$, where $I$ is the defining ideal of $R = k[\xx ,\yy
,\aa ,\bb ]/I$.  Since $\tau _{x}(x)=0$, we have an induced ring
homomorphism $\tau _{x} \colon R/(x)\rightarrow R$.  We see
immediately that if $\phi \colon R\rightarrow R/(x)$ is the canonical
projection, then $\phi \circ \tau _{x}$ is the identity map on
$R/(x)$.  Hence $R/(x)$ is isomorphic to a subring of $R$.  Since $R$
is reduced, so is $R/(x)$, so we have $(x)=\rad{(x)}$.
\end{proof}


We will need the second lemma of the following pair for the proof of
Lemma \ref{l:poly/build-Y(1,1,t-1)}.

\begin{lem}\label{l:poly/I(1,1,t-1)+(a_t)-n=2}
For $n=2$ and $t\in [l]$ the ideal $I = I(1,1,t-1)+(a_{t})$ in
$R(2,l)$ is equal to its radical.
\end{lem}

\begin{proof}
By Corollary \ref{c:poly/n=2-I(m,r,k)-gens}, we have $I =
J+(x_{1},a_{t})$, where $J = \sum _{i < t}(a_{i}-x_{2},b_{i}-y_{2})$.
By Lemma \ref{l:poly/elimination}, we have $R(2,l)/J\cong R(2,L)$
where $L = [l]\setminus [t-1]$.  Replacing $[l]$ with $L$, this
reduces the problem to showing that $I = \rad{I}$ in $R(2,l)$ for $I =
(x_{1},a_{1})$.  For this $I$ we have set-theoretically
\begin{equation}\label{e:poly/V(I)-big-unions}
V(I) = \left( \bigcup _{f(1)=1} V(x_{1})\cap W_{f} \right) \cup \left(
\bigcup _{f(1)=2} V(x_{1},x_{2})\cap W_{f} \right).
\end{equation}

By Lemma \ref{l:poly/elimination}, the ideal
$(a_{1}-x_{1},b_{1}-y_{1})$ is equal to its radical, and its zero set
is $\bigcup _{f(1)=1} W_{f}$.  Applying Lemma \ref{l:poly/translation}
with $x = x_{1}$ to the union of subspaces $\bigcup _{f(1)=1} W_{f}$,
we see that the ideal $(x_{1})+(a_{1}-x_{1},b_{1}-y_{1}) =
I+(b_{1}-y_{1})$ is equal to its radical, and hence to
$(\rad{I})+(b_{1}-y_{1})$.  The zero set of this ideal is
$V(x_{1})\cap \bigcup _{f(1)=1} W_{f}$, which is the first big union
in \eqref{e:poly/V(I)-big-unions}.


Again by Lemma \ref{l:poly/elimination}, the ideal
$(a_{1}-x_{2},b_{1}-y_{2})$ is equal to its radical, with zero set
$\bigcup _{f(1)=2} W_{f}\cong Z(2,l-1)$.  Observing that
$V(x_{1},x_{2}) = Y(2,2,0)$ in $Z(2,l-1)$, it follows from Corollary
\ref{c:poly/n=2-I(m,r,k)-gens} that the ideal of $V(x_{1},x_{2})\cap
\bigcup _{f(1)=2} W_{f}$, that is, of the second big union in
\eqref{e:poly/V(I)-big-unions}, is $(\xx ,\aa
)+(a_{1}-x_{2},b_{1}-y_{2}) = (\xx ,\aa )+(b_{1}-y_{2})$.  Since
$b_{1}-y_{1}$ does not vanish identically on any component of the
second big union, the ideal $(\rad{I})\icolon (b_{1}-y_{1})$ is equal
to the ideal of the latter.  Hence we have $(\rad{I})\icolon
(b_{1}-y_{1}) = (\xx ,\aa )+(b_{1}-y_{2})$.

We now claim that $I\icolon (b_{1}-y_{1})$ contains $(\xx ,\aa
)+(b_{1}-y_{2})$, and hence $I\icolon (b_{1}-y_{1}) = (\rad{I})\icolon
(b_{1}-y_{1})$.  For the claim, observe that
$(b_{1}-y_{1})(b_{1}-y_{2}) = 0$ in $R(2,l)$, which shows that
$b_{1}-y_{2}\in I\icolon (b_{1}-y_{1})$.  Observe also that
$(b_{1}-y_{1})(a_{1}-x_{2}) = 0$ in $R(2,l)$, which shows that
$x_{2}\in I\icolon (b_{1}-y_{1})$.  Finally observe that the
determinant
\begin{equation}\label{e:poly/det}
\det \begin{bmatrix}
	a_{i}&	b_{i}&	1\\
	a_{1}&	b_{1}&	1\\
	x_{1}&	y_{1}&	1
\end{bmatrix}
\end{equation}
from \eqref{e:poly/I(Z(2,l))-gens} reduces modulo $(x_{1},a_{1})$ to
$a_{i}(b_{1}-y_{1})$, which shows that $a_{i}\in I\icolon
(b_{1}-y_{1})$ for all $i$.

Multiplying the identity of ideals $I\icolon (b_{1}-y_{1}) =
(\rad{I})\icolon (b_{1}-y_{1})$ by $b_{1}-y_{1}$ yields $I\cap
(b_{1}-y_{1}) = (\rad{I})\cap (b_{1}-y_{1})$, and we showed above that
$I+(b_{1}-y_{1}) = (\rad{I})+(b_{1}-y_{1})$.  By the modular law for
ideals, since $I\subseteq \rad{I}$, these imply $I = \rad{I}$.
\end{proof}


\begin{lem}\label{l:poly/I(1,1,t-1)+(a_t)}
For all $n$ and for $t\in [l]$, if $I$ is the ideal
$I(1,1,t-1)+(a_{t})$ in $R(n,l)$, then $V(I)\cap U_{2}$ is reduced.
\end{lem}

\begin{proof}
Note that by Lemma \ref{l:poly/translation}, $V(a_{t})$ is reduced and
equal to the union over all $f$ of $V(x_{f(t)})\cap W_{f}$.  Thus
$(a_{t})$ belongs to the lattice of ideals in Lemma
\ref{l:poly/Y-reduced-U1}, as does $I(1,1,t-1)$.  Therefore $V(I)\cap
U_{1}$ is reduced.

For $P\in U_{2}\setminus U_{1}$, $Z(n,l)$ coincides locally with the
arrangement $Z = E^{[n]\setminus N}\times Z(N,L)$ from Lemma
\ref{l:poly/U_2-local-pix}, and we may as well replace $Z(n,l)$ with
$Z$.

First suppose $1\not \in N$.  Fixing some $h$ such that $P\in W_{h}$
as in Lemma \ref{l:poly/U_2-local-pix}, either $Y(1,1,t-1)$ is locally
empty, if $1\in h([t-1])$, or else $Y(1,1,t-1)$ coincides locally with
$V(x_{1})$.  If $t\not \in L$ then $V(a_{t}) = V(x_{j})$ for $j =
h(t)\not \in N$.  Otherwise $V(a_{t}) = E^{[n]\setminus
N}\times Y$, where $Y = V(a_{t})\cap Z(N,L)$, which is reduced by Lemma
\ref{l:poly/translation}.  In every case, the scheme-theoretic
intersection of $Y(1,1,t-1)$ and $V(a_{t})$ either is empty or
coincides locally with a product of reduced schemes, so $V(I)$ is
locally reduced when $1\not \in N$.

Suppose instead that $1\in N$.  Then $Y(1,1,t-1)$ coincides locally
with $E^{[n]\setminus N}\times Y(1,1,t'-1)$, where $t'-1 = |[t-1]\cap
L|$.  If $t\not \in L$ we have $V(a_{t}) = V(x_{h(t)})$ and the result
is immediate, as before.  Otherwise it reduces to Lemma
\ref{l:poly/I(1,1,t-1)+(a_t)-n=2}.
\end{proof}


\subsection{Further description of $Y(m,r,k)$}
\label{ss:poly/Y(m,r,k)}

In \S \S \ref{ss:poly/lifting-II} and \ref{ss:poly/basis} we will
prove a series of lemmas enabling us to extend bases of the type
specified in Theorem~\ref{t:common-basis} from the coordinate rings of
smaller arrangements to those of larger ones.  Since the
characteristic property of our bases is their compatibility with the
ideals $I(m,r,k)$, we will of course have to use some geometric facts
about the corresponding subspace arrangements $Y(m,r,k)$.  Our purpose
in this subsection is to take note of some such facts for later use.
We begin by refining our description of the arrangements $Y(m,r,k)$.

\begin{lem}\label{l:poly/Y(m,r,k)}
Let $P$ be a point of $Z(n,l)$ and let $f$ be the pointwise maximum of
all functions $h$ such that $P$ lies on $W_{h}$.  Then $P$ lies on
$W_{f}$, and $P$ belongs to $Y(m,r,k)$ if and only if the set $T =
\{j:P\in V(x_{j}) \}$ satisfies $|T\cap [r]\setminus f([k])|\geq m$.
\end{lem}


\begin{proof}
Define an equivalence relation $\sim '_{P}$ on $[n]$ by the rule
$i\sim '_{P} j$ if and only if $P\in V(x_{i}-x_{j},y_{i}-y_{j})$.  It
induces an equivalence relation $\sim _{P}$ on functions $h\colon
[l]\rightarrow [n]$ by the rule $g\sim _{P} h$ if and only if
$g(i)\sim' _{P} h(i)$ for all $i\in [l]$ (this differs from the
equivalence relation in Lemma \ref{l:poly/U_2-local-pix} in that it
depends on both the $\xx $ and the $\yy $ coordinates).  It is easy to
see that the functions $h$ such that $P$ lies on $W_{h}$ form an
equivalence class, and that the pointwise maximum of two equivalent
functions is equivalent to each of them.  This shows that $P$ lies on
$W_{f}$.

The function $f$ takes at most one value in each $\sim
'_{P}$-equivalence class, namely, the maximum element in that class.
Furthermore, that value is only in $[r]$ if the whole equivalence
class is contained in $[r]$.  Hence, among all functions in its $\sim
_{P}$-equivalence class, $f$ maximizes the size of the intersection of
$[r]\setminus f([k])$ with every $\sim '_{P}$-equivalence class.  The
set $T$ is a union of $\sim '_{P}$-equivalence classes, so it follows
that $f$ maximizes $|T\cap [r]\setminus f([k])|$.  But $P$ belongs to
$Y(m,r,k)$ if and only if the maximum of this quantity over all $h$
such that $P\in W_{h}$ is at least $m$.
\end{proof}


At several points later on we will need to know that intersections of
the arrangements $Y(m,r,k)$ consist of spaces on which the $\yy $
coordinates are independent.  Note that this is not obvious from the
definition, although it has to be true if Theorem~\ref{t:common-basis}
is to hold.  Using the preceding lemma, we can prove a precise version
of the fact we require.

\begin{lem}\label{l:poly/Y-intersect-set}
Set-theoretically, $Y(m,r,k)\cap Y(m',r',k')$ is the union of subspaces
\begin{equation}\label{e:poly/Y-intersection}
\bigcup _{f,T} V(x_{j}:j\in T)\cap W_{f},
\end{equation}
over all $f:[l]\rightarrow [n]$ and $T\subseteq [n]$ satisfying
\begin{equation}\label{e:poly/Y-intersect-T-conditions}
|T\cap [r]\setminus f([k])|\geq m \quad \text{and}\quad |T\cap
[r']\setminus f([k'])|\geq m'.
\end{equation}
\end{lem}

\begin{proof}
It is clear that $Y(m,r,k)\cap Y(m',r',k')$ contains the union in
\eqref{e:poly/Y-intersection}.  Given a point $P$ of $Y(m,r,k)\cap
Y(m',r',k')$, let $T = \{j:P\in V(x_{j}) \}$.  By the definitions of
$Y(m,r,k)$ and $Y(m',r',k')$, we have $P\in W_{f}$ for some $f$ such
that $|T\cap [r]\setminus f([k])|\geq m$, and $P\in W_{f'}$ for some
$f'$ such that $|T\cap [r']\setminus f'([k'])|\geq m'$.  {\it A
priori}, $f$ and $f'$ might be different, but Lemma
\ref{l:poly/Y(m,r,k)} allows us to take them both to be the maximum
$f$ for which $W_{f}$ contains $P$.  Hence $Y(m,r,k)\cap Y(m',r',k')$
is contained in the union in \eqref{e:poly/Y-intersection}.
\end{proof}


In the proof of Lemma \ref{l:poly/build-Y(1,1,t-1)}, we will need the
following result.  To state it we need to introduce an automorphism
$\theta $ of $E^{n}\times E^{l}$ which applies a right circular shift
to the $E^{n}$ coordinates and fixes the $E^{l}$ coordinates, that is,
\begin{equation}\label{e:poly/theta}
\theta (P_{1},\ldots,P_{n},Q_{1},\ldots,Q_{l}) =
(P_{n},P_{1},\ldots,P_{n-1},Q_{1},\ldots,Q_{l}).
\end{equation}
This shift operation turns out to play a role in all three stages of
our basis construction and will be discussed further in \S
\ref{ss:poly/basis}.  Note that the definition of the polygraph is
symmetric in the $P_{i}$'s, so $\theta $ maps $Z(n,l)$ onto itself.

\begin{lem}\label{l:poly/Y-vs-Z0}
Given $n$, $l$ and $t\in [l]$, set $L = [l]\setminus \{t \}$ and let
$\pi \colon Z(n,l)\rightarrow Z(n,L)$ be the projection on the
coordinates other than $a_{t},b_{t}$.  Let $\theta $ be as above.
Then for $r>0$ and $k<t$ we have
\begin{multline}\label{e:poly/Y-vs-Z0-k<t}
Y(1,1,t)\cup (Y(m,r,k)\cap Y(1,1,t-1))\\
	 = Y(1,1,t)\cup (\pi ^{-1}\theta Y_{n,L}(m-1,r-1,k)\cap
Y(1,1,t-1)),
\end{multline}
while for $r>0$ and $k\geq t$ we have
\begin{multline}\label{e:poly/Y-vs-Z0-k>=t}
Y(1,1,t)\cup (Y(m,r,k)\cap Y(1,1,t-1)) \\
	= Y(1,1,t)\cup (\pi ^{-1}\theta Y_{n,L}(m,r-1,k-1)\cap Y(1,1,t-1)).
\end{multline}
\end{lem}
 
\begin{proof}
Note that $Y(1,1,t)$ is the union of the subspaces $V(x_{1})\cap
W_{f}$, over all $f$ such that $1\not \in f([t])$.
Let $Z_{0}\subseteq Z(n,l)$ be the closed subset
\begin{equation}\label{e:poly/Z0-first}
Z_{0} = \bigcup _{f}V(x_{1})\cap  W_{f}:\quad f(t) = 1,\quad 1\not \in
f([t-1]).
\end{equation}
Thus $Z_{0}$ is the union of the components of $Y(1,1,t-1)$ that are
not components of $Y(1,1,t)$.  For any subset $X\subseteq Z(n,l)$, it
follows that we have
\begin{equation}\label{e:poly/Y-vs-Z0-X}
Y(1,1,t)\cup (X\cap Y(1,1,t-1))  = Y(1,1,t)\cup (X\cap Z_{0}).
\end{equation}

For $k<t$ we show that $Y(m,r,k)\cap Z_{0} = \pi ^{-1}\theta
Y_{n,L}(m-1,r-1,k)\cap Z_{0}$.  Fix a point $P\in Z_{0}$ and let $f$
be the maximum $f$ such that $W_{f}$ contains $P$.  Note that a
component $W_{h}\subseteq Z(n,L)$ contains $\pi (P)$ if and only if $h
= g|_{L}$ for some $W_{g}$ containing $P$.  Hence $h = f|_{L}$ is the
maximum $h$ such that $W_{h}$ contains $\pi (P)$.  Let $T = \{j:P\in
V(x_{j}) \}$.  By Lemma \ref{l:poly/Y(m,r,k)}, we have $P\in Y(m,r,k)$
if and only if $|T\cap [r]\setminus f([k])|\geq m$.  Since $P\in
Z_{0}$, we have $1\in T$, and since $f$ is maximum we also have $1\not
\in f([t-1])$.  Hence we have $1\in T\cap [r]\setminus f([k])$ and
$|T\cap [r]\setminus f([k])| = 1+|T\cap ([r]\setminus \{1 \})\setminus
f([k])|$.  Note that $[r]\setminus \{1 \} = \theta [r-1]$ and that
$f([k]) = h([k])$, since $k<t$.  By Lemma \ref{l:poly/Y(m,r,k)} again,
we have $\pi (P)\in \theta Y_{n,L}(m-1,r-1,k)$ if and only if $|T\cap
\theta [r-1]\setminus h([k])| \geq m-1$.  This establishes
\eqref{e:poly/Y-vs-Z0-k<t}.

For $k\geq t$, a subtle point arises, because it is not true in this
case that $Y(m,r,k)\cap Z_{0} = \pi ^{-1}\theta Y_{n,L}(m,r-1,k-1)\cap
Z_{0}$.  We can show, however, that $Y(m,r,k)$ and $\pi
^{-1}\theta Y_{n,L}(m,r-1,k-1)$ have the same intersection with
$Z_{0}\setminus Y(1,1,t)$, which suffices to establish
\eqref{e:poly/Y-vs-Z0-k>=t}.  Therefore let $P$ be a point of $Z_{0}$,
assume $P\not \in Y(1,1,t)$ and define $f$, $h$ and $T$ as above.  We
claim that $f(t)=1$.  Since $P\in Z_{0}$ we certainly have $P\in
V(a_{t}-x_{1},b_{t}-y_{1})$ but we could also have $P\in
V(a_{t}-x_{j},b_{t}-y_{j})$ for some $j>1$, in which case $f(t)$ would
be the largest such $j$.  If this occurs, however, then we have $1\sim
'_{P} j$, in the notation used in the proof of Lemma
\ref{l:poly/Y(m,r,k)}.  This implies $f(i)>1$ for all $i\in [l]$ and
hence $P\in Y(1,1,l)\subseteq Y(1,1,t)$, contrary to assumption.

Having established that $f(t)=1$, we have $1\not \in T\cap
[r]\setminus f([k])$ and $|T\cap [r]\setminus f([k])| = |T\cap
([r]\setminus \{1 \})\setminus f([k])| = |T\cap \theta [r-1]\setminus
h([k]\setminus \{t \})|$.  Now we conclude as before that $P$ belongs
to $Y(m,r,k)$ if and only if $\pi (P)$ belongs to $\theta
Y_{n,L}(m,r-1,k-1)$.
\end{proof}


\subsection{Lifting common ideal bases}
\label{ss:poly/lifting-I}

We now develop a technique for lifting a common free $k[\yy ]$-module
basis for certain ideals in a $k[\yy ]$-algebra $R$ to a basis of a
finite $R$-algebra $S$ generated by one extra coordinate variable $b$.
We give the general method here.  In \S \ref{ss:poly/lifting-II} we
apply it with $R=R(n,l-1)$ to obtain bases for the coordinate rings of
certain special arrangements $Z'$ contained in $Z(n,l)$.  This basis
lifting from $Z(n,l-1)$ to $Z'\subseteq Z(n,l)$ will be crucial for
one of the three stages of our basis construction procedure in \S
\ref{ss:poly/basis}.

\begin{defn}\label{d:poly/V_m(Z')}
Let $Z=\Spec R$ be a Noetherian affine scheme.  Let $S = R[b]/J$ be a
finite $R$-algebra generated by one variable $b$, {\it i.e.}, $Z'=
\Spec S$ is a closed subscheme of $\AA ^{1}(Z)$ finite over $Z$.  We
denote by $I_{m}(Z')$ the ideal in $R$ consisting of all elements
$a_{m-1}$ for some
\begin{equation}\label{e:poly/leading-form}
a_{m-1}b^{m-1}+\cdots +a_{1}b+a_{0}\in J\quad (a_{i}\in R).
\end{equation}
We denote by $V_{m}(Z')$ the corresponding closed subscheme
$V(I_{m}(Z'))\subseteq Z$.
\end{defn}

Note that the elements in \eqref{e:poly/leading-form} for a given $m$
constitute an $R$-submodule of $J$, so $I_{m}(Z')$ is in fact an
ideal.  If $F_{m}$ denotes the the $R$-submodule of $S$ generated by
the elements $\{1,b,\ldots,b^{m-1} \}$, then $F_{m}/F_{m-1}$ is a
cyclic $R$-module generated by $b^{m-1}$ and isomorphic to
$R/I_{m}(Z')$.  In other words, we have $Rb^{m-1}\cap F_{m-1} =
b^{m-1}I_{m}(Z')$.  Since $S$ is assumed to be finitely generated as
an $R$-module, we have $F_{r} = S$ for some $r$ and hence
$I_{m}(Z')=(1)$ for all $m>r$.

Multiplying \eqref{e:poly/leading-form} by $b$, we see that
$I_{m}(Z')\subseteq I_{m+1}(Z')$ for all $m$.  If we agree to set
$I_{0}(Z') = 0$, then we have
\begin{equation}\label{e:poly/I_m-increasing}
0=I_{0}(Z')\subseteq I_{1}(Z')\subseteq \cdots \subseteq I_{r+1}(Z') =
(1).
\end{equation}


The ideals $I_{m}(Z')$ are called {\it partial elimination ideals} by
Green \cite{GreM98}, who studied them in connection with the theory
of lexicographic generic initial ideals.  Green observed that the
geometric interpretation of the the loci $V_{m}(Z')$ has to do with
the lengths of the fibers of the morphism $Z'\rightarrow Z$, as we now
explain.

\begin{defn}\label{d:poly/mu_P(S)}
Let $S$ be a finitely generated module over a Noetherian commutative
ring $R$.  For $P\in \Spec R$ we set
\begin{equation}\label{e:poly/def-of-mu_P(S)}
\mu _{P}(S) = \dim _{K_{P}}(K_{P}\otimes _{R}S),
\end{equation}
where $K_{P} = R_{P}/PR_{P}$ is the residue field of the local ring at
$P$.
\end{defn}

By Nakayama's lemma, $\mu _{P}(S)$ is the number of elements of any
minimal set of generators of $S_{P}$ as an $R_{P}$-module.  Also, if
$S$ is an $R$-algebra, then $\mu _{P}(S)$ is the length of the fiber
of $\Spec S$ over the $K_{P}$-valued point $\Spec K_{P}\rightarrow \Spec R$
with image $P$.

\begin{lem}[\cite{GreM98}]\label{l:poly/V_m={mu_P(S)>=m}}
Let $Z$ and $Z'$ be as in Definition \ref{d:poly/V_m(Z')}.  As sets,
we have
\begin{equation}\label{e:poly/V_m={mu_P(S)>=m}}
V_{m}(Z') = \{P\in Z:\mu _{P}(S)\geq m \}.
\end{equation}
\end{lem}

\begin{proof}
Let $P$ be a point of $Z = \Spec R$.  The Artin $K_{P}$-algebra
$K_{P}\otimes _{R} S$ is a quotient of the polynomial ring $K_{P}[b]$,
and therefore has a $K_{P}$-basis $\{1,b,\ldots,b^{m_{0}-1} \}$, where
$m_{0} = \mu _{P}(S)$.  Equivalently, $\{1,b,\ldots, b^{m_{0}-1} \}$
minimally generates $S_{P}$ as an $R_{P}$-module, so with $F_{m}$ as
in the remarks following Definition \ref{d:poly/V_m(Z')}, we have
$(F_{m_{0}})_{P}=S_{P}$, $(F_{m_{0}-1})_{P}\not =S_{P}$.  This implies
$I_{m_{0}}(Z')_{P}\not =(1)$, $I_{m_{0}+1}(Z')_{P} = (1)$, that is,
$P\in V_{m_{0}}(Z')$, $P\not \in V_{m_{0}+1}(Z')$.  By
\eqref{e:poly/I_m-increasing} it follows that $P\in V_{m}(Z')$ if and
only if $m\leq \mu _{P}(S)$.
\end{proof}


The next two lemmas give some further elementary facts about partial
elimination ideals.

\begin{lem}\label{l:poly/V_m-times}
Let $Z$ and $Z'$ be as in Definition \ref{d:poly/V_m(Z')}.  Assume
they are schemes over $k$, and let $X$ be an affine scheme flat over
$k$.  Then we have
\begin{equation}\label{e:poly/V_m-times}
V_{m}(Z'\times X) = V_{m}(Z')\times X
\end{equation}
as closed subschemes of $Z\times X$.
\end{lem}

\begin{proof}
Let $Z = \Spec R$, $Z'= \Spec S$, and $X = \Spec T$, where $S =
R[b]/J$.  Then $S\otimes T = (R\otimes T)[b]/(J\otimes T)$, so the
defining ideal of $Z'\times X$ is $J\otimes T$.  Since $T$ is flat, we
have $(J\otimes T)\cap (R\otimes T)\{1, b, \ldots, b^{m-1} \}$ =
$(J\cap R \{ 1, b, \ldots, b^{m-1} \})\otimes T$.  From the definition
it then follows that $I_{m}(Z'\times X) = I_{m}(Z')\otimes T$, which
is the algebraic equivalent of \eqref{e:poly/V_m-times}.
\end{proof}

Note that the lemma as stated holds over an arbitrary ground ring $k$.
When $k$ is a field the hypothesis that $X$ is flat over $k$ is
satisfied automatically.

\begin{lem}\label{l:poly/V_m-union}
With $Z$ and $Z'$ as in Definition \ref{d:poly/V_m(Z')}, suppose that
$Z'$ is the scheme-theoretic union of closed subschemes $Z'_{1}$ and
$Z'_{2}$, that is, the defining ideal of $Z'$ is an intersection $J =
J_{1}\cap J_{2}$.  Then we have
\begin{equation}\label{e:poly/V_m-union}
I_{s}(Z'_{1})I_{t}(Z'_{2})\subseteq I_{s+t-1}(Z')
\end{equation}
for all $s$, $t$.
\end{lem}

\begin{proof}
Let $p\in I_{s}(Z'_{1})$ be the leading coefficient of an element
$f\in J_{1}$ with leading form $pb^{s-1}$, as in
\eqref{e:poly/leading-form}.  Similarly let $q\in I_{t}(Z'_{2})$ be
the leading coefficient of $g\in J_{2}$ with leading form $qb^{t-1}$.
Then we have $fg\in J_{1}J_{2}\subseteq J$, and the leading form of
$fg$ is $pqb^{m-1}$ where $m = s+t-1$.
\end{proof}


For us, the main significance of the partial elimination ideals
$I_{m}(Z')$ is their usefulness in lifting a free $k[\yy ]$-module
basis from $R$ to $S$.  Supposing that $R$ has a basis which is a
common basis for all the partial elimination ideals, it is
straightforward to construct from it a free $k[\yy ]$-module basis of
$S$.  In keeping with our program of reducing everything to the case
$n=2$, we require a stronger form of the construction, in which the
hypotheses are limited to things we can check locally on $U_{2}$ or
$U_{1}$.  We remark, however, that when the hypotheses in the
following lemma hold locally, it actually follows that they hold
everywhere---in other words, the ideals $I_{m}$ in
\eqref{e:poly/I_m(Z')=I_m-locally} below are equal to the partial
elimination ideals.

\begin{lem}\label{l:poly/I_m-basis}
Let $R$ be a $k[\yy ]$-algebra and let $Z=\Spec R$, $Z' = \Spec S$ be
as in Definition \ref{d:poly/V_m(Z')}.  Assume that $R$ is free and
$S$ is torsion-free as $k[\yy ]$-modules.  Suppose that $R$ has a
basis $B$ which is a union of disjoint subsets $B=\bigcup
_{j=0}^{r}B_{j}$ such that for each $m$, the submodule $I_{m} = k[\yy
]\bigcup _{j<m}B_{j}$ is an ideal in $R$ satisfying
\begin{equation}\label{e:poly/I_m(Z')=I_m-locally}
I_{m}(Z')_{P} = (I_{m})_{P}
\end{equation}
for all $P\in U_{2}\cap Z$.  Then $S$ is a free $k[\yy ]$-module with
basis
\begin{equation}\label{e:poly/basis-of-S}
B' = \bigcup _{m=1}^{r} B_{m}\cdot \{1,b,\ldots,b^{m-1} \} .
\end{equation}
Furthermore, let $I\subseteq R$ be an ideal spanned as a $k[\yy
]$-module by a subset $A\subseteq B$, and let $Y' = \Spec S/IS$.
Suppose that
\begin{equation}\label{e:poly/I_m(Y')-local}
I_{m}(Y')_{P} = I_{P}+I_{m}(Z')_{P}
\end{equation}
for all $m>0$ and all $P\in U_{1}\cap Z$.  Then the ideal $IS$ is
spanned by the subset $A' = A\cdot \{1,b,b^{2},\ldots \}\cap B'$ of
$B'$.
\end{lem}


\begin{proof}
Consider the filtration $0=F_{0}\subseteq F_{1}\subseteq \cdots $ of
$S$, where $F_{m} = R\cdot \{1,b,\ldots,b^{m-1} \}$.  Without loss of
generality we can assume that $F_{r} = S$, since the hypotheses and
conclusion of the Lemma are unaltered if we replace $r$ by a larger
integer $r_{1}$ and define the extra subsets $B_{j}$ for $r<j\leq
r_{1}$ to be empty.  We have $F_{m}/F_{m-1}\cong R/I_{m}(Z')$ and
therefore, by \eqref{e:poly/I_m(Z')=I_m-locally} and Lemma
\ref{l:poly/hat(U)-vs-U}, $(F_{m}/F_{m-1})_{Q}\cong (R/I_{m})_{Q}$ for
all $Q\in \hat{U}_{2}$.  It follows that $(F_{m}/F_{m-1})_{Q}$ is a
free $k[\yy ]_{Q}$-module with basis $b^{m-1}\bigcup _{j \geq
m}B_{j}$, and since $B'$ is the union of these, $S_{Q}$ is a free
$k[\yy ]_{Q}$-module with basis $B'$.  By Lemma
\ref{l:poly/basis-on-U2} this implies that $S$ is a free $k[\yy
]$-module with basis $B'$.

By \eqref{e:poly/I_m(Y')-local} and Lemma \ref{l:poly/hat(U)-vs-U}, we
have $I_{m}(Y')_{Q} = I_{Q}+I_{m}(Z')_{Q}$ for $Q\in \hat{U}_{1}$.  By
the definition of $Y'$, multiplication by $b^{m-1}$ is an isomorphism
$R/I_{m}(Y')\rightarrow (F_{m}+IS)/(F_{m-1}+IS)\cong
F_{m}/(F_{m-1}+IS\cap F_{m})$. The $R$-module
$F_{m}/(F_{m-1}+b^{m-1}I)$ is generated by $b^{m-1}$ and is therefore
equal to $Rb^{m-1}/(Rb^{m-1}\cap (F_{m-1}+b^{m-1}I))$, or to
$Rb^{m-1}/(b^{m-1}I_{m}(Z')+b^{m-1}I)$, by the modular law.
Multiplication by $b^{m-1}$ is thus an isomorphism
$R/(I_{m}(Z')+I)\rightarrow F_{m}/(F_{m-1}+b^{m-1}I)$.  Hence for
$Q\in \hat{U}_{1}$ we have $(F_{m-1}+IS\cap F_{m})_{Q} =
(F_{m-1}+b^{m-1}I)_{Q}$ and therefore $(IS\cap F_{m})_{Q}\subseteq
(b^{m-1}I + F_{m-1} \cap IS)_{Q}$.  It follows that the map
$I_{Q}/(I_{m}(Z')\cap I)_{Q}\rightarrow (IS\cap F_{m})_{Q}/(IS\cap
F_{m-1})_{Q}$ given by multiplication by $b^{m-1}$ is surjective.  It
is injective by the definition of $I_{m}(Z')$, so it is an
isomorphism.

Setting $A_{j} = A\cap B_{j}$, this implies that $b^{m-1}\bigcup _{j
\geq m}A_{j}$ is a basis of $(IS\cap F_{m})_{Q}/(IS\cap F_{m-1})_{Q}$.
Since $A'$ is the union of these, $A'$ is a basis of $IS_{Q}$.  In
particular this holds for $Q = 0$, which by Lemma
\ref{l:poly/basis-of-submodule} suffices to show that $A'$ spans $IS$.
\end{proof}


\subsection{Lifting for special arrangements}
\label{ss:poly/lifting-II}

The most difficult stage in the basis construction procedure to be
presented in \S \ref{ss:poly/basis} will involve the construction of a
free $k[\yy ]$-module basis in the ring $R(n,l)/I(1,1,t-1)$, given a
suitable basis in the ring $R(n,l)/I(1,1,t)$.  To pass from one ring
to the other we will need a basis of $I(1,1,t)/I(1,1,t-1)$.  We will
get this basis from a basis of $R(n,l-1)$ by applying the lifting
theory from \S \ref{ss:poly/lifting-I} to certain special
arrangements, which we now define.


\begin{defn}\label{d:poly/special}
Given $n$, $l$, $r$ and $k$, with $r\in [n]$ and $k\in [l]$, the {\it
special arrangement} $Z' = Z'(r,k)$ is the subspace arrangement
\begin{equation}\label{e:poly/def-of-Z'}
Z' = \bigcup V(a_{k})\cap W_{f},\quad \text{over $f\colon
[l]\rightarrow [n]$ such that $f(k)\in [r]\setminus f([k-1])$},
\end{equation}
a reduced closed subscheme of $Z(n,l)$.  We denote by $\pi \colon
Z'\rightarrow Z(n,L)$ the projection on the coordinates other than
$a_{k}$, $b_{k}$, mapping $Z'$ into the polygraph $Z(n,L)$ in indices
$[n]$ and $L = [l]\setminus \{k \}$.  For consistency, we define
$Z'(0,k)$ to be empty when $r = 0$.
\end{defn}

The significance of the special arrangements $Z'(r,k)$ is twofold.  On
the one hand, they are related to the arrangements $Y_{n,L}(m,r,k-1)$
in the ground scheme $Z(n,L)$ via our basis lifting theory.  On the
other hand, in $Z(n,l)$ there is a geometric relationship between the
closed subsets $Z'(n,t)$, $Z'(n-1,t)$, $Y(1,1,t)$ and $Y(1,1,t-1)$
that makes possible the step from from $R(n,l)/I(1,1,t-1)$ to
$R(n,l)/I(1,1,t)$.  The latter aspect will be explained in detail in
\S \ref{ss:poly/basis}.  Here we treat the basis lifting aspect.  To
start with, we need the following purely set-theoretic fact, which
should give a strong hint as to where we are headed.

\begin{lem}\label{l:poly/V_m>Y(m,r,k)}
Let $Z' = Z'(r,k)$ be the special arrangement over $Z=Z(n,L)$, where
$L = [l]\setminus \{k \}$, and let $\pi \colon Z'\rightarrow Z(n,L)$
be the coordinate projection.  Then $Y_{n,L}(m,r,k-1)\subseteq Z(n,L)$
is the closure of the locus consisting of points $P\in U_{1}\cap
Z(n,L)$ for which the fiber $\pi ^{-1}(P)$ has size $|\pi
^{-1}(P)|\geq m$.
\end{lem}

\begin{proof}
Let $P$ be a point of $U_{1}\cap Z(n,L)$, and let $Q\in \pi ^{-1}(P)$
be a point of the fiber.  Note that this implies $Q\in U_{1}\cap
Z(n,l)$.  By Lemma \ref{l:poly/U_1-local-pix} there is a unique
$W_{f}$ containing $Q$ and a unique $W_{h}$ containing $P$.  Clearly
we must have $h = f|_{L}$.  By the definition of $Z'$, we
have $f(k)\in [r]\setminus h([k-1])$, and $x_{f(k)}$ vanishes at $Q$ and
hence at $P$.  Note that $f$ determines $Q$, since $W_{f}$ projects
isomorphically on $E^{n}$, and this projection must map $Q$ to the
image of $P$.  Conversely, for every $j\in [r]\setminus h([k-1])$ such
that $x_{j}$ vanishes at $P$, there is a point $Q\in \pi ^{-1}(P)$
lying on $W_{f}$, where $f(k) = j$ and $f|_{L} = h$.

This shows that the number of points in the fiber $\pi ^{-1}(P)$ is
equal to the size of the set $T \cap [r]\setminus h([k-1])$, where $T
= \{j :P\in V(x_{j}) \}$.  By Lemma \ref{l:poly/Y(m,r,k)},
$Y_{n,L}(m,r,k-1)\cap U_{1}$ consists of those points $P$ for which
$|T \cap [r]\setminus h([k-1])|\geq m$.  Since $Y_{n,L}(m,r,k-1)\cap
U_{1}$ is dense in $Y_{n,L}(m,r,k-1)$, the result follows.
\end{proof}

We will now apply the theory developed in \S \ref{ss:poly/lifting-I}
to the case where $Z' = \Spec S$ is the special arrangement $Z'(r,k)$
over $Z = Z(n,L)$, with $L = [l]\setminus \{k \}$.  The coordinate
ring of the ground scheme is thus $R = R(n,L)$.  The map $\pi \colon
Z'\rightarrow Z$ is the projection on the coordinates other than
$a_{k},b_{k}$.  Note that $a_{k}$ vanishes on $Z'$ by definition, so
$Z'$ is a closed subscheme of $\AA ^{1}\times Z(n,L)$, where the extra
coordinate $b$ on $\AA ^{1}$ is $b_{k}$.  Note also that the
product $\prod _{j\in [n]} (b_{k}-y_{j})$ is a monic polynomial in
$b_{k}$ which vanishes identically on $Z(n,l)$ and therefore also on
$Z'$, so $Z'$ is finite over $Z(n,L)$.  Thus $Z'$ and $Z = Z(n,L)$ are
as in Definition \ref{d:poly/V_m(Z')}.


\begin{lem}\label{l:poly/I_m<I(m,r,k)}
Let $Z' = Z'(r,k)$ be the special arrangement over $Z=Z(n,L)$, where
$L = [l]\setminus \{k \}$. Then we have
\begin{equation}\label{e:poly/I(m,r,k-1)>=I_m(Z')}
I_{m}(Z')\subseteq I_{n,L}(m,r,k-1),
\end{equation}
where $I_{n,L}(m,r,k-1)\subseteq R(n,L)$ is the ideal of
$Y_{n,L}(m,r,k-1)$.
\end{lem}

\begin{proof}
Let $\pi \colon Z'\rightarrow Z(n,L)$ be the projection.  If $P$ is a
point of $Z(n,L)$ such that the fiber $\pi ^{-1}(P)$ has at least $m$
points then in the notation of Definitions \ref{d:poly/V_m(Z')} and
\ref{d:poly/mu_P(S)} we have $\mu _{P}(S)\geq m$.

By Lemmas \ref{l:poly/V_m={mu_P(S)>=m}} and \ref{l:poly/V_m>Y(m,r,k)}
it follows that $Y_{n,L}(m,r,k-1)\subseteq V_{m}(Z')$.  Since
$Y_{n,L}(m,r,k-1)$ is reduced, this implies $I_{m}(Z')\subseteq
I_{n,L}(m,r,k-1)$.
\end{proof}


\begin{lem}\label{l:poly/I_m=I(m,r,k)-case-n=2}
For $n=2$, with the notation of Lemma \ref{l:poly/I_m<I(m,r,k)}, we
have equality in \eqref{e:poly/I(m,r,k-1)>=I_m(Z')}.
\end{lem}

\begin{proof}
By Lemma \ref{l:poly/I_m<I(m,r,k)} we have $I_{m}(Z') \subseteq
I(m,r,k-1)$ in $R(2,L)$.  This containment implies equality trivially
if $I_{m}(Z') = (1)$ or $I(m,r,k-1)=0$.  For $m=0$ we have
$I(m,r,k-1)=0$.  For $r=0$, $Z'$ is empty, so $I_{m}(Z')=(1)$ for $m>0$.
Hence we may assume $m,r>0$.

For $m=1$, the definition yields that $I_{1}(Z')$ is the kernel of the
ring homomorphism $R\rightarrow S$ corresponding to the projection
$\pi \colon Z'\rightarrow Z$, that is, the ideal of the image of $Z'$
in $Z$.  The latter coincides with $Y(1,r,k-1)$ by construction.

Since $(b_{k}-y_{1})(b_{k}-y_{2})$ vanishes on $Z'$, we have
$I_{m}(Z')=(1)$ for $m>2$.  If $r=1$, then $b_{k}-y_{1}$ vanishes on
$Z'$, and if $k>1$, then $b_{k}+b_{1}-y_{1}-y_{2}$ vanishes on $Z'$.
Hence we have $I_{m}(Z')=(1)$ for $r=1$ or $k>1$, and $m>1$.  This
leaves only the case $m=r=2$, $k=1$, in which we are to show
\begin{equation}\label{e:poly/I(2,2,0)<=I_2}
I(2,2,0)\subseteq I_{2}(Z').
\end{equation}


We have already seen (Corollary \ref{c:poly/n=2-I(m,r,k)-gens}) that
$I(2,2,0) = (\xx ,\aa )$.  Thus we have only to show that $x_{1}$,
$x_{2}$, and the $a_{i}$ for $i\in L = [l]\setminus \{1 \}$ belong to
$I_{2}(Z')$.  For this we need only observe that $x_{1}$, $x_{2}$, and
$a_{i}$ are the $b_{1}$ coefficients of the following polynomials,
which vanish on $Z'$:
\begin{equation}\label{e:poly/exhibit-pols}
\begin{gathered}
(b_{1}-y_{2})x_{1}\\
(b_{1}-y_{1})x_{2}\\
(b_{1}-y_{1})a_{i}+(b_{i}-y_{2})x_{1}.
\end{gathered}
\end{equation}
For the vanishing of the last of these, note that it is congruent
modulo $(a_{1},(b_{1}-y_{2})x_{1})$ to the determinant
\begin{equation}\label{e:poly/det-again}
\det \begin{bmatrix}
	a_{i}&	b_{i}&	1\\
	a_{1}&	b_{1}&	1\\
	x_{1}&	y_{1}&	1
\end{bmatrix}
\end{equation}
from \eqref{e:poly/I(Z(2,l))-gens}, which vanishes on $Z(2,l)$.
\end{proof}


\begin{lem}\label{l:poly/I_m=I(m,r,k)-on-U2}
For all $n$, with the notation of Lemma \ref{l:poly/I_m<I(m,r,k)},
we have equality in \eqref{e:poly/I(m,r,k-1)>=I_m(Z')} locally on
$U_{2}$, that is,
\begin{equation}\label{e:poly/I(m,r,k-1)>=I_m(Z')-on-U2}
I_{m}(Z')_{P} = I_{n,L}(m,r,k-1)_{P}
\end{equation}
for all $P\in U_{2}\cap Z(n,L)$.
\end{lem}

\begin{proof}
First consider a point $P\in U_{1}\cap Z(n,L)$. By Lemma
\ref{l:poly/U_1-local-pix}, $P$ lies on a unique component $W_{f}$,
for some $f\colon L\rightarrow [n]$, and $Z(n,L)$ coincides locally
with $Z_{1} = W_{f}\cong E^{n}$.  Since $\pi \colon Z'\rightarrow
Z(n,L)$ maps $W_{g}$ into $W_{g|_{L}}$, every point of $Z'$ lying over
$P$ belongs to a component $V(a_{k})\cap W_{g}$ with $g|_{L} = f$ and
hence $g(k)\in T = [r]\setminus f([k-1])$.  We can replace $Z'$ with
the union $Z_{1}'$ of only these components, without changing $S_{P}$
and therefore without changing the local ideals $I_{m}(Z')_{P}$.

Identifying $Z_{1}$ with $E^{n} = \Spec k[\xx ,\yy ]$, $Z'_{1}$ is the
union of subspaces $V(x_{j},b_{k}-y_{j})\subseteq \Spec k[\xx ,\yy
,b_{k}]$ over all $j\in T$.  Thus the defining ideal $J_{1}$ of
$Z'_{1}$ contains
\begin{equation}\label{e:poly/J}
J = \prod _{j\in T }(x_{j},b_{k}-y_{j}).
\end{equation}
From this we see that $I_{m}(Z'_{1})$ contains all square-free
monomials of degree $|T|+1-m$ in the variables $\{x_{j}:j\in T \}$.
These monomials generate the ideal of the union of coordinate
subspaces
\begin{equation}\label{e:poly/sq-free-monom-locus}
\bigcup V(x_{j}:j\in T_{1})\cap Z_{1},\quad \text{over subsets
$T_{1}\subseteq T$ of size $|T_{1}|\geq m$},
\end{equation}
which coincides locally with $Y_{n,L}(m,r,k-1)$ at $P$.  Hence we have
$I_{n,L}(m,r,k-1)_{P}\subseteq I_{m}(Z')_{P}$, and we already have the
reverse containment by Lemma \ref{l:poly/I_m<I(m,r,k)}.


For $P\in U_{2}\setminus U_{1}$, the local picture of $Z(n,L)$ is
given by Lemma \ref{l:poly/U_2-local-pix} as $E^{[n]\setminus N}\times
Z(N,L')$ for some $N =\{p,q \}\subseteq [n]$, and $L' \subseteq L$.
Let $F$ be the equivalence class of functions $f\colon L\rightarrow
[n]$ defined in Lemma \ref{l:poly/U_2-local-pix}, containing all $f$
such that $P\in W_{f}$.  We put $Z_{1} = E^{[n]\setminus N}$ and
$Z_{2} = Z(N,L')$.

As before, $Z'$ coincides locally over $P$ with the union of its
components $V(a_{k})\cap W_{g}$, where $g|_{L} = f$ for some $f\in F$.
We divide these into two classes.  The first consists of components
with $g(k)\not \in N$.  These have a point lying over $P$ only if
$g(k)$ belongs to
\begin{equation}\label{e:poly/T-case-U2}
T = ([r]\setminus N)\setminus f([k-1]).
\end{equation}
Note that $T$ is well-defined, since $f([k-1])\cup N$ depends only on
$F$.  The union of these components is $Z'_{1}\times Z_{2}$, where
$Z'_{1} \subseteq \AA ^{1}\times Z_{1}$ is the union of the subspaces
$V(x_{j},b_{k}-y_{j})$ over $j\in T$.  Repeating our analysis of the
case $P\in U_{1}$, we see that for each $s$, $I_{s}(Z'_{1})$ contains
the ideal $I'_{s}\subseteq R_{1}$ generated by all square-free
monomials of degree $|T|+1-s$ in the variables $\{x_{j}:j\in T \}$.
Here $R_{1} = \Spec k[\xx ',\yy ']$ is the coordinate ring of $Z_{1}$.


The second class consists of components $V(a_{k})\cap W_{g}$ of $Z'$
with $g|_{L}\in F$ and $g(k)\in N$.  The union of these components is
$Z_{1}\times Z'_{2}$, where $Z'_{2}$ is the special arrangement
$Z'(r',k')$ over $Z_{2}$, with extra coordinate $b_{k}$, and $r' =
|N\cap [r]|$, $k'-1 = |L'\cap [k-1]|$.  By Lemma
\ref{l:poly/I_m=I(m,r,k)-case-n=2} we have $I_{N,L'}(t,r',k'-1)=
I_{t}(Z'_{2})$ for all $t$.  Note that this is still correct when
$r'=0$ and $Z_{2}'$ is empty.

Using Lemmas \ref{l:poly/V_m-times} and \ref{l:poly/V_m-union}, we see
that for $s+t=m+1$ we have $I'_{s}\otimes I_{N,L'}(t,r',k'-1)\subseteq
I_{m}(Z')$.  By Lemma \ref{l:poly/V(xT) x Y(n=2)-reduced}, the ideal
\begin{equation}\label{e:poly/sum-over-s,t}
I = \sum _{s+t=m+1} I'_{s}\otimes I_{N,L'}(t,r',k'-1)\subseteq
I_{m}(Z')
\end{equation}
is equal to its radical.  Here we used the identity $I'_{s}\otimes
I_{N,L'}(t,r',k'-1) = (I'_{s}\otimes R(N,L'))\cap (R_{1}\otimes
I_{N,L'}(t,r',k'-1))$, which holds because $R_{1}/I'_{s}$ and
$R(N,L')/I_{N,L'}(t,r',k'-1)$ are flat over the field $k$.  We have
\begin{equation}\label{e:poly/V(sum)}
\begin{aligned}
V(I) &	= \bigcap _{s+t=m+1}(V(I'_{s})\times Z_{2})\cup(Z_{1}\times
Y_{N,L'}(t,r',k'-1)) \\
& 	= \bigcup _{s+t=m}V(I'_{s})\times Y_{N,L'}(t,r',k'-1).
\end{aligned}
\end{equation}
The second equality here follows because the closed subsets
$V(I'_{s})\subseteq Z_{1}$ and $Y_{N,L'}(t,r',k'-1)\subseteq Z_{2}$
decrease as $s$ and $t$ increase.  Unraveling the local isomorphism of
$Z(n,L)$ with $Z_{1}\times Z_{2}$ we see that the union in
\eqref{e:poly/V(sum)} coincides locally at $P$ with
$Y_{n,L}(m,r,k-1)$.  This shows that $I_{n,L}(m,r,k-1)_{P}\subseteq
I_{m}(Z')_{P}$, and we already have the reverse containment by Lemma
\ref{l:poly/I_m<I(m,r,k)}, as before.
\end{proof}

The next three lemmas contain specific consequences of the
preceding results that we will use in \S \ref{ss:poly/basis}.
Throughout, the special arrangements are over $Z = Z(n,L)$, where $L =
[l]\setminus \{k \}$, and all ideals $I(m',r',k')$ considered are in
the ring $R(n,L)$.  For readability, we suppress the subscripts,
writing $I(m',r',k')$ for $I_{n,L}(m',r',k')$ and $Y(m',r',k')$ for
$Y_{n,L}(m',r',k')$.


\begin{lem}\label{l:poly/Z'}
Let $Z' = Z'(r,k)$ be the special arrangement over $Z = Z(n,L)$, where
$L = [l]\setminus \{k \}$.  Assume that the coordinate ring $R(n,L)$
is a free $k[\yy ]$-module with basis $B$ such that each ideal
$I(m',r',k')\subseteq R(n,L)$ is spanned by a subset of $B$.  In
particular, for $r' = r$, $k' = k-1$, we have $B = \bigcup _{j =
0}^{r}B_{j}$, where the $B_{j}$ are disjoint, such that $I(m,r,k-1) =
k[\yy ]\cdot \bigcup _{j< m} B_{j}$ for all $m$.  Then the coordinate
ring $S=\Ocal (Z')$ is a free $k[\yy ]$-module with basis
\begin{equation}\label{e:poly/basis-B'-of-Z'}
B' = \bigcup_{m=1}^{r} B_{m}\cdot \{1,b_{k},\ldots,b_{k}^{m-1} \}.
\end{equation}
Moreover, for all $m'$, $r'$, $k'$, if $A\subseteq B$ spans
$I(m',r',k')$, then $A' = A\cdot \{1,b_{k},b_{k}^{2},\ldots \}\cap B'$
spans $I(m',r',k')S$.
\end{lem}

\begin{proof}
By Lemma \ref{l:poly/torsion-free}, the coordinate ring $S$ of $Z'$ is
a torsion-free $k[\yy ]$-module.  By Lemma
\ref{l:poly/I_m=I(m,r,k)-on-U2} we have $I_{m}(Z')_{P} =
I(m,r,k-1)_{P}$ for all $m$ and all $P\in U_{2}$.  Thus the hypotheses
of Lemma \ref{l:poly/I_m-basis} are satisfied and we have
\eqref{e:poly/basis-B'-of-Z'} immediately.  For Lemma
\ref{l:poly/I_m-basis} also to yield us the last conclusion we need to
have $I(m',r',k')_{P}+I(m,r,k-1)_{P} = I_{m}(Y')_{P}$ for $P\in U_{1}$
and $m>0$, where $Y'$ is the scheme-theoretic preimage $\pi
^{-1}Y(m',r',k')$ in $Z'$, that is, $Y' = V(J')$ where $J' =
I(m',r',k')S$.


As in the proof of Lemma \ref{l:poly/V_m>Y(m,r,k)}, if $P$ is a point
of $U_{1}\cap Z(n,L)$ belonging to the unique component $W_{f}$, then
$Z'$ contains $|T\cap [r]\setminus f([k-1])|$ distinct points lying
over $P$, where $T = \{j:P\in V(x_{j}) \}$.  If $P$ belongs to
$U_{1}\cap Y(m,r,k-1)$, then (since $W_{f}$ is unique) we have $|T\cap
[r]\setminus f([k-1])|\geq m$.  If $P$ also belongs to $Y(m',r',k')$
then the points lying over $P$ belong to $Y'$, which shows that
$V_{m}(Y'_{\text{red}})\cap U_{1}\supseteq Y(m',r',k')\cap
Y(m,r,k-1)\cap U_{1}$ as sets.  By Corollary
\ref{c:poly/Y-reduced-U1}, the last intersection is
scheme-theoretically reduced, so we have $I_{m}(Y'_{\text{red}})_{P}
\subseteq I(m',r',k')_{P} + I(m,r,k-1)_{P}$ for all $P\in U_{1}$.

Let $J\subseteq R(n,L)[b_{k}]$ be the defining ideal of $Z'$, so $S =
R(n,L)[b_{k}]/J$.  We defined $J'$ as an ideal in $S$, but we can also
regard it as an ideal in $R(n,L)[b_{k}]$ which contains $J$.
Obviously $J'$ is contained in the defining ideal $\rad{J'}$ of
$Y'_{\text{red}}$.  This implies $I_{m}(Y')\subseteq
I_{m}(Y'_{\text{red}})$ and hence $I_{m}(Y')_{P}\subseteq
I(m',r',k')_{P}+I(m,r,k-1)_{P}$.  For the reverse containment, we
clearly have $I_{m}(Y')\supseteq I_{1}(Y') = J'\cap R(n,L)\supseteq
I(m',r',k')$ for $m>0$, and we have $I_{m}(Y')_{P}\supseteq
I(m,r,k-1)_{P}$ because $I(m,r,k-1)_{P} = I_{m}(Z')_{P}$ and $J'$
contains $J$.
\end{proof}


\begin{lem}\label{l:poly/reduced-preimage}
Under the hypotheses of Lemma \ref{l:poly/Z'}, each ideal
$I(m',r',k')S$ is equal to its radical, so $S/I(m',r',k')S$ is the
coordinate ring of the reduced preimage $\pi ^{-1}Y(m',r',k')$
in $Z'$, and this coordinate ring is a free $k[\yy ]$-module.
\end{lem}

\begin{proof}
Lemma \ref{l:poly/Z'} implies that $S/I(m',r',k')S$ is a free $k[\yy
]$-module and therefore a subring of $k(\yy )\otimes S/I(m',r',k')S$.
Let $J$ be the ideal of $Z'$ as a subscheme of $Z(n,l)$, so $S \cong
R(n,l)/J$.  Both $J$ and $\rad{I(m',r',k')R(n,l)}$ belong to the
lattice in Lemma \ref{l:poly/Y-reduced-U1}.  Moreover
$I(m',r',k')R(n,l)_{P} = \rad{I(m',r',k')R(n,l)}_{P}$ for all $P\in
U_{1}$, as follows easily from Lemma \ref{l:poly/U_1-local-pix}.
Hence Lemma \ref{l:poly/Y-reduced-U1} implies that
$I(m',r',k')R(n,l)_{P}+J_{P}$ is equal to its radical for all $P\in
U_{1}$, and by Lemma \ref{l:poly/hat(U)-vs-U} this implies that
$(S/I(m',r',k')S)_{Q}$ is a reduced ring for all $Q\in \hat{U}_{1}$.
In particular, taking $Q=0$, $k(\yy )\otimes S/I(m',r',k')S$ is a
reduced ring.
\end{proof}


\begin{lem}\label{l:poly/Z''}
Let $Z' = Z'(r,k)$ and $Z''=Z'(r-1,k)$ be special arrangements over $Z
= Z(n,L)$, where $L = [l]\setminus k$.  Let $Z' = \Spec S$, and let
$J'\subseteq S$ be the ideal of $Z''$ as a closed subscheme of $Z'$.
Assume that $R(n,L)$ is a free $k[\yy ]$-module with a basis $B$ such
that each ideal $I(m',r',k')$ is spanned by a subset of $B$.  Then
$J'$ is a free $k[\yy ]$-module with a basis $B'$ such that every
ideal $I(m',r',k')S\cap J'$ is spanned by a subset of $B'$.
\end{lem}

\begin{proof}
For each basis element $p\in B$, let $m(p)$ be the unique non-negative
integer $m$ such that $p$ belongs to $I(m+1,r,k-1)$ but not to
$I(m,r,k-1)$.  Similarly, let $m'(p)$ be the unique $m'$ such that $p$
belongs to $I(m'+1,r-1,k-1)$ but not to $I(m',r-1,k-1)$.  Thus by
Lemma \ref{l:poly/Z'}, the set of elements
\begin{equation}\label{e:poly/p b^j}
B_{1} = \{pb_{k}^{j}: p\in B,\; j<m(p) \}
\end{equation}
is a free $k[\yy ]$-module basis of $S$ and the set
\begin{equation}\label{e:poly/p b^j'}
B'_{1} = \{pb_{k}^{j}: p\in B,\; j<m'(p) \}
\end{equation}
is a  free $k[\yy ]$-module basis  of $S/J'$.  From  the definition we
have   $Y(m+1,r-1,k-1)\subseteq  Y(m+1,r,k-1)\subseteq  Y(m,r-1,k-1)$,
hence  $I(m,r-1,k-1)\subseteq  I(m+1,r,k-1)\subseteq  I(m+1,r-1,k-1)$.
Therefore we have $m'(p)\leq m(p)\leq m'(p)+1$ for all $p\in B$.

For each $p\in B$ such that $m(p) = m'(p)+1$, the element $p_{0} =
pb_{k}^{m(p)-1}$ belongs to $B_{1}$ but not to $B'_{1}$.  There is,
however, a unique $k[\yy ]$-linear combination $p_{1}$ of the elements
in $B'_{1}$ such that $p_{0}\equiv p_{1} \pmod{J'}$.  Set $q =
p_{0}-p_{1}$, so we have $q\in J'$, and $q$ is congruent to $p_{0}$
modulo $k[\yy ]B_{1}'$.  Let $B'$ be the set of all such elements $q$.

The union $B'\cup B_{1}'$ is again a basis of $S$, while $B'_{1}$ is a
basis of $S/J'$ and $B'$ is contained in $J'$.  This implies that $B'$
spans $J'$, that is, $J'$ is a free $k[\yy ]$-module with basis $B'$.


By Lemma \ref{l:poly/Z'}, $I(m',r',k')S$ is spanned by the subset
$A_{1} = \{pb_{k}^{j}\in B_{1} : p\in I(m',r',k') \}$.  Consider an
element $p\in I(m',r',k')\cap B$ with $p_{0} = pb_{k}^{m(p)-1}\in
B_{1}\setminus B'_{1}$, and the corresponding element $q$ constructed
as above.  By Lemma \ref{l:poly/Z'} applied to $Z''$, the set $A'_{1}
= \{p'b_{k}^{j}\in B'_{1} : p'\in I(m',r',k') \}$ spans
$I(m',r',k')S/(I(m',r',k')S\cap J')$.  Hence there is a $k[\yy
]$-linear combination of the elements of $A'_{1}$ congruent to $p_{0}$
modulo $I(m',r',k')S\cap J'$.  By uniqueness this linear combination
must be $p_{1}$, so $q$ belongs to $I(m',r',k')S\cap J'$.  Let
$A'\subseteq B'$ be the set of these elements $q$.  As before, $A'\cup
A'_{1}$ spans $I(m',r',k')S$, while $A'_{1}$ is a basis of
$I(m',r',k')S/(I(m',r',k')S\cap J')$ and $A'$ is contained in
$I(m',r',k')S\cap J'$.  This implies that $A'$ spans $I(m',r',k')S\cap
J'$.
\end{proof}


\subsection{Basis construction}
\label{ss:poly/basis}

We are now ready to establish the lemmas which we will use inductively
to construct the basis required by Theorem~\ref{t:common-basis}.  The
induction proceeds in stages of three types.  These involve the
construction of a basis of $R(n,l)/I(1,1,l)$ from a basis of
$R(n-1,l)$, the construction of a basis of $R(n,l)/I(1,1,t-1)$ from a
basis of $R(n,l)/I(1,1,t)$, and the construction of a basis of
$R(n,l)$ from a basis of $R(n,l)/I(1,1,0)$.  Because we work with the
rings $R(n,l)/I(1,1,t)$ at intermediate stages, we need to define the
type of basis we want in these rings, as well as in $R(n,l)$.

\begin{defn}\label{d:poly/common-basis}
A {\it common ideal basis} of $R(n,l)$ is a set $B$ for which the
conclusion of Theorem~\ref{t:common-basis} holds: namely, $R(n,l)$ is
a free $k[\yy ]$-module with basis $B$ and every ideal $I(m,r,k)$ is
spanned by a subset of $B$.

Given $t\in \{0 \}\cup [l]$, let $J=I(1,1,t)$.  A {\it common ideal
basis} of $R(n,l)/J$ is a set $B$ such that $R(n,l)/J$ is a free
$k[\yy ]$-module with basis $B$ and every ideal
\begin{equation}\label{e:poly/ideal-mod-J}
\rad{(I(m,r,k)+J)}/J
\end{equation}
is spanned as a $k[\yy ]$-module by a subset of $B$.
\end{defn}

For the sake of clarity let us point out that if a common ideal basis
of $R(n,l)$ exists, then $R(n,l)/(I(m,r,k)+J)$ is a free $k[\yy
]$-module and hence a subring of $k(\yy )\otimes R(n,l)/(I(m,r,k)+J)$.
By Corollary \ref{c:poly/Y-reduced-U1}, the latter ring is reduced, so
$I(m,r,k)+J = \rad{(I(m,r,k)+J)}$.  In light of
Theorem~\ref{t:common-basis}, the radical sign in
\eqref{e:poly/ideal-mod-J} is therefore superfluous.  However, in the
construction of our basis in stages, we will have to construct common
ideal bases for the rings $R(n,l)/J$ first, before we reach $R(n,l)$
itself.  For this reason we have taken care to define a common ideal
basis of $ R(n,l)/J = \Ocal (Y(1,1,t))$ to be a common basis of the
ideals of the {\it reduced} arrangements $Y(m,r,k)\cap Y(1,1,t)$,
without assuming in advance that these intersections are
scheme-theoretically reduced.


Each stage of our basis construction procedure involves the action of
the cyclic shift $\theta $ defined in \eqref{e:poly/theta}.  Let us
briefly note what this action means in various contexts.
Geometrically, $\theta $ acts on $E^{n}$ as the right circular shift
\begin{equation}\label{e:poly/theta-on-E^n}
\theta (P_{1},\ldots,P_{n}) = (P_{n},P_{1},\ldots,P_{n-1}),
\end{equation}
and on $E^{l}$ as the identity.  In coordinates,
$\theta $ acts by
\begin{equation}\label{e:poly/theta-on-polynomials}
\theta x_{i} = x_{i+1},\; \theta y_{i} = y_{i+1}\;(i<n), \quad \theta
x_{n} = x_{1},\; \theta y_{n} = y_{1},
\end{equation}
and fixes the coordinates $\aa ,\bb $.  The shift $\theta $ also acts
on various indexing data, {\it e.g.}, for subsets $S\subseteq [n]$ we
have $\theta S = \theta (S)$, and for functions $f:[l]\rightarrow [n]$
we have $\theta f = \theta \circ f$.  The various actions are mutually
consistent: we have $\theta W_{f} = W_{\theta f}$, for example, and
the ideal of $\theta Y(m,r,k)$ is $\theta I(m,r,k)$.

The ring homomorphism defined by \eqref{e:poly/theta-on-polynomials}
actually corresponds to the {\it inverse} of the morphism $\theta $.
However, it is this definition that makes the notation consistent, as
the reader can verify.  (This is a special case of a general
phenomenon.  If a group $G$ acts on a scheme $X = \Spec R$, then $g\in
G$ should be defined to act on $R$ by the ring homomorphism
corresponding to the morphism $g^{-1}$.  Otherwise the left action of
$G$ on $X$ would become a right action on $R$.)


We turn now to the actual stages of the basis construction procedure,
beginning with the easiest.  Throughout we will maintain the
hypothesis that the bases we construct consist of homogeneous
polynomials.  Strictly speaking, we only need to maintain homogeneity
with respect to the $x$-degree, although it is not hard to see on closer
examination that everything in the construction also maintains
homogeneity with respect to the $y$-degree.

\begin{lem}\label{l:poly/build-Y(1,1,l)}
Suppose that $R(n-1,l)$ has a homogeneous common ideal basis.  Then so
does $R(n,l)/I(1,1,l)$.
\end{lem}

\begin{proof}
We can assume we are given a common ideal basis of $R(n-1,l)$,
represented by a set of homogeneous polynomials $B'\subseteq k[\xx
,\yy ,\aa ,\bb ]$ not involving the variables $x_{n}$, $y_{n}$.  We
will show that $B = \theta B'$ is the required basis of
$R(n,l)/I(1,1,l)$.

Let $J = I(1,1,l)$.  From the definition it is immediate that
$Y(1,1,l) \cong Z(N,l)\times \Spec k[y_{1}]$, where $N = [n]\setminus
\{1 \} = \theta [n-1]$.  This implies that $R(n,l)/J$ is a free
$k[\yy ]$-module with basis $B$.  It follows from Lemma
\ref{l:poly/Y-intersect-set} that for $r>0$, the image of
$Y(m,r,k)\cap Y(1,1,l)$ under the isomorphism is
$Y_{N,l}(m-1,r-1,k)\times \Spec k[y_{1}]$, so $\rad{(I(m,r,k)+J)}/J =
I_{N,l}(m-1,r-1,k)\otimes k[y_{1}]$.  This shows that $B$ is a common
ideal basis as far as the ideals $\rad{(I(m,r,k)+J)}/J$ with $r>0$ are
concerned, and those with $r=0$ are trivial.
\end{proof}


The subsequent stages in the basis construction procedure require the
following technical lemma.

\begin{lem}\label{l:poly/J-and-I's}
Let $R$ be a  $k[\yy ]$-algebra and let $J$, $I_{1},\ldots, I_{m}$ be
ideals of $R$ such that $R/I_{i}$ is a torsion-free $k[\yy
]$-module for all $i$, and for every intersection $I$ of some of the
ideals $I_{i}$, we have $(I+J)_{P} = \rad{(I+J)}_{P}$ for all $P\in
U_{2}$.  Assume that $J$ is a free $k[\yy ]$-module with basis
$B''$ such that each ideal $I_{i}\cap J$ is spanned by a subset of
$B''$, and that $R/J$ is a free $k[\yy ]$-module with basis $B'$ such
that each ideal $\rad{(I_{i}+J)}/J$ is spanned by a subset of $B'$.
Then $R$ is a free $k[\yy ]$-module with a basis $B$ such that $J$ and
all the ideals $I_{i}$ are spanned by subsets of $B$.
\end{lem}

\begin{proof}
Choosing an arbitrary representative $q\in R$ of each $J$-coset
belonging to $B'$, we can assume that $B'$ is given as a subset of
$R$.  Then $R$ is a free $k[\yy ]$-module with basis $B'\cup B''$.
This basis need not be a common basis of the ideals $I_{i}$; to
construct $B$ we have to modify it.

Given $q\in B'$, let $\alpha $ be the set of indices $i$ for which
$\rad{(I_{i}+J)}$ contains $q$, and let $I = \bigcap _{i\in \alpha
}I_{i}$.  Then $\rad{(I+J)} = \bigcap _{i\in \alpha }\rad{(I_{i}+J)}$
contains $q$.  For $Q\in \hat{U}_{2}$ we have $(I+J)_{Q} =
\rad{(I+J)}_{Q}$ by Lemma \ref{l:poly/hat(U)-vs-U}, hence $q\in
(I+J)_{Q}$.  Thus there exists an element $q_{1}\in R_{Q}$ satisfying
\begin{equation}\label{e:poly/q1-conditions}
q_{1}\in I_{Q},\quad q_{1}- q \in J_{Q}.
\end{equation}
These conditions determine $q_{1}$ and $q_{1}-q$ uniquely modulo
$(I\cap J)_{Q}$.  Now, it follows from the hypotheses of the Lemma
that if $C\subseteq B''$ is the subset consisting of elements not
belonging to $I\cap J$, then $J/I\cap J$ is a free $k[\yy ]$-module
with basis $C$.  Expanding the unique element $q_{1}-q\in J_{Q}/(I\cap
J)_{Q}$ satisfying \eqref{e:poly/q1-conditions} in terms of the basis
$C$, we have
\begin{equation}\label{e:poly/q1-unique-expression}
q_{1} = q+\sum _{\gamma} p_{\gamma }b_{\gamma },
\end{equation} 
with $p_{\gamma }\in k[\yy ]_{Q}$ and $b_{\gamma }\in C$.  The
coefficients $p_{\gamma }$ are uniquely determined by
\eqref{e:poly/q1-conditions} and \eqref{e:poly/q1-unique-expression},
and the solution for any $Q$ is also the solution for $Q = 0$.  As in
the proof of Lemma \ref{l:poly/basis-on-U2}, it follows that as
elements of $k(\yy )$, the coefficients $p_{\gamma }\in k[\yy ]_{Q}$
do not depend on $Q$ and hence belong to $k[\yy ]$.  In particular,
the polynomial $q_{1}$ defined in \eqref{e:poly/q1-unique-expression}
is an element of $R$ satisfying \eqref{e:poly/q1-conditions} for all
$Q\in \hat{U_{2}}$.


We now define $B$ to be the union of $B''$ and the set of elements
$q_{1}$ constructed as above, for all $q\in B'$.  By
\eqref{e:poly/q1-unique-expression} we have $q_{1} \equiv q \pmod{J}$.
This implies that the elements $q_{1}$ form a basis of $R/J$, so $B$
is a basis of $R$.

The hypothesis that $R/I_{i}$ is a torsion-free $k[\yy ]$-module means
that the canonical map $R/I_{i}\rightarrow k(\yy )\otimes (R/I_{i})$
is injective.  Since $k(\yy )\otimes (R/I_{i}) = (R/I_{i})_{Q}$ for
$Q=0$, and we have $q_{1}\in I_{(0)}\subseteq (I_{i})_{(0)}$, this
implies $q_{1}\in I_{i}$ for all $i\in \alpha $.

For each $i$, the elements $q\in B'\cap \rad{(I_{i}+J)}$ span the
$k[\yy ]$-module $\rad{(I_{i}+J)}/J$ by hypothesis, and therefore span
the vector space $k(\yy )\otimes \rad{(I_{i}+J)}/J = k(\yy )\otimes
((I_{i}+J)/J)$.  Since the corresponding elements $q_{1}$ belong to
$I_{i}$ and satisfy $q_{1}\equiv q \pmod{J}$, it follows from the
canonical isomorphism $I_{i}/(I_{i}\cap J)\rightarrow (I_{i}+J)/J$
that they span $k(\yy )\otimes (I_{i}/(I_{i}\cap J))$.  By hypothesis
the elements of $B''\cap I_{i}$ span $I_{i}\cap J$, so the elements of
$B\cap I_{i}$ all together span $k(\yy )\otimes I_{i}$.  By Lemma
\ref{l:poly/basis-of-submodule}, this implies they span $I_{i}$ as a
$k[\yy ]$-module.
\end{proof}

To maintain homogeneity, we need a graded version of the preceding
lemma.  For the stage going from $R(n,l)/I(1,1,0)$ to $R(n,l)$ in our
basis construction procedure we will want a bit more: the graded
version should hold in each degree separately.  To avoid cluttering
the notation, we did not mention the grading in the proof of Lemma
\ref{l:poly/J-and-I's}, but everything we need is already implicit
there.  To see this, consider the situation in which $R = \bigoplus
_{d}R_{d}$ is a graded $k[\yy ]$-algebra with $k[\yy ]\subseteq
R_{0}$, the ideals $I_{i}$ and $J$ are homogeneous, and the elements
of $B'$ and $B''$ are also homogeneous.  Then the construction of the
basis $B$ in the proof of the lemma can be carried out in each degree
separately.  More precisely, the proof yields a common basis for the
homogeneous components of the ideals $J$ and $I_{i}$ in a given degree
$d$, assuming only that we have bases $B''_{d}$ of $J_{d}$ and
$B'_{d}$ of $(R/J)_{d}$ satisfying the hypotheses of the lemma in
degree $d$. Thus we have the following corollary to the proof of Lemma
\ref{l:poly/J-and-I's}.


\begin{cor}\label{c:poly/J-and-I's}
Let $R$, $J$ and $I_{1},\ldots,I_{m}$ be as in Lemma
\ref{l:poly/J-and-I's}, but assume further that $R = \bigoplus
_{d}R_{d}$ is graded, with $k[\yy ]\subseteq R_{0}$, and the ideals
$J$, $I_{1},\ldots, I_{m}$ are homogeneous.  For a given degree $d$,
suppose that $J_{d}$ is a free $k[\yy ]$-module with basis $B''_{d}$
such that each $(I_{i}\cap J)_{d}$ is spanned by a subset of
$B''_{d}$, and that $(R/J)_{d}$ is a free $k[\yy ]$-module with basis
$B'_{d}$ such that each $(\rad{(I_{i}+J)}/J)_{d}$ is spanned by a
subset of $B'_{d}$.  Then $R_{d}$ is a free $k[\yy ]$-module with a
basis $B_{d}$ such that $J_{d}$ and all $(I_{i})_{d}$ are spanned by
subsets of $B_{d}$.
\end{cor}

We remark that the torsion-freeness hypothesis and the local
reducedness condition $(I+J)_{P} = \rad{(I+J)}_{P}$ in Lemma
\ref{l:poly/J-and-I's} and Corollary \ref{c:poly/J-and-I's} are
automatically satisfied in the situations where we will apply them.
The ideals $I_{i}$ will always be ideals of (reduced) subspace
arrangements $Y(m,r,k)$ or $Y(m,r,k)\cap Y(1,1,t)$.  By Lemmas
\ref{l:poly/torsion-free} and \ref{l:poly/Y-intersect-set} their
coordinate rings $R/I_{i}$ are torsion-free $k[\yy ]$-modules.  The
ideal $J$ will always be $I(1,1,t)$ for some $t$, so it follows from
Lemma \ref{l:poly/Y-reduced-U2} that $V(I+J)$ is locally reduced on
$U_{2}$. 

It follows {\it a posteriori} from the conclusion of Lemma
\ref{l:poly/J-and-I's} that if the reducedness hypotheses $(I+J)_{P} =
\rad{(I+J)}_{P}$ hold locally on $U_{2}$, we must actually have $(I+J)
= \rad{(I+J)}$.  This is similar to the situation we encountered
earlier with Lemma \ref{l:poly/I_m-basis}.  As in that case, we proved
the lemma with the weaker local hypotheses so that we can verify them
by reduction to the case $n=2$.


We come now to the next stage in the basis construction procedure.
This is actually the third stage in the inductive process, but we
present it before the second stage because for the latter we still
need some further geometric results.

\begin{lem}\label{l:poly/build-Z(n,l)}
Suppose that $R(n,l)/I(1,1,0)$ has a homogeneous common ideal basis.
Then so does $R(n,l)$.
\end{lem}

\begin{proof}
Let $J = I(1,1,0)$.  By definition, $Y(1,1,0) = V(x_{1})$, so $J =
\rad{(x_{1})}$.  By Lemma \ref{l:poly/translation}, $\rad{(x_{1})} =
(x_{1})$.  Since $x_{1}$ does not vanish identically on any $W_{f}$,
it is a non-zero-divisor in $R(n,l)$.  Hence multiplication by $x_{1}$
is an isomorphism of $R(n,l)$ onto $J = (x_{1})$.

We may assume we are given a homogeneous common ideal basis $B'$ of
$R(n,l)/J$.  Suppose that in a given $x$-degree $d$, we can find a
free $k[\yy ]$-module basis $B_{d}$ of $R(n,l)_{d}$ such that every
$I(m,r,k)_{d}$ is spanned by a subset of $B_{d}$.  We claim that
$x_{1}\theta B_{d}$ is then a basis of $J_{d+1}$ with subsets spanning
each $(I(m,r,k)\cap J)_{d+1}$.  Granting this claim for the moment, we
can apply Corollary \ref{c:poly/J-and-I's} with $B''_{d+1} =
x_{1}\theta B_{d}$ and $B'_{d+1}$ the degree $d+1$ part of the given
basis $B'$ to obtain a free $k[\yy ]$-module basis $B_{d+1}$ of
$R(n,l)_{d+1}$ which is a common basis for every $I(m,r,k)_{d+1}$.  In
degree zero, we can take $B_{0} = B'_{0}$, since $R(n,l)_{0} =
(R(n,l)/J)_{0}$.  We can then construct $B_{d}$ by induction for all
degrees $d$, obtaining a common ideal basis $B = \bigcup _{d}B_{d}$ of
$R(n,l)$. 

It remains to prove the claim.  If $B_{d}$ is a free module basis of
$R(n,l)_{d}$ then so is $\theta B_{d}$, and therefore $x_{1}\theta
B_{d}$ is a basis of $J_{d+1}$, since multiplication by $x_{1}$ is an
isomorphism of $R(n,l)_{d}$ onto $J_{d+1}$.  To complete the proof,
observe that for any ideal $I\subseteq R(n,l)$ we have $I\cap J =
I\cap (x_{1}) = x_{1}(I\icolon (x_{1}))$.  If $I$ is radical, then so
is $I\icolon (x_{1})$, and $V(I\icolon (x_{1}))$ is the union of those
components of $V(I)$ on which $x_{1}$ does not vanish identically.
Applying this to $I=I(m,r,k)$ for $r>0$, we get $I(m,r,k)\icolon
(x_{1}) = \theta I(m,r-1,k)$, and hence $I(m,r,k)\cap J = x_{1}\theta
I(m,r-1,k)$.  For $r>0$, this shows that if $B_{d}$ has a
subset spanning $I(m,r-1,k)_{d}$, then $x_{1}\theta B_{d}$ has a
subset spanning $(I(m,r,k)\cap J)_{d+1}$.  This suffices, since
$I(m,r,k)$ is trivially equal to $0$ or $(1)$ for $r=0$.
\end{proof}


The remaining and most subtle stage of the basis construction
procedure is the one going from $R(n,l)/I(1,1,t)$ to
$R(n,l)/I(1,1,t-1)$, for $t\in [l]$.  For this stage we will also
assume we have a common ideal basis of $R(n,l-1)$, and apply the basis
lifting technique developed in \S \S \ref{ss:poly/lifting-I} and
\ref{ss:poly/lifting-II}.  Recall that the basis lifting involves the
special arrangements $Z'(r,k)$.  In order to make use of them we need
a lemma relating the ideal of $Y(1,1,t)$ as a closed subscheme of
$Y(1,1,t-1)$ to the ideal of $Z''$ as a closed subscheme of $Z'$, for
a suitable pair of special arrangements $Z''\subseteq Z'$.  Before
stating the lemma we fix some notation.  Let $n>1$, $l>0$ and $t\in
[l]$ be given.  Let $Z' = Z'(n,t)$ be the special arrangement with $r
= n$ and $k=t$, and let $Z'' = Z'(n-1,t)$.

Note that $Z''$ is a subarrangement of $Z'$.  The components of $Z'$
that do not belong to $Z''$ are the subspaces
\begin{equation}\label{e:poly/Z'-not-Z''-components}
V(a_{t})\cap W_{f} = V(x_{n})\cap W_{f}
\end{equation}
for $f:[l]\rightarrow [n]$ such that $f(t) = n$, $f(i)\not =n$ for
$i\in [t-1]$.  Let $Z_{1}$ be the union of these components.  Then we
have 
\begin{equation}\label{e:poly/Z'=Z''-U-Z0}
Z' = Z''\cup Z_{1},
\end{equation}
with $Z''$ and $Z_{1}$ having no component in common.

Similarly, $Y(1,1,t)$ is a subarrangement of $Y(1,1,t-1)$, and the
components of $Y(1,1,t-1)$ not belonging to $Y(1,1,t)$ are the subspaces
\begin{equation}\label{e:poly/Y(t-1)-not-Y(t)-components}
V(x_{1})\cap W_{f}
\end{equation}
for $f:[l]\rightarrow [n]$ such that $f(t)=1$, $f(i)\not =1$ for $i\in
[t-1]$.  Letting $Z_{0}$ be their union, 
we have as above
\begin{equation}\label{e:poly/Y(t-1)=Y(t) U theta-Z0}
Y(1,1,t-1) = Y(1,1,t)\cup  Z_{0},
\end{equation}
with $Y(1,1,t)$ and $Z_{0}$ having no component in common.  (Not
accidentally, this is the same $Z_{0}$ that appeared in the proof of
Lemma \ref{l:poly/Y-vs-Z0}.)


Observe that $Z_{0} = \theta Z_{1}$.  The key fact about this set-up
is as follows.

\begin{lem}\label{l:poly/scheme-intersection}
With $n$, $l$, $t$, $Z'$, $Z''$ and $Z_{0}$ as above, we have 
\begin{equation}\label{e:poly/scheme-intersection}
\theta Z''\cap Z_{0} = Y(1,1,t)\cap Z_{0},
\end{equation}
scheme-theoretically on $U_{2}$.  In other words, for all $P\in U_{2}$
we have equality of ideals $\theta I(Z'')_{P}+I(Z_{0})_{P} =
I(1,1,t)_{P}+ I(Z_{0})_{P}$.
\end{lem}

\begin{proof}
We caution the reader immediately that the intersections in
\eqref{e:poly/scheme-intersection} are {\it not} scheme-theoretically
reduced in general, so it is not enough to check the result
set-theoretically.  Instead, we must use our knowledge of the local
picture on $U_{2}$ to write down equations.  If $P$ is a point of the
intersection on either side of \eqref{e:poly/scheme-intersection},
then we have $P\in W_{f}\cap W_{g}$ for some $f$, $g$ with $f(t)\not
=1$, $g(t)=1$, so $P$ is not in $U_{1}$.  Thus we need only consider
points $P\in U_{2}\setminus U_{1}$.

The schemes in question are subschemes of $Z(n,l)$.  Fix a point $P\in
U_{2}\setminus U_{1}$, and let $N$, $L$, $\sim $, $F$, and $h\in F$ be
as in Lemma \ref{l:poly/U_2-local-pix}.  We may replace $Z(n,l)$ with
$E^{[n]\setminus N}\times Z(N,L)$ without changing any of the local
ideals at $P$.

First consider the case $1\not \in N$.  If $P\in Z_{0}$ we must have
$h(t)=1$.  But if $P\in \theta Z''$ or $P\in Y(1,1,t)$, we must have
$h(t)\not =1$.  Hence both intersections in
\eqref{e:poly/scheme-intersection} are locally empty at $P$.

For the case $1\in N$ we may assume $h(t)\in N$, and hence $t\in L$,
as otherwise $P\not \in Z_{0}$ and the result is trivial.  Under the
local isomorphism of $Z(n,l)$ with $E^{[n]\setminus N}\times Z(N,L)$,
$Z_{0}$ coincides locally with $E^{[n]\setminus N}\times
\tilde{Z}_{0}$, where $\tilde{Z}_{0}\subseteq Z(N,L)$ is the subspace
arrangement
\begin{equation}\label{e:poly/Z-tild_0}
\tilde{Z}_{0} = \bigcup _{f}V(x_{1})\cap W_{f},
\end{equation}
over $f:L\rightarrow N$ such that $f(t)=1$, $f(i)\not =1$ for $i\in
[t-1]\cap L$.


Similarly, $Y(1,1,t)$ coincides locally with $E^{[n]\setminus N}\times
Y_{N,L}(1,1,s)$, where $s=|[t]\cap L|$, and $\theta
Z''$ coincides with $E^{[n]\setminus N}\times \tilde{Z}$, where
\begin{equation}\label{e:poly/Z-tild}
\tilde{Z} = \bigcup _{f}V(a_{t})\cap W_{f},
\end{equation}
over $f:L\rightarrow N$ such that $f(t)\not =1$, $f(i)=1$ for $i\in
[t-1]\cap L$.

We are to show that $\tilde{Z}\cap \tilde{Z}_{0} =
Y_{N,L}(1,1,s)\cap \tilde{Z}_{0}$ as subschemes of $Z(N,L)$.
Renaming the indices so that $N$ becomes $[2]$, $L$ becomes
$[\tilde{l}\,]$, and (hence) $t$ becomes $s$, our subschemes
become subschemes of $Z(2,\tilde{l}\,)$:
\begin{equation}\label{e:poly/Z's-in-Z(2,l)}
\begin{aligned}
\tilde{Z}_{0}&	 = \bigcup V(x_{1})\cap W_{f}:\quad f(s)=1,\;
\text{$f(i) = 2$ for $i<s$},\\
\tilde{Z}&	 = \bigcup V(x_{2})\cap W_{f}:\quad f(s)=2,\;
\text{$f(i) = 1$ for $i<s$},\\
Y(1,1,s)&	 = \bigcup V(x_{1})\cap W_{f}:\quad f(s)=2,\;
\text{$f(i) = 2$ for $i<s$}.
\end{aligned}
\end{equation}


Let $I\subseteq R(2,\tilde{l}\,)$ be given by
\begin{equation}\label{e:poly/I-for-Z'-schem}
I = (a_{s}-x_{1},b_{s}-y_{1})+\sum _{i<s}
(a_{i}-x_{2},b_{i}-y_{2}).
\end{equation}
By Lemma \ref{l:poly/elimination}, $I=\rad{I}$ and $V(I)\cong
Z(2,\tilde{L})$, where $\tilde{L} = [\tilde{l}]\setminus [s]$.  By
Corollary \ref{c:poly/n=2-I(m,r,k)-gens} or Lemma
\ref{l:poly/translation}, we have $I(1,1,0) = (x_{1})$ in
$Z(2,\tilde{L})$, which shows that $(x_{1})+I$ is a radical ideal.
Since we clearly have $V((x_{1})+I) = \tilde{Z}_{0}$
set-theoretically, it follows that
\begin{equation}\label{e:poly/I(Z-tild-0)}
I(\tilde{Z}_{0}) =
(x_{1})+(a_{s}-x_{1},b_{s}-y_{1})+\sum _{i<s}
(a_{i}-x_{2},b_{i}-y_{2}).
\end{equation}
By symmetry, we have
\begin{equation}\label{e:poly/I(Z-tild'')}
I(\tilde{Z}) = (x_{2})+(a_{s}-x_{2},b_{s}-y_{2})+\sum _{i<s}
(a_{i}-x_{1},b_{i}-y_{1}),
\end{equation}
and by Corollary \ref{c:poly/n=2-I(m,r,k)-gens} we have
\begin{equation}\label{e:poly/I(1,1,s)}
I(1,1,s) =
(x_{1})+(a_{s}-x_{2},b_{s}-y_{2})+\sum _{i<s}
(a_{i}-x_{2},b_{i}-y_{2}).
\end{equation}
Now we see immediately that both $I(\tilde{Z})+I(\tilde{Z}_{0})$ and
$I(1,1,s)+I(\tilde{Z}_{0})$ contain $(x_{1}-x_{2},y_{1}-y_{2})$ and
hence they both reduce to
\begin{equation}\label{e:poly/both-ideals}
(x_{1},x_{2},y_{1}-y_{2})+\sum _{i\leq s} (a_{i},b_{i}-y_{1}).
\end{equation}
\end{proof}

As was the case with Lemma \ref{l:poly/Y-reduced-U2}, one can show
using Theorem~\ref{t:common-basis} that the conclusion of Lemma
\ref{l:poly/scheme-intersection} actually holds everywhere, and not
just on $U_{2}$.  As before, the result on $U_{2}$ suffices for our
purposes, and the restriction to $U_{2}$ enables us to prove it by
reduction to the case $n=2$.


\begin{lem}\label{l:poly/build-Y(1,1,t-1)} Given $n>1$, $l>0$ and
$t\in [l]$, suppose that $R(n,l-1)$ and $R(n,l)/I(1,1,t)$ each have a
homogeneous common ideal basis.  Then so does $R(n,l)/I(1,1,t-1)$.
\end{lem}

\begin{proof}
To simplify notation let $R = R(n,l)/I(1,1,t-1)$ and let $J$ be the
ideal $I(1,1,t)/I(1,1,t-1)$ in $R$.

We first prove that $R$ is a free $k[\yy ]$-module.  The ideal
$(a_{t})\subseteq R$ is isomorphic to $R/(0\icolon (a_{t}))$.  Since
$I(1,1,t-1)\icolon (a_{t}) = I(1,1,t)$ in $R(n,l)$, we have
$R/(0\icolon (a_{t}))\cong R(n,l)/I(1,1,t)$, and the latter is a free
$k[\yy ]$-module by hypothesis.  The locus $V(a_{t})\cap Y(1,1,t-1)$
is isomorphic to the preimage $\pi ^{-1}Y_{n,L}(1,1,t-1)$ in the
special arrangement $Z'(n,1)$ over $Z(n,L)$, where $L = [l]\setminus
\{1 \}$.  The isomorphism is given by transposing the indices $1$ and
$t$ in $[l]$.  Lemma \ref{l:poly/reduced-preimage} then implies that
$R/\rad{(a_{t})}$ is a free $k[\yy ]$-module.  Let $B_{1}$ be a free
$k[\yy ]$-module basis of $(a_{t})$ and let $B_{2}$ be a basis of
$R/\rad{(a_{t})}$.  By Lemmas \ref{l:poly/hat(U)-vs-U} and
\ref{l:poly/I(1,1,t-1)+(a_t)} we have $\rad{(a_{t})}_{Q} =
(a_{t})_{Q}$ for $Q\in \hat{U}_{2}$, so $B_{1}\cup B_{2}$ is a basis
of $R_{Q}$.  Then $R$ is a free $k[\yy ]$-module by Lemma
\ref{l:poly/basis-on-U2}.


Now take $Z'$ , $Z''$, and $Z_{0}$ as in the preamble to Lemma
\ref{l:poly/scheme-intersection}, so $Z'$ and $Z''$ are special
arrangements over $Z(n,L)$, where $L = [l]\setminus \{t \}$.  Set
$R_{0} = \Ocal (Z_{0})$, $R' = \Ocal (Z')$, and let $J'' =
\theta I(Z'')/I(Z')$ be the ideal of $\theta Z''$ as a closed
subscheme of $Z' = \theta Z'$.

By \eqref{e:poly/Y(t-1)=Y(t) U theta-Z0} we have $I(Z_{0})\cap J = 0$
in $R$ so the canonical map $J \rightarrow JR_{0}$ is an isomorphism.
Similarly, applying $\theta $ to \eqref{e:poly/Z'=Z''-U-Z0}, we have
$I(Z_{0})\cap \theta I(Z'')=0$ in $R'$, hence $J''\rightarrow
J''R_{0}$ is an isomorphism as well.

Since $R$ and $R/J$ are both free $k[\yy ]$-modules, finitely
generated and $y$-graded in each $x$-degree, it follows that $J$ is a
free $k[\yy ]$-module.  Hence $JR_{0}$ is a free $k[\yy ]$-module.  By
Lemma \ref{l:poly/Z''}, $I(Z'')/I(Z')$ is a free $k[\yy ]$-module with
a common basis $B$ for the ideals $I\cap I(Z'')/I(Z')$, where
$I=I_{n,L}(m,r,k)R'$.  In particular, this implies that $J''$ and
$J''R_{0}$ are free $k[\yy ]$-modules.  By Lemma
\ref{l:poly/torsion-free}, $R_{0}$ is a torsion-free $k[\yy ]$-module.
By Lemmas \ref{l:poly/hat(U)-vs-U} and
\ref{l:poly/scheme-intersection}, the free submodules $JR_{0}$ and
$J''R_{0}$ of $R_{0}$ coincide when localized at $Q\in \hat{U}_{2}$,
so by Corollary \ref{c:poly/2-free-submodules} they are equal.  Let
$J_{0}=JR_{0} = J''R_{0}$.

With $B$ as in the previous paragraph, $\theta B$ is a free $k[\yy
]$-module basis of the ideal $J''$ in $R'$, with subsets spanning the
ideals $J''\cap \theta I_{n,L}(m,r,k)R'$.  We have canonical
isomorphisms $J''\cong J_{0}\cong J$.  Let $B''\subseteq J$ be the
image of $\theta B$ under the composite isomorphism $J''\cong J$.


The canonical isomorphism $J''\rightarrow J_{0}$ is given by
restriction to $Z_{0}$ of functions on $Z'$ vanishing on $\theta Z''$.
Lemma \ref{l:poly/reduced-preimage} implies that each ideal $\theta
I_{n,L}(m,r,k)R'$ is equal to its radical, namely the ideal $I(\pi
^{-1}\theta Y_{n,L}(m,r,k))$ of the reduced preimage of $\theta
Y_{n,L}(m,r,k)$ in $Z'$.  Therefore, since every function $p\in J''$
vanishes on $\theta Z''$, $p$ belongs to $\theta I_{n,L}(m,r,k)R'$ if
and only if its restriction to $Z_{0}$ vanishes on $Z_{0}\cap \pi
^{-1}\theta Y_{n,L}(m,r,k)$.  Similarly, writing $\pi _{1}$ for the
coordinate projection $Y(1,1,t-1)\subseteq Z(n,l)\rightarrow Z(n,L)$,
a function $p\in J$ belongs to the ideal $J\cap I( \pi _{1}^{-1}\theta
Y_{n,L}(m,r,k))$ in $R$ if and only if the same criterion holds (note
that $\pi $ and $\pi _{1}$ have the same restriction to $Z_{0}$).
This shows that the isomorphism $J''\cong J$ carries $J''\cap \theta
I_{n,L}(m,r,k)R'$ onto $J\cap I( \pi _{1} ^{-1}\theta
Y_{n,L}(m,r,k))$.  Hence the latter ideals are spanned by subsets of
$B''$.

Given $r>0$, $k$ and $m$, set $k' = k$, $m' = m-1$ if $k<t$, or else
$k' = k-1$, $m' = m$ if $k\geq t$.  By Lemma \ref{l:poly/Y-vs-Z0}, we
have $Y(1,1,t)\cup (Y(m,r,k)\cap Y(1,1,t-1)) = Y(1,1,t)\cup \pi_{1}
^{-1}\theta Y_{n,L}(m',r-1,k')$.  The ideal in $R$ of $Y(1,1,t)\cup
\pi_{1} ^{-1}\theta Y_{n,L}(m',r-1,k')$ is $J\cap I( \pi _{1}
^{-1}\theta Y_{n,L}(m',r-1,k'))$, and we have shown above that it is
spanned by a subset of $B''$.  But this ideal is equal to the ideal of
$Y(1,1,t)\cup (Y(m,r,k)\cap Y(1,1,t-1))$, namely, $J\cap
\rad{I(m,r,k)R}$.  This shows that all the ideals $J\cap
\rad{I(m,r,k)R}$ are spanned by subsets of $B''$, the cases with $r =
0$ being trivial.  The conclusion now follows from Corollary
\ref{c:poly/J-and-I's}.
\end{proof}


\subsection{Proof and consequences of Theorem~\ref{t:common-basis}}
\label{ss:poly/proof}

We have established all the stages in our basis construction procedure
in \S \ref{ss:poly/basis}.  To complete the proof of
Theorem~\ref{t:common-basis} and its corollary
Theorem~\ref{t:polygraphs-are-free}, we have only to assemble the
pieces.

\begin{proof}[Proof of Theorem~\ref{t:common-basis}] We prove the
theorem by induction on $n$ and $l$.  The base case $l=0$ is given by
Lemma \ref{l:poly/l=0}.  The base case for $n$ is $n=1$.  Note that
for $n=1$ we have $Z(1,l) \cong Z(1,0) = E$, for all $l$, and that the
only non-trivial $Y(m,r,k)$ is $Y(1,1,0)$, which already appears in
$Z(1,0)$.  Thus the case $n=1$ is included in the case $l=0$.

For $n>1$ and $l>0$ we can assume by induction that $R(n-1,l)$ has a
homogeneous common ideal basis.  Then by Lemma
\ref{l:poly/build-Y(1,1,l)}, so does $R(n,l)/I(1,1,l)$.  We can also
assume by induction that $R(n,l-1)$ has a homogeneous common ideal
basis.  Applying Lemma \ref{l:poly/build-Y(1,1,t-1)} repeatedly, with
$t$ descending from $l$ to $1$, we conclude that $R(n,l)/I(1,1,0)$ has
a homogeneous common ideal basis.  Then by Lemma
\ref{l:poly/build-Z(n,l)}, so does $R(n,l)$.
\end{proof}

The bridge to Hilbert schemes is supplied by the following consequence
of Theorem~\ref{t:polygraphs-are-free}, which we restate here for the
reader's convenience.  We remark that the bridge is not really quite
as narrow as we are making it appear.  The polygraph rings $R(n,l)$
carry geometric information about the tensor powers $\rho
^{*}B^{\otimes l}$ of the tautological bundle over the isospectral
Hilbert scheme $X_{n}$.  Theorem~\ref{t:polygraphs-are-free} can be
further exploited to obtain vanishing theorems for these vector
bundles.  This subject will be taken up elsewhere, along with its
application to the determination of the character formula for diagonal
harmonics, which was discussed briefly in the introduction.


{
\def\theprop{\ref{p:hilbert/J-is-free}}
\begin{prop}\label{p:hilbert/J-is-free-A}
Let $J = \CC [\xx ,\yy ] A$ be the ideal generated by the space of
alternating polynomials $A = \CC [\xx ,\yy ]^{\epsilon }$.  Then
$J^{d}$ is a free $\CC [\yy ]$-module for all $d$.
\end{prop}
}

\begin{proof}
Set $l = nd$, and let $Z(n,l)$ be the polygraph over $\CC $, a
subspace arrangement in $(\CC ^{2})^{n}\times (\CC ^{2})^{l}$.  Let
$G=S_{n}^{d}$ be the Cartesian product of $d$ copies of the symmetric
group $S_{n}$, acting on $(\CC ^{2})^{n}\times (\CC ^{2})^{l}$ by
permuting the factors in $(\CC ^{2})^{l}$ in $d$ consecutive blocks of
length $n$.  In other words each $w \in G$ fixes the coordinates $\xx
,\yy $ on $(\CC ^{2})^{n}$, and for each $k=0,\ldots,d-1$ permutes the
coordinate pairs $a_{kn+1}, b_{kn+1}$ through $a_{kn+n}, b_{kn+n}$
among themselves.

Let $R(n,l) = \CC [\xx ,\yy ,\aa ,\bb ]/I(n,l)$ be the coordinate ring
of $Z(n,l)$.  By Theorem~\ref{t:polygraphs-are-free}, $R(n,l)$ is a
free $\CC [\yy ]$-module.  By the symmetry of its definition, $I(n,l)$
is a $G$-invariant ideal, so $G$ acts on $R(n,l)$.  We claim that
$J^{d}$ is isomorphic as a $\CC [\xx ,\yy ] $-module to the space
$R(n,l)^{\epsilon }$ of $G$-alternating elements of $R(n,l)$.  Each
$x$-degree homogeneous component of $R(n,l)$ is a finitely generated
$y$-graded free $\CC [\yy ]$-module.  Since $R(n,l)^{\epsilon }$ is a
graded direct summand of $R(n,l)$, it is a free $\CC [\yy ]$-module,
so the claim proves the Lemma.


Let $f_{0}\colon [l]\rightarrow [n]$ be defined by $f_{0}(kn+i) = i$
for all $0\leq k < d$, $1\leq i \leq n$.  Restriction of regular
functions from $Z(n,l)$ to its component subspace $W_{f_{0}}$ is given
by the $\CC [\xx ,\yy ]$-algebra homomorphism $\psi :R(n,l)\rightarrow
\CC [\xx ,\yy ]$ mapping $a_{kn+i},b_{kn+i}$ to $x_{i},y_{i}$.
Observe that $\psi $ maps $R(n,l)^{\epsilon }$ surjectively onto
$\CC [\xx ,\yy ] A^{d} = J^{d}$.

Let $p$ be an arbitrary element of $R(n,l)^{\epsilon }$.  Since $p$ is
$G$-alternating, $p$ vanishes on $W_{f}$ if $f(kn+i)=f(kn+j)$ for some
$0\leq k<d$ and some $1\leq i<j\leq n$.  Thus the regular function
defined by $p$ on $Z(n,l)$ is determined by its restriction to those
components $W_{f}$ such that for each $k$, the sequence
$f(kn+1),\ldots,f(kn+n)$ is a permutation of $\{1,\ldots,n \}$.
Moreover, for every such $f$ there is an element $w\in G$ carrying
$W_{f}$ onto $W_{f_{0}}$.  Hence $p$ is determined by its restriction
to $W_{f_{0}}$.  This shows that $p$ vanishes on $Z(n,l)$ if $\psi
(p)=0$, that is, the kernel of the map $\psi \colon R(n,l)^{\epsilon
}\rightarrow J^{d}$ is zero.
\end{proof}


\subsection{Arbitrary ground rings}
\label{ss:poly/ground-ring}

For convenience of exposition we have assumed that $k$ is a field of
characteristic zero.  For completeness we now show that the results on
polygraphs hold over any ground ring.  This generalization is not
needed elsewhere in the paper, where we always have $k=\CC $.

To begin with, observe that the proof of
Theorems~\ref{t:polygraphs-are-free} and \ref{t:common-basis} goes
through with minor modifications for $k=\ZZ $.  As a matter of
notation, we must replace $k(\yy )$ throughout with $\QQ (\yy )$.  In
Lemma \ref{l:poly/numerical-criterion-of-freeness} and elsewhere, we
must read ``generates as a $\ZZ $-module'' for ``spans as a $k$-vector
space.''

In the proofs of Lemmas \ref{l:poly/V(xT) x Y(n=2)-reduced},
\ref{l:poly/I(1,1,t-1)+(a_t)} and \ref{l:poly/I_m=I(m,r,k)-on-U2} we
implicitly used the facts that every scheme over a field $k$ is flat
over $k$, and that the product of reduced schemes over a field of
characteristic zero is reduced.  In our case, the schemes in question
are (reduced) subspace arrangements in $E^{n}\times E^{l}$.  If $Y =
\bigcup Y_{\alpha }$ is a subspace arrangement defined over $\ZZ $,
then its coordinate ring $R = \Ocal (Y)$ is a subring of a direct sum
$\bigoplus \Ocal (Y_{\alpha })$ of polynomial rings over $\ZZ $.
Therefore $R$ is a torsion-free Abelian group, that is, a flat $\ZZ
$-module.  Every subspace arrangement over $\ZZ $ is thus flat, and a
product of reduced schemes flat over $\ZZ $ is reduced.

In the proofs of Lemmas \ref{l:poly/basis-on-U2} and
\ref{l:poly/J-and-I's} we also used the fact that a rational function
regular outside a subset of codimension two in $\Spec k[\yy ]$ is
regular.  This holds with any normal integral domain in place of
$k[\yy ]$, and in particular for $\ZZ [\yy ]$.

The theorem for $k={\ZZ }$ implies the following fully general result.

\begin{thm}\label{t:polygraph-over-any-ring}
Let $E=\AA ^{2}(k)$, where $k$ is any commutative ring with unit.
Let $Z(n,l)\subseteq E^{n}\times E^{l}$ be a polygraph over $k$
(defined as the union of the closed subschemes $W_{f}$, just as in
\ref{d:poly/Z(n,l)}).  Then the coordinate ring $R(n,l)=\Ocal (Z(n,l))$ is a
free $k[\yy ]$-module.
\end{thm}


\begin{proof}
By the preceding remarks, the Theorem holds for $k=\ZZ $.  Tensoring
the exact sequence
\begin{equation}\label{e:poly/ZZ-ideal-sequence}
0\rightarrow I_{\ZZ }(n,l)\rightarrow \ZZ [\xx ,\yy ,\aa ,\bb
]\rightarrow R_{\ZZ }(n,l)\rightarrow 0,
\end{equation}
over $\ZZ [\yy ]$ with $k[\yy ]$, we get an exact sequence
\begin{equation}\label{e:poly/k-ideal-sequence}
0\rightarrow I = k[\yy ] I_{\ZZ }(n,l)\rightarrow k[\xx ,\yy ,\aa
,\bb] \rightarrow R\rightarrow 0
\end{equation}
with $R$ a free $k[\yy ]$-module.  Since $I_{\ZZ }(n,l)\subseteq
I_{\ZZ }(W_{f})$ for all $f$, and $I_{f} = k[\yy ]I_{\ZZ }(W_{f})$
by \eqref{e:poly/I_f}, we have $I\subseteq I(n,l)$.

Let $\delta = \prod _{i<j}(y_{i}-y_{j})$.  Note that $U_{1}$ is the
affine open set $U_{\delta }$, so inverting $\delta $ is the same
thing as localizing to $U_{1}$.  By Lemma \ref{l:poly/U_1-local-pix},
therefore, we have $R_{\ZZ }(n,l)_{\delta } = \bigoplus _{f}\Ocal_{\ZZ }
(W_{f})_{\delta }$, and hence $R_{\delta } = \bigoplus _{f}\Ocal
(W_{f})_{\delta }$.  By definition $I(n,l)$ is the kernel of the
canonical homomorphism
\begin{equation}\label{e:poly/I(n,l)-as-kernel}
k[\xx ,\yy ,\aa ,\bb ]\rightarrow \bigoplus _{f} \Ocal (W_{f}).
\end{equation}
Hence inverting $\delta $ in \eqref{e:poly/k-ideal-sequence} gives
$I_{\delta } = I(n,l)_{\delta }$.

Now $\delta $ is a non-zero-divisor in $k[\yy ]$, so it is also a
non-zero-divisor on the free $k[\yy ]$-module $R$ and its submodule
$I(n,l)/I$.  But since $(I(n,l)/I)_{\delta } = 0$, every element of
$I(n,l)/I$ is annihilated by a power of $\delta $.  Therefore we have
$I(n,l)/I=0$ and $R=R(n,l)$.
\end{proof}


\section{Applications, conjectures, and problems}
\label{s:other}

\subsection{The $G$-Hilbert scheme}
\label{ss:other/G-hilb}

Let $G\subseteq GL(V)$ be a finite group acting faithfully on a
complex vector space $V = \CC ^{m}$.  For all vectors $v$ in a
non-empty open subset of $V$, the orbit $Gv$ consists of $N = |G|$
distinct points and thus represents an element of $\Hilb ^{N}(V)$.
The closure in $\Hilb ^{N}(V)$ of the set of such orbits is a
component of the $G$-fixed locus of $\Hilb ^{N}(V)$, called the {\it
Hilbert scheme of $G$-orbits} or {\it $G$-Hilbert scheme} and denoted
$\Ghilb{G}{V}$.  The universal family over $\Ghilb{G}{V}$ has a $G$
action, in which the coordinate ring of each fiber affords the regular
representation of $G$.

There is a canonical Chow morphism
\begin{equation}\label{e:other/G-Hilb-chow}
\Ghilb{G}{V}\rightarrow V/G
\end{equation}
making the $G$-Hilbert scheme $\Ghilb{G}{V}$ projective and birational
over $V/G$.  In the case $V = (\CC ^{2})^{n}$, $G = S_{n}$, as shown
in \cite{Hai99}, Section 4, the morphism in
\eqref{e:other/G-Hilb-chow} factors through a morphism
\begin{equation}\label{e:other/HilbG->Hn}
\phi \colon \Ghilb{S_{n}}{(\CC ^{2})^{n}}\rightarrow H_{n}
\end{equation}
of schemes over $S^{n}\CC ^{2}$.  In fact, if $Y$ denotes the
universal family over $\Ghilb{S_{n}}{(\CC ^{2})^{n}}$, then
$Y/S_{n-1}$ is a flat family of degree $n$, and it can be naturally
identified with a family of subschemes of $\CC ^{2}$.  The universal
property of $H_{n}$ then gives rise to $\phi $.

By Theorem~\ref{t:main}, the isospectral Hilbert scheme is a flat
family $X_{n}\subseteq H_{n}\times (\CC ^{2})^{n}$ of degree $n!$.
Hence the universal property gives rise to a morphism
$H_{n}\rightarrow \Hilb ^{n!}((\CC ^{2})^{n})$ whose image is clearly
contained in $\Ghilb{S_{n}}{(\CC ^{2})^{n}}$.  This morphism is
inverse to $\phi $, since it is so on the generic locus.  By
construction, $X_{n}$ is pulled back via this $\phi ^{-1}$ from the
universal family $Y$ over $\Ghilb{S_{n}}{(\CC ^{2})^{n}}$.  Hence we
have the following result, which is in fact equivalent to the flatness
of $\rho $ and thus to the Cohen-Macaulay property of $X_{n}$.


\begin{thm}\label{t:G-hilb}
The morphism $\phi $ in \eqref{e:other/HilbG->Hn} is an isomorphism of
$\Ghilb{S_{n}}{(\CC ^{2})^{n}}$ onto $H_{n}$, and it identifies the
universal family over $\Ghilb{S_{n}}{(\CC ^{2})^{n}}$ with the
isospectral Hilbert scheme $X_{n}$.
\end{thm}

In light of Theorem~\ref{t:G-hilb}, the vector bundle $P = \rho
_{*}\Ocal _{X_{n}}$ considered in the proof of Proposition
\ref{p:hilbert/n!->main} may be identified with the tautological
bundle on $\Ghilb{S_{n}}{(\CC ^{2})^{n}}$, a vector bundle whose
fibers afford the regular representation of $S_{n}$.  To each
irreducible representation $V^{\lambda }$ of $S_{n}$ is associated a
vector bundle $C_{\lambda } = \Hom _{S_{n}}(V^{\lambda },P)$ on
$\Ghilb{S_{n}}{(\CC ^{2})^{n}}$ known as a {\it character sheaf}.  We
have $P = \bigoplus _{\lambda } C_{\lambda }\otimes V^{\lambda }$, and
since $P$ affords the regular representation, $C_{\lambda }$ has rank
$\dim (V^{\lambda }) = \chi ^{\lambda }(1)$.  It follows from
Theorem~\ref{t:main}, Proposition \ref{p:hilbert/equiv}, and
Theorem~\ref{t:MPK} that the Kostka-Macdonald coefficient
$\tilde{K}_{\lambda \mu }(q,t)$ may be interpreted as the doubly
graded Hilbert series of the fiber of $C_{\lambda }$ at the $\TT
^{2}$-fixed point $I_{\mu }$.  Equivalently, $\tilde{K}_{\lambda \mu
}(q,t)$ is the character of the $\TT ^{2}$ action on this fiber.

We remark that the geometric picture described by
Theorems~\ref{t:main} and \ref{t:G-hilb} applies more generally with
any smooth complex quasiprojective surface $E$ replacing $\CC ^{2}$.
Specifically, we have that $\Ghilb{S_{n}}{E^{n}}$ is isomorphic to the
Hilbert scheme $\Hilb ^{n}(E)$, and the universal family over
$\Ghilb{S_{n}}{E^{n}}$ may be identified with the isospectral Hilbert
scheme $X_{n}(E)$, which is also the blowup of $E^{n}$ with center the
union of the pairwise diagonals.  To see this, note that Lemma
\ref{l:hilbert/product} holds in this situation and reduces the
question to a local one over points of the form $(P,P,\ldots,P)\in
E^{n}$, for $P\in E$.  Passing to the completion with respect to the
ideal of such a point, we see that the Cohen-Macaulay property of
$X_{n}(E)$ and its identification with the blowup are equivalent to
the corresponding facts for $X_{n}(\CC ^{2})$.

Returning to the general situation, assume that $G$ is a subgroup of
$\SL (V)$.  Then $V/G$ is Gorenstein and has rational singularities.  A
desingularization
\begin{equation}\label{e:other/H->V/G}
\sigma \colon H\rightarrow V/G
\end{equation}
is said to be {\it crepant} if $\omega _{H} = \Ocal _{H}$.  The
generalized McKay correspondence, conjectured by Reid
\cite{Rei97} and proved by Batyrev \cite{Bat99}, asserts that if $H$
is a crepant resolution, then the sum of the Betti numbers of $H$ is
equal to the number of conjugacy classes, or the number of irreducible
characters, of $G$.  For $V = (\CC ^{2})^{n}$, $G = S_{n}$, the
Hilbert scheme $H_{n}$ is a crepant resolution of $S^{n}\CC ^{2}$, by
Proposition \ref{p:hilbert/omega-Hn}.  The computation of the homology
of $H_{n}$ by Ellingsrud and Str{\o}mme \cite{ElSt87} verifies the
McKay correspondence in this case.

In dimension $2$, it develops that $\Ghilb{G}{V}$ is a crepant
resolution of $V/G$ for every $G$.  Nakamura conjectured that this
should hold in dimension $3$ as well, and proved it for $G$ abelian
\cite{Nak00}.  Recently, Bridgeland, King and Reid \cite{BrKiRe99}
established Nakamura's conjecture by proving that a certain fiber
dimension condition on the map $\Ghilb{G}{V}\rightarrow V/G$ is
sufficient to imply that $\Ghilb{G}{V}$ is a crepant resolution of
$V/G$.  In this situation, moreover, the McKay correspondence holds in
a strong form, expressed as an equivalence of derived categories.  By
Theorem~\ref{t:G-hilb}, we know that $\Ghilb{G}{V}$ is a crepant
resolution of $V/G$ in the case $G = S_{n}$, $V = (\CC ^{2})^{n}$.  As
things stand, this is the only known higher-dimensional family of
groups for which $\Ghilb{G}{V}$ is a crepant resolution, excluding
examples built up as products of lower-dimensional cases.  It is
natural to ask whether Theorem~\ref{t:G-hilb} might generalize to
yield other such families.


\begin{prob}\label{prob:other/complex-reflection}
Let $G$ be a complex reflection group with defining representation
$W$, and let $V = W\oplus W^{*}$.  For which $G$ is $\Ghilb{G}{V}$ a
crepant resolution of $V/G$?
\end{prob}

The related problem of whether $(W\oplus W^{*})/G$ admits a crepant
resolution at all has been considered by Kaledin \cite{Kal99} and
Verbitsky \cite{Ver99}, who prove that this can happen {\it only} for
complex reflection groups.

We have limited information about the solution to Problem
\ref{prob:other/complex-reflection}.  By our results here, the
symmetric groups, or Weyl groups of type $A$, are examples for which
$\Ghilb{G}{V}$ is a crepant resolution.  By the classical McKay
correspondence for subgroups of $\SL (2)$, the cyclic groups $\ZZ _{n}$,
regarded as one-dimensional complex reflection groups, are also
examples.  The hyperoctahedral groups, or Weyl groups of type $B$, are
counterexamples.  This can be seen by explicit computation for
$B_{2}$.  Since $B_{2}$ occurs as the stabilizer of a vector in the
reflection representation of $B_{n}$, it follows that $B_{n}$ is a
counterexample for all $n\geq 2$.

For the hyperoctahedral groups $G = B_{n}$, and more generally for $G$
the wreath product of the symmetric group $S_{n}$ with the cyclic
group $\ZZ _{m}$, crepant resolutions of $(W\oplus W^{*})/G$ do exist.
In these cases, the action of $G$ on $V = W\oplus W^{*}$ is the wreath
product action on $(\CC ^{2})^{n}$, where the cyclic group $\ZZ _{m}$
acts on $\CC ^{2} = \CC \oplus \CC ^{*}$ as the direct sum of a
one-dimensional representation and its dual.  As observed by Wang
\cite{Wan99}, whenever $G$ acts on $V=(\CC ^{2})^{n}$ as the wreath
product of $S_{n}$ with a finite subgroup $\Gamma \subseteq \SL (2)$,
there are two natural crepant resolutions.  One is a component of the
$\Gamma $-fixed locus in the Hilbert scheme $H_{n|\Gamma |}$.  The
other is the Hilbert scheme $\Hilb ^{n}(X_{\Gamma })$ of the minimal
resolution $X_{\Gamma } = \Ghilb{\Gamma }{\CC ^{2}}$ of $\CC
^{2}/\Gamma $.  In general (and for $G = B_{2}$ in particular) these
two resolutions can be different, and neither of them need coincide
with $\Ghilb{G}{V}$.


\subsection{Higher dimensions}
\label{ss:other/d>2}

For $d>2$ the Hilbert scheme $\Hilb ^{n}(\CC ^{d})$ has in general
multiple irreducible components, frequently of dimension exceeding
$dn$, and bad singularities.  Nevertheless there is a distinguished
{\it principal component} $H_{n}(\CC ^{d})$, the closure of the open
subset parametrizing reduced subschemes $S\subseteq \CC ^{d}$ with $n$
distinct points.  The isospectral Hilbert scheme $X_{n}(\CC ^{d})$
over $H_{n}(\CC ^{d})$ can be defined as in \ref{d:hilbert/Xn}, and
the analogs of Propositions \ref{p:hilbert/blowup-H} and
\ref{p:hilbert/blowup} hold also for $d>2$.

\begin{conj}\label{conj:other/CM?}
For all $d$ and $n$, the isospectral Hilbert scheme $X_{n}(\CC ^{d})$
over the principal component $H_{n}(\CC ^{d})$ of $\Hilb ^{n}(\CC
^{d})$ is normal and Cohen-Macaulay.
\end{conj}

Since $H_{n}(\CC ^{d}) = X_{n}(\CC ^{d})/S_{n}$ the conjecture implies
that the principal component $H_{n}(\CC ^{d})$ is itself normal and
Cohen-Macaulay.  Note that it does not imply the analog of the $n!$
conjecture in $d$ sets of variables, since in general $H_{n}(\CC
^{d})$ is singular and thus the projection $\rho $ need not be flat.
In other words, Theorem~\ref{t:G-hilb} may fail, and $H_{n}(\CC ^{d})$
may no longer coincide with the $G$-Hilbert scheme $\Ghilb{S_{n}}{(\CC
^{d})^{n}}$ for $d>2$.

An example exhibiting the failure of the higher-dimensional analogs of
the $n!$ conjecture and Theorem~\ref{t:G-hilb} occurs in dimension
$d=3$, with $n=4$.  The analog of the $n!$ conjecture would require
the determinant
\begin{equation}\label{e:other/false-n!}
\Delta (\xx ,\yy ,\zz ) = \det \begin{bmatrix} 
1&	x_{1}&	y_{1}&	z_{1}\\
1&	x_{2}&	y_{2}&	z_{2}\\
1&	x_{3}&	y_{3}&	z_{3}\\
1&	x_{4}&	y_{4}&	z_{4}\\
\end{bmatrix}
\end{equation}
to have $4!=24$ linearly independent partial derivatives, but in fact
it has only $20$.  More generally, Tesler \cite{Tes99} has shown that
the analogous determinant in dimension $d$ has $\binom{2d}{d}$
independent derivatives.  

Nevertheless, Conjecture \ref{conj:other/CM?} does hold for the case
$d=3$, $n=4$, as we have verified using the computer algebra system
Macaulay \cite{BaSt94}.  Both $X_{4}(\CC ^{3})$ and $H_{4}(\CC ^{3})$
turn out to be normal and Gorenstein.  (Strictly speaking, we have
verified this over a field of characteristic 31991, Macaulay's
default, and of course modulo the correctness of the program.)


Conjecture \ref{conj:other/CM?} is rather speculative, and the author
would not be greatly surprised if it turned out to be wrong.  A more
definite conjecture, whose failure he would indeed consider a
surprise, is the following higher-dimensional version of the polygraph
theorem.

\begin{conj}\label{conj:other/d-variate}
The coordinate ring of the polygraph $Z(n,l)$ over $E=\AA ^{d}$ is a
free $k[\zz]$-module, where $\xx ,\yy, \ldots,\zz =
x_{1},y_{1},\ldots,z_{1}, \; \ldots ,\; x_{n},y_{n},\ldots,z_{n}$ are
the coordinates on $E^{n}$.
\end{conj}

This conjecture implies the analogs of Proposition
\ref{p:hilbert/J-is-free} and Corollary
\ref{c:hilbert/free->what-J-is} in more than two sets of variables,
with essentially the same proof.  Hence it implies that $X_{n}(\CC
^{d})$ is the blowup of $(\CC ^{d})^{n}$ with center the reduced union
of the pairwise diagonals, that it is arithmetically normal in its
embedding as a blowup, and that it is flat over the coordinate space
$\CC ^{n} = \Spec \CC [\zz ]$ in any one set of the variables.


\subsection{Bases}
\label{ss:other/bases}

The proof of the $n!$ conjecture via Proposition \ref{p:hilbert/equiv}
and Theorem~\ref{t:main} does not yield an explicit basis of the
space $D_{\mu }$ or the ring $R_{\mu }$.  We still have an important
combinatorial open problem.

\begin{prob}\label{prob:other/K-lam-mu-basis}
Find a doubly homogeneous basis of $R_{\mu }$, compatible with (or
triangular with respect to) the decomposition into isotypic components
as an $S_{n}$-module, and indexed combinatorially so as to yield a
purely combinatorial interpretation of the Kostka-Macdonald
coefficients $\tilde{K}_{\lambda \mu }(q,t)$.
\end{prob}

A related doubly graded $S_{n}$-module is the the ring $R_{n} = \CC
[\xx ,\yy ]/I_{n}$, where $I_{n}$ is the ideal generated by all
$S_{n}$-invariant polynomials without constant term.  This ring has
the same graded character as the space of {\it diagonal harmonics},
and is the subject of a series of conjectures relating its character
to various combinatorial objects and to Macdonald polynomials
\cite{GaHa96, Hai94}.  In particular it is conjectured that $\dim _{\CC
}R_{n} = (n+1)^{n-1}$.

Each ring $R_{\mu }$ for $\mu $ a partition of
$n$ is a quotient ring of $R_{n}$, so bases of $R_{\mu }$ may be
realized as subsets of a basis of $R_{n}$.


\begin{prob}\label{prob:other/diag-basis}
Find a doubly homogeneous basis of $R_{n}$ compatible with the
$S_{n}$-module structure, and for each partition $\mu $ of $n$, a
distinguished subset which is a basis of $R_{\mu }$.
\end{prob}

Algebraically, the relationship between $R_{n}$ and $R(n,l)$ can be
described as follows.  Let $l = n$.  In the polygraph $Z(n,n)$ we have
a subarrangement 
\begin{equation}\label{e:other/perm-arrangement}
Z = \bigcup W_{f},\quad \text{over permutations $f\colon
[n]\rightarrow [n]$}.
\end{equation}
Let $J\subseteq R(n,n)$ be the ideal of $Z$.  For any
$S_{n}$-invariant polynomial $p(\xx ,\yy )$ we clearly have $p(\aa
,\bb )-p(\xx ,\yy )\in J$.  Hence if $p\in I_{n}$ is an
$S_{n}$-invariant polynomial without constant term, then we have
$p(\aa ,\bb )\in J+(\xx ,\yy )$.  This yields a ring homomorphism
$\psi \colon R_{n}\rightarrow R(n,n)/(J+(\xx ,\yy ))$, sending $\xx
,\yy $ in $R_{n}$ to $\aa ,\bb $ in $R(n,n)/(J+(\xx ,\yy ))$.  We
conjecture that $\psi $ is an isomorphism.  Postulating this for the
moment, we have a surjection of $R(n,n)/(\xx ,\yy )$ onto a copy of
$R_{n}$ in the variables $\aa ,\bb $.

From the proof of Theorem~\ref{t:common-basis}, we have an at least
somewhat explicit free $k[\yy ]$-module basis of $R(n,n)$, and
therefore a vector space basis of $R(n,n)/(\yy )$.  This basis is not
compatible with the $S_{n}$-action, nor does a subset of it span the
ideal $(\xx )$.  Nevertheless, some subset of it must be a basis of
$R(n,n)/(\xx ,\yy )$, and some further subset must be a basis of
$R(n,n)/(J+(\xx ,\yy ))$.


\begin{prob}\label{prob:other/R(n,l)/(x,y)-basis}
(1) Find a distinguished subset of the free $k[\yy ]$-module basis of
$R(n,l)$ constructed in the proof of Theorem~\ref{t:common-basis}
which is a vector space basis of $R(n,l)/(\xx ,\yy )$.

(2) For $l = n$, find a further distinguished subset of the above
basis, with $(n+1)^{n-1}$ elements, such that setting $\xx = \yy = 0$
and then replacing $a_{1},b_{1},\ldots,a_{n},b_{n}$ with
$x_{1},y_{1},\ldots,x_{n},y_{n}$, we get a basis of $R_{n}$.  To agree
with the conjectures on $R_{n}$ discussed in \cite{Hai94}, this basis
should be indexed by {\it parking functions} in such a way that the
$x$-degree of the basis element is the weight of the parking function.

(3) For each $\mu $, find a still further distinguished subset of
the above basis of $R_{n}$, with $n!$ elements, whose image in $R_{\mu
}$ is a basis.  For each $\lambda $, distinguish yet another subset
whose projection on the isotypic component of $R_{\mu }$ affording the
$S_{n}$ character $\chi ^{\lambda }$ is a basis of that isotypic
component, so that the enumerator of its elements by $x$- and
$y$-degrees is equal to $\chi ^{\lambda }(1)\tilde{K}_{\lambda \mu
}(q,t)$.
\end{prob}


\bibliographystyle{hamsplain}
\bibliography{references}

\end{document}